\documentclass[10pt,a4paper,reqno]{amsart}
\pdfoutput=1
\usepackage{amsmath,amsfonts,amssymb,amsthm,epsfig,amscd,psfrag,latexsym,comment}
\usepackage{caption}
\usepackage{graphicx}
\usepackage{array,subfigure}
\usepackage{multirow}
\usepackage[all]{xy}
\input xy
\xyoption{all}

\usepackage[colorlinks=true, pdfpagemode=none, pdfmenubar=false, pdfstartview=FitH, linkcolor=blue, citecolor=blue, urlcolor=blue, pdffitwindow=false]{hyperref}

\title{On adjunctions for Fourier-Mukai transforms}
\author{Rina Anno}
\email{anno@math.uchicago.edu}
\address{Department of Mathematics\\
University of Chicago\\
5734 S. University Avenue\\
Chicago, Illinois 60637\\
USA}
\author{Timothy Logvinenko} 
\email{T.Logvinenko@warwick.ac.uk} 
\address{Mathematics Institute\\ 
University of Warwick\\
Coventry, CV4 7AL\\
UK}

\DeclareMathOperator{\iden}{Id}
\DeclareMathOperator{\homm}{Hom}
\DeclareMathOperator{\shhomm}{{\it\mathcal{H}om\rm}}
\DeclareMathOperator{\eend}{End}

\DeclareMathOperator{\picr}{Pic}

\DeclareMathOperator{\cl}{Cl}

\DeclareMathOperator{\spec}{Spec\;}

\DeclareMathOperator{\tor}{Tor}
\DeclareMathOperator{\supp}{Supp}

\DeclareMathOperator{\ev}{ev}

\DeclareMathOperator{\cohcat}{Coh}
\DeclareMathOperator{\modd}{\bf Mod}
\DeclareMathOperator{\lder}{\bf L}
\DeclareMathOperator{\rder}{\bf R}

\DeclareMathOperator{\id}{Id}

\DeclareMathOperator{\evmap}{ev}

\begin{document}

\def\bv{\mathbf{v}}
\def\kgc_{K^*_G(\mathbb{C}^n)}
\def\kgchi_{K^*_\chi(\mathbb{C}^n)}
\def\kgcf_{K_G(\mathbb{C}^n)}
\def\kgchif_{K_\chi(\mathbb{C}^n)}
\def\gpic_{G\text{-}\picr}
\def\gcl_{G\text{-}\cl}
\def\trch_{{\chi_{0}}}
\def\regring{{R}}
\def\regrep{{V_{\text{reg}}}}
\def\givrep{{V_{\text{giv}}}}
\def\lbar{{(\mathbb{Z}^n)^\vee}}
\def\genpx_{{p_X}}
\def\genpy_{{p_Y}}
\def\genpcn_{p_{\mathbb{C}^n}}
\def\gnat{gnat}
\def\twalg{{\regring \rtimes G}}
\def\L{{\mathcal{L}}}
\def\O{{\mathcal{O}}}
\def\gcd{\mbox{gcd}}
\def\lcm{\mbox{lcm}}
\def\tf{{\tilde{f}}}
\def\tD{{\tilde{D}}}

\def\mckquiv{\mbox{Q}(G)}
\def\C{{\mathbb{C}}}
\def\sF{{\mathcal{F}}}
\def\sW{{\mathcal{W}}}
\def\sL{{\mathcal{L}}}
\def\O{{\mathcal{O}}}
\def\Z{{\mathbb{Z}}}
\def\hmone{{\mathcal{W}}}

\theoremstyle{definition}
\newtheorem{defn}{Definition}[section]
\newtheorem*{defn*}{Definition}
\newtheorem{exmpl}[defn]{Example}
\newtheorem*{exmpl*}{Example}
\newtheorem{exrc}[defn]{Exercise}
\newtheorem*{exrc*}{Exercise}
\newtheorem*{chk*}{Check}
\newtheorem*{remarks*}{Remarks}
\theoremstyle{plain}
\newtheorem{theorem}{Theorem}[section]
\newtheorem*{theorem*}{Theorem}
\newtheorem{conj}[defn]{Conjecture}
\newtheorem*{conj*}{Conjecture}
\newtheorem{prps}[defn]{Proposition}
\newtheorem*{prps*}{Proposition}
\newtheorem{cor}[defn]{Corollary}
\newtheorem*{cor*}{Corollary}
\newtheorem{lemma}[defn]{Lemma}
\newtheorem*{claim*}{Claim}
\newtheorem{Specialthm}{Theorem}
\renewcommand\theSpecialthm{\Alph{Specialthm}}
\numberwithin{equation}{section}
\renewcommand{\textfraction}{0.001}
\renewcommand{\topfraction}{0.999}
\renewcommand{\bottomfraction}{0.999}
\renewcommand{\floatpagefraction}{0.9}
\setlength{\textfloatsep}{5pt}
\setlength{\floatsep}{0pt}
\setlength{\abovecaptionskip}{2pt}
\setlength{\belowcaptionskip}{2pt}
\begin{abstract}
We show that the adjunction counits of a Fourier-Mukai transform 
$\Phi\colon D(X_1) \rightarrow D(X_2)$ arise from maps 
of the kernels of the corresponding Fourier-Mukai transforms. 
In a very general setting of proper separable schemes 
of finite type over a field we write down these maps of kernels
explicitly~-- facilitating the computation of the twist (the cone of
an adjunction counit) of $\Phi$.  We also give another description of
these maps, better suited to computing cones if the kernel of 
$\Phi$ is a pushforward from a closed subscheme $Z \subset X_1 \times X_2$.
Moreover, we show that we can replace the condition of properness 
of the ambient spaces $X_1$ and $X_2$ by that of $Z$ being proper over 
them and still have this description apply as is. This can be used, 
for instance, to compute spherical twists on non-proper varieties 
directly and in full generality. 
 \end{abstract}

\maketitle

\section{Introduction} \label{section-intro}

The bounded derived category $D(X)$ of coherent sheaves on a variety
$X$ had long been recognized as a crucial invariant of $X$ which holds 
a wealth of information about its geometry. In order to work
conveniently with functors between the derived categories of two varieties
the language of Fourier-Mukai transforms was developed by Mukai, 
Bondal and Orlov, Bridgeland and many others. In brief, we can 
define a functor $D(X_1) \rightarrow D(X_2)$ by specifying 
an object in the derived category of $D(X_1 \times X_2)$. A morphism 
between such defining objects induces a natural transformation 
between the functors. In this paper we write down the adjunction counit 
of a general Fourier-Mukai transform in this language --- as morphisms 
of defining objects.  

 Let $X_1$ and $X_2$ be a pair of smooth projective varieties. 
We have the following commutative diagram:
\begin{align} \label{eqn-big-projection-tree-1}
\xymatrix{
& & X_1 \times X_2 \times X_1 \ar[ld]_{\pi_{12}} \ar[d]^{\pi_{13}} \ar[rd]^{\pi_{23}} & & \\
& X_1 \times X_2 \ar[ld]_{\pi_1} \ar[rd]^>>>>>>>{\pi_2} & X_1 \times X_1
\ar[lld]_>>>>>>>>>>>>>>{\tilde{\pi}_1}
\ar[rrd]^>>>>>>>>>>>>>>{\tilde{\pi}_2} & X_2 \times X_1 \ar[ld]_>>>>>>>{\pi_2} \ar[rd]^{\pi_1} & \\
X_1 & & X_2 & & X_1
}
\end{align}

Let $E \in D(X_1 \times X_2)$. The \em Fourier-Mukai 
transform from $X_1$ to $X_2$ with kernel $E$ \rm 
is the functor 
\begin{align}
\Phi_E(-) = \pi_{2*}\left(E \otimes \pi^*_1(-)\right).
\end{align}
Here and throughout the paper all the functors are derived unless 
mentioned otherwise. It is well-known (e.g. \cite{BonOr95}, Lemma 1.2) 
that the left adjoint of $\Phi_E$ is the Fourier-Mukai transform 
from $D(X_2)$ to $D(X_1)$ with kernel 
$E^\vee \otimes \pi_1^!(\mathcal{O}_{X_1})$ 
where $\pi_1^!(\mathcal{O}_{X_1}) = \pi_2^*(\omega_{X_2})[\dim X_2]$. 
Denote this adjoint by $\Phi^{\text{ladj}}_E$. 
A composition of Fourier-Mukai transforms 
is again a Fourier-Mukai transform (\cite{Muk81}, Prop. 1.3). 
In particular, $\Phi^{\text{ladj}}_E \Phi_E$ is 
the Fourier-Mukai transform $D (X_1) \rightarrow D(X_1)$
with kernel
\begin{align}
\label{eqn-intro-kernel-of-Phi^ladj_E-Phi_E}
Q = \pi_{13 *}\left(\pi_{12}^* E \otimes \pi_{23}^* E^\vee
\otimes \pi_{23}^* \pi^!_1(\mathcal{O}_{X_1})\right). 
\end{align}
On the other hand, the identity functor $\id$ is 
the Fourier-Mukai transform $D(X_1) \rightarrow D(X_1)$ with kernel 
$\mathcal{O}_{\Delta} = \Delta_* \mathcal{O}_{X_1}$ where $\Delta$
is the diagonal inclusion $X_1 \hookrightarrow X_1 \times X_1$.

Consider now the left adjunction counit 
\begin{align} \label{eqn-canonical-fmt-adjunction-morphism}
\Phi^{\text{ladj}}_E \Phi_E \rightarrow \id. 
\end{align}
In general, morphisms between Fourier-Mukai kernels 
map neither injectively nor surjectively to natural 
transformations between the Fourier-Mukai transforms. 
Thus there is no \em a priori \rm reason for 
\eqref{eqn-canonical-fmt-adjunction-morphism} to come 
from some morphism $Q \rightarrow \mathcal{O}_{\Delta}$.
In this paper we construct explicitly a natural choice 
of such morphism, working in a much greater generality 
of separated schemes of finite type over a field.  

The principal application is to compute, and even define,
\em spherical twists\rm. These are an important class of 
auto-equivalences of the derived category $D(X)$ of a variety $X$. 
They are first examples of genuinely derived
auto-equivalences, in a sense that they are neither shifts, 
nor come from auto-equivalences of the underlying
abelian category $\cohcat X$. In brief, \em a spherical 
twist \rm is an auto-equivalence of $D(X)$ produced from 
a \em spherical object \rm in $D(X)$ or, more generally,  
\em a spherical functor \rm $D(Y) \rightarrow D(X)$.
Spherical objects were introduced by Seidel and Thomas in 
\cite{SeidelThomas-BraidGroupActionsOnDerivedCategoriesOfCoherentSheaves}
as mirror symmetry analogues of Lagrangian spheres on a
symplectic manifold. Their defining properties ensure 
that the twist by a spherical object is an auto-equivalence 
of $D(X)$. This was generalised in \cite{Anno-SphericalFunctors} 
to exact functors between triangulated categories 
in such a way that Seidel-Thomas spherical objects are 
precisely the (Fourier-Mukai kernels of) spherical functors 
$D(\spec k) \rightarrow D(X)$, where $k$ is the base field.

Taking the twist of a functor is completely general and
does not in itself rely on the fact that the functor is spherical. 
The ideal definition would be the following:

\vskip0.2cm
\bf ``Definition'': \rm
Let $C_1$ and $C_2$ be triangulated categories. Let $S$
be an exact functor $C_1 \rightarrow C_2$ which has a right (resp. 
left) adjoint $R$ (resp. $L$). The \em twist \rm (resp. 
the \em dual co-twist\rm) of $S$ is the functor $T_S\colon C_2
\rightarrow C_2$ (resp. $F'_S\colon C_1 \rightarrow C_1$)
which is the functorial cone of the adjunction counit 
$SR \rightarrow \id$ (resp. $LS \rightarrow \id$).
\vskip0.2cm

The problem with this definition is the well-known fact that 
cones in triangulated categories are not functorial. 
The cone of a morphism between two objects is uniquely defined 
(up to an isomorphism), but a cone of a morphism between two functors 
might not exist or might not be unique. This is usually fixed 
by restricting to a setting where the cone of 
a morphism of functors is well-defined, cf.~\cite{Anno-SphericalFunctors}, \S1. 
One way is to consider only the functors which are Fourier-Mukai 
transforms and only the natural transformations which come from 
morphisms of Fourier-Mukai kernels. But then to define 
a twist of a Fourier-Mukai transform we need a natural 
choice of the morphism of Fourier-Mukai kernels underlying 
the corresponding adjunction counit, while to compute 
the twist we need an efficient way of computing the cone 
of this morphism. This paper addresses both of these issues. 

The construction of the natural morphism of Fourier-Mukai 
kernels underlying the adjunction counit of a general 
Fourier-Mukai transform is carried out 
in Section \ref{section-fmt-adjunction-morphism}.
Thanks to the recent advances in Grothendieck duality 
machinery summarised in Section \ref{section-preliminaries}
we can work with separated schemes of finite type over a field 
and with derived categories $D_{qc}(-)$ of unbounded 
complexes with quasi-coherent cohomology. 
So let $X_1$ and $X_2$ be two separated schemes of finite type,
$E$ a perfect object of $D(X_1 \times X_2)$ and 
$\Phi_E$ the Fourier-Mukai transform $D(X_1) \rightarrow D(X_2)$ 
with kernel $E$. Let $X_2$ be proper, so that 
the left adjoint $\Phi^{\text{ladj}}_E$ of $\Phi_E$ is 
again a Fourier-Mukai transform. Then 
the left adjunction counit $\Phi^{\text{ladj}}_E
\Phi_E \rightarrow \id$ is induced by the morphism 
$Q = \pi_{13 *}\left(\pi_{12}^* E \otimes \pi_{23}^* E^\vee
\otimes \pi_{23}^* \pi^!_1(\mathcal{O}_{X_1})\right) \rightarrow 
\mathcal{O}_\Delta$ 
which roughly is the composition of the following:

\begin{align}
\label{eqn-intro-derived-restriction-to-diagonal-in_X1X2X1}
\pi_{13 *}\left( 
\text{The adjunction unit $\id \rightarrow \Delta_{13*} \Delta_{13}^*$
for the diagonal $X_1 \times X_2 \xrightarrow{\Delta_{13}}  
X_1 \times X_2 \times X_1$} \right)
\\   
\label{eqn-intro-evaluation-map-on-X1X2}
\Delta_* \pi_{1 *} \left(\text{The evaluation map } E \otimes E^\vee \rightarrow
\mathcal{O}_{X_1 \times X_2} \text{ on } X_1 \times X_2\right) \\
\label{eqn-intro-adjunction-morphisms-for-pi_1}
\Delta_* \left(\text{The adjunction counit } \pi_{1 *}
\pi^!_1(\mathcal{O}_{X_1}) \rightarrow \mathcal{O}_{X_1} \right)
\end{align}

For the precise formulas see Theorem
\ref{theorem-left-adjunction-counit-morphism}. 
When $X_1$ is also proper $\Phi_E$, $\Phi^{ladj}_E$ and 
\eqref{eqn-intro-derived-restriction-to-diagonal-in_X1X2X1}-\eqref{eqn-intro-adjunction-morphisms-for-pi_1}
restrict to the full subcategories of $D_{qc}(-)$
consisting of bounded complexes with coherent cohomologies. 
If $X_2$ is smooth  
$\pi^!_1(\mathcal{O}_{X_1}) =  \pi_2^*(\omega_{X_2})[\dim X_2]$ as before. 
Theorem \ref{theorem-right-adjunction-counit-morphism} give 
the analogous result for the right adjunction counit. 

This allows us to define the twist and the dual co-twist of 
any Fourier-Mukai transform. Section \ref{section-the-pushforward-case}
deals with the issue of computing them. Anyone trying to compute
the cone of the decomposition 
\eqref{eqn-intro-derived-restriction-to-diagonal-in_X1X2X1}-\eqref{eqn-intro-adjunction-morphisms-for-pi_1}
will find it ill-suited to the task if the support of $E$ has high
codimension in $X_1 \times X_2$. We give an example 
in Section \ref{section-the-global-intersection-example}
with $E$ the structure sheaf $\mathcal{O}_Z$ 
of a complete intersection subscheme $Z$ in 
$X_1 \times X_2$ of codimension $d > 0$ which satisfies 
certain transversality conditions. Then morphisms  
\eqref{eqn-intro-derived-restriction-to-diagonal-in_X1X2X1}
and \eqref{eqn-intro-evaluation-map-on-X1X2} both have huge cones
with non-zero cohomologies in all degrees from $-d$ to $0$. However
these two cones mostly annihilate each other and the cone 
of 
composition \eqref{eqn-intro-derived-restriction-to-diagonal-in_X1X2X1}-\eqref{eqn-intro-evaluation-map-on-X1X2} is actually 
quite small. This suggests an alternative decomposition of
\eqref{eqn-intro-derived-restriction-to-diagonal-in_X1X2X1}-\eqref{eqn-intro-evaluation-map-on-X1X2}
better suited to computing cones, 
cf.~\eqref{eqn-gci-alternative-decomposition}. 

In the rest of Section \ref{section-the-pushforward-case} 
we make this into a general argument. The key idea
is to take the decomposition 
\eqref{eqn-intro-derived-restriction-to-diagonal-in_X1X2X1}-\eqref{eqn-intro-adjunction-morphisms-for-pi_1}
obtained in Section \ref{section-fmt-adjunction-morphism} and
apply to it the base change for K{\"u}nneth maps. If $E$ 
is a pushforward of an object from a closed 
subscheme $Z \overset{\iota_Z}{\hookrightarrow} X_1 \times X_2$, 
then the evaluation map 
$E \otimes E^\vee \rightarrow \mathcal{O}_{X_1 \times X_2}$
involves the derived self-intersection of $Z$ inside 
$X_1 \times X_2$. In precise terms, it involves
the K{\"u}nneth map (see Section  
\ref{section-kunneth-maps-and-base-change} for the definition)
for the fiber square $\sigma_\Delta$ depicted 
on the left in \eqref{eqn-intro-two-fiber-squares}: 

\begin{align} 
\label{eqn-intro-two-fiber-squares}
\sigma_\Delta\colon
\vcenter{\xymatrix{
Z \ar[d] \ar[r] & Z \ar[d]^{\iota_Z} \\
Z \ar[r]_<<<<<{\iota_Z} & X_1 \times X_2
}} 
\quad  
\xleftarrow{\text{Restriction to
$X_1 \times X_2 \xrightarrow{\Delta_{13}}  
X_1 \times X_2 \times X_1$}}
\quad 
\sigma\colon 
\vcenter{\xymatrix{
Z' \ar[d] \ar[r] &  X_1 \times Z \ar[d]^{\iota_{Z23}} \\
Z \times X_1 \ar[r]_<<<<<{\iota_{Z12}} & X_1 \times X_2 \times X_1
}}
\end{align}

Thus in \eqref{eqn-intro-derived-restriction-to-diagonal-in_X1X2X1}-\eqref{eqn-intro-evaluation-map-on-X1X2} 
we first restrict fiber square $\sigma$ to the diagonal 
$X_1 \times X_2$ in $X_1 \times X_2 \times X_1$ 
which turns it into $\sigma_\Delta$ and then we do the K{\"u}nneth 
map on $\sigma_\Delta$. Given two subschemes, 
the cone of the K{\"u}nneth map for the fiber square of their intersection 
reflects, roughly, how far this intersection is from transverse. 
In $\sigma_\Delta$ we have the self-intersection 
of $Z$ in $X_1 \times X_2$ which is the opposite of transverse. 
This suggests first doing the K{\"u}nneth map on $\sigma$, 
as the intersection of $Z \times X_1$ with $X_1 \times Z$ 
in $X_1 \times X_2 \times X_1$ may be more transverse, and 
then restricting to the diagonal $Z$ in $Z'$. 

Write $\pi_{Z1}$ for the composition 
$Z \overset{\iota_Z}{\hookrightarrow} X_1 \times X_2 \xrightarrow{\pi_1} X_1$.
In Prop. \ref{prps-base-change-for-kunneth-maps} we 
prove that K{\"u}nneth maps commute with arbitrary base change. Then in  
Theorem \ref{theorem-pushfwd-left-adjunction-counit}
we show that the composition
\eqref{eqn-intro-derived-restriction-to-diagonal-in_X1X2X1}-\eqref{eqn-intro-adjunction-morphisms-for-pi_1}
is isomorphic to roughly the following 
(cf.~Theorem \ref{theorem-pushfwd-left-adjunction-counit} for precise 
formulas):
\begin{align}
\label{eqn-intro-kunneth-map-upstairs}
\pi_{13 *}  \left( \text{The K{\"u}nneth map for } \sigma \right) \\
\label{eqn-intro-derived-restriction-to-diagonal-in-D'}
\pi_{13 *} \iota_{Z'*} 
\left(\text{The adjunction unit $\id \rightarrow \Delta'_{*} \Delta'^{*}$
for the diagonal $Z \xrightarrow{\Delta'} Z'$} \right) 
\\ 
\label{eqn-intro-evaluation-map-on-D-v2}
\Delta_* \pi_{Z1 *}  
\left(\text{The evaluation map for } E \text{ on } Z \right) \\
\label{eqn-intro-adjunction-counit-for-pi_D1-v2}
\Delta_* \left(\text{The adjunction counit } \pi_{Z1 *}
\pi^!_{Z1}(\mathcal{O}_{X_1}) \rightarrow \mathcal{O}_{X_1} \right)
\end{align}
This is our preferred decomposition of morphism 
$Q \rightarrow \mathcal{O}_\Delta$. Theorem 
\ref{theorem-pushfwd-right-adjunction-counit} gives 
the analogous statement for the right adjunction counit. 

One advantage of decomposition 
\eqref{eqn-intro-kunneth-map-upstairs}-\eqref{eqn-intro-adjunction-counit-for-pi_D1-v2} is that most of the morphisms in it can become 
isomorphisms under fairly
reasonable assumptions on $E$ and $Z$. Indeed, while the K{\"u}nneth map 
for square $\sigma_\Delta$ is never an isomorphism unless $Z$ 
is the whole of $X_1 \times X_2$, the K{\"u}nneth map for $\sigma$ 
is an isomorphism whenever 
the intersection of $Z \times X_1$ with $X_1 \times Z$ in 
$X_1 \times X_2 \times X_1$ is transverse. 
The evaluation map for $E$ on $Z$ is an isomorphism 
whenever $E$ is a line bundle 
or any invertible object of $D(Z)$. 
The adjunction counit in \eqref{eqn-intro-adjunction-counit-for-pi_D1-v2}
is an isomorphism whenever 
$Z \xrightarrow{\pi_{Z1}} X_1$ is such that $\pi_{Z1 *}
\mathcal{O}_Z = \mathcal{O}_{X_1}$, e.g. $Z$ is a blowup of
$X_1$ or a Fano fibration over it. This allows for a number
of scenarios where the twist or the dual co-twist of $\Phi_E$
can be written down fairly easily, as we demonstrate 
in Cor. \ref{cor-best-case-scenario}.

Another advantage of decomposition
\eqref{eqn-intro-kunneth-map-upstairs}-\eqref{eqn-intro-adjunction-counit-for-pi_D1-v2}
is that it moves the action away from ambient spaces
$X_1 \times X_2 \times X_1$ and $X_1 \times X_2$ to their subschemes
$Z'$ and $Z$. This allows us to replace
the assumption of $X_2$ being proper by the assumption 
of $Z$ being proper over $X_1$ and $X_2$ (see Theorem 
\ref{theorem-pushfwd-left-adjunction-counit}). 
Something to be appreciated by those who
want to do spherical twists on non-compact varieties, e.g.
total spaces of cotangent bundles of projective varieties.

Finally, in Section \ref{section-an-example} we give an example 
of an explicit computation using Theorem 
\ref{theorem-pushfwd-left-adjunction-counit}. 
We consider the naive derived category transform induced 
by a Mukai flop. This transform is not an equivalence - it was proved by 
Namikawa in \cite{Namikawa-MukaiFlopsAndDerivedCategories} 
by direct comparison of $\homm$ spaces. We demonstrate how 
its dual co-twist can be computed quickly and efficiently by our methods.

\bf Acknowledgements: \rm We would like to thank Alexei Bondal and Paul
Bressler for useful discussions. We would like to thank an anonymous referee 
for major improvements to all aspects of the paper. 
The first author would like to thank the Department of 
Mathematics of the University of Chicago for their support. 
The second author would like to thank the University of
Liverpool, the Max-Planck-Institut f{\"u}r Mathematik and the Steklov
Mathematical Institute for their hospitality during his work on this paper. 

\section{Preliminaries} \label{section-preliminaries}

Let $k$ be an algebraically closed field of characteristic $0$. 
The level of generality we choose to work at in the main 
body of this paper is that of separated schemes of finite type
over $k$. These assumptions are necessary for the Grothendieck
duality machinery which ensures that the direct image functor
in the definition of a Fourier-Mukai transform has a right adjoint. 
Without them we cannot expect a general Fourier-Mukai transform 
to have a right and a left adjoint. 

Some of the auxiliary results we prove along the way 
hold in a greater generality than the one above. We would 
like to think of these results as being of potential interest 
to others who find themselves in an unfortunate situation of 
having to show a complicated diagram of derived functors to commute. 
We try therefore to state these results in maximal generality they 
hold at. 

By a ringed space we always mean a commutative ringed space. 
By a \em concentrated \rm map of schemes we mean a map which is quasi-compact 
and quasi-separated. A scheme $X$ is said to be
\em concentrated \rm if it is concentrated over $\spec \mathbb{Z}$. If $Y$
is a concentrated scheme, then a map $X \rightarrow Y$ is 
concentrated
if and only if $X$ is concentrated \cite[\S1.2]{Grothendieck-EGA-IV-1}.

We make frequent use of a notion of a \em perfect \rm map of schemes $X
\xrightarrow{f} Y$, cf. \cite[\S4]{IllusieConditionsDeFinitudeRelative}. 
For maps of finite type between noetherian schemes $f$ is perfect
if and only if it is of finite $\tor$-dimension, i.e. the derived 
functor of $f^*$ is cohomologically bounded. 

Given an adjoint pair of functors $(F,G)$, by the \em right adjoint
with respect to $F$ \rm of some natural transformation 
$FH_1 \rightarrow H_2$, we mean the natural transformation 
$H_1 \rightarrow GH_2$ induced by the adjunction. 
Similarly, by \em the left adjoint with respect to $G$ \rm
of some $H_1 \rightarrow GH_2$ we mean the $FH_1 \rightarrow H_2$
induced by the adjunction.  

Throughout the paper we employ a variety of greek letters 
to denote an assortment of natural maps which exist 
between compositions of standard derived functors. 
These are defined at length over the course of Sections 
\ref{section-derived-functors-and-derived-categories}-\ref{section-standard-relations-between-derived-functors},
but for the convenience of our readers we have also compiled 
a brief index: 

\noindent
\begin{footnotesize}
\begin{minipage}{0.5\linewidth}
\begin{align*}
\begin{array}{|c|c|c|}
\hline
\alpha_f
&
\overset{\text{the projection formula}}{f_* A \otimes B \rightarrow f_* (A \otimes f^*B)} 
&
\eqref{eqn-definition-of-projection-formula-morphism}
\\
\hline
\beta_f
&
\id \rightarrow f_* f^* 
&
\eqref{eqn-adjunction-unit-and-counit-for-f^*-f_*}
\\
\hline
\gamma_f 
&
f^* f_* \rightarrow \id
&
\eqref{eqn-adjunction-unit-and-counit-for-f^*-f_*}
\\
\hline
\delta_f
&
\overset{\text{the sheafified Grothendieck duality}}{
\text{\footnotesize$f_* \rder\shhomm_X(A, f^\times B) \rightarrow
\rder\shhomm_Y(f_*A,B)$} 
}
&
\eqref{eqn-sheafified-Grothendieck-duality-morphism}
\\
\hline
\epsilon_f 
&
f_* f^\times \rightarrow \id 
&
\eqref{eqn-adjunction-unit-and-counit-for-f_*-f^x}
\\
\hline
\zeta_{g,f} 
&
f^* g^* \xrightarrow{\sim} (g \circ f)^* 
&
\eqref{eqn-pseudofunctoriality-inverse-image}
\\
\hline
\eta_{g,f}
&
(g \circ f)_* \xrightarrow{\sim} g_*f_*
&
\eqref{eqn-pseudofunctoriality-direct-image}
\\
\hline
\theta_{A,B}
& 
A \longrightarrow 
\rder\shhomm_X\left(\rder\shhomm_X\left(A,B\right),B \right)
&
\eqref{eqn-reflexivity-morphism}
\\
\hline
\theta_E 
& 
E \rightarrow E^{\vee\vee}
&
\eqref{eqn-dualizing-morphism} 
\\
\hline
\kappa_f 
&
f_*A \otimes f_*B \rightarrow f_*(A \otimes B)
&
\eqref{eqn-symmetric-monoidal-functor-struct-for-direct-image}
\\
\hline 
\end{array}
\end{align*}
\end{minipage}
\begin{minipage}{0.5\linewidth}
\begin{align*}
\begin{array}{|c|c|c|}
\hline 
\kappa_\sigma
&
\overset{\text{the K\"unneth map}}{\text{\footnotesize $f_{1 *}
(A_1) \otimes f_{2 *} (A_2) \rightarrow h_*\left(g_1^*(A_1) \otimes
g_2^*(A_2)\right)$}}
&
\eqref{eqn-kunneth-map-definition}
\\
\hline 
\lambda_f 
&
\id \rightarrow f^\times f_*
&
\eqref{eqn-adjunction-unit-and-counit-for-f_*-f^x}
\\
\hline 
\mu_{\sigma}
& 
\overset{\text{the base change}}{g^* f_* \rightarrow f'_* g'^*}
&
\eqref{eqn-definition-of-base-change-morphism}
\\
\hline 
\nu_f 
&
f^*(A \otimes B) \xrightarrow{\sim} f^*(A) \otimes f^* B 
&
\eqref{eqn-definition-of-morphism-nu}
\\
\hline 
\xi
&
\text{\footnotesize $\rder\shhomm_X\left(A,B\right) \otimes C
\rightarrow 
\rder\shhomm_X\left(A,B \otimes C \right)$}
&
\eqref{eqn-bringing-an-object-into-rhom-bracket}
\\
\hline 
\xi_E
&
E^\vee \otimes (-) \xrightarrow{\sim} \rder \shhomm_X(E, -) 
&
\eqref{eqn-perfect-project-rhom-is-tensor-with-dual}
\\
\hline 
\rho
&
(A \otimes B) \otimes C 
\xrightarrow{\sim} 
A \otimes (B \otimes C)
&
\eqref{eqn-associativity-of-tensor-product-definition}
\\
\hline 
\tau_f 
&
\text{\footnotesize $f_* \rder\shhomm_X(f^*A,B) \xrightarrow{\sim} 
\rder\shhomm_Y(A, f_*B)$}
&
\eqref{eqn-definition-of-sheafified-f^*-f_*-adjunction}
\\
\hline 
\upsilon_A
&
\text{\tiny $\rder\shhomm(A\otimes B,C)
\xrightarrow{\sim}
\rder\shhomm(B,\rder\shhomm(A,C))$}
&
\eqref{eqn-sheafified-tensor-product-adjunction}
\\
\hline 
\chi_f
&
f^\times A \otimes f^* B \rightarrow f^\times(A \otimes B) 
&
\eqref{eqn-definition-of-f^!-f!xf^*-morphism}
\\
\hline
\end{array}
\end{align*}
\end{minipage}
\end{footnotesize}

\subsection{Derived categories and derived functors} 
\label{section-derived-functors-and-derived-categories}

Let $X$ be a scheme or a ringed space. We 
denote by $D(\mathcal{O}_X\text{-}\modd)$ 
the unbounded derived category of the abelian category $\mathcal{O}_X$-$\modd$. 
We denote by $D_{\text{qc}}(X)$ (resp. $D(X)$) the full 
subcategory of $\mathcal{O}_X$-$\modd$ consisting of 
complexes with quasi-coherent (resp. bounded and coherent) cohomology.  
We denote by $D_\text{perf}(X)$ the full subcategory of $D(X)$
consisting of the objects which are locally quasi-isomorphic 
to a bounded complex of free $\mathcal{O}_X$-modules of finite rank.

For a reference text on derived categories and derived functors
we recommend \cite{Hartshorne-Residues-and-Duality}, 
for the traditional approach, and 
\cite{Lipman-NotesOnDerivedFunctorsAndGrothendieckDuality}, 
for a more modern approach. One should also mention the expositions
in \cite{KashiwaraSchapira-CategoriesAndSheaves} and 
\cite{AmnonNeeman-TriangulatedCategories}. A key feature of the
modern approach is that thanks to 
the results of \cite{Spaltenstein-ResolutionsOfUnboundedComplexes}
we can now work freely with unbounded complexes.
The authors of this paper adhere to a general principle that 
wherever possible general results on derived functors and 
isomorphisms between them should first be proved in the 
setting of $D_\text{qc}(-)$, and then shown to restrict 
to the usual setting of $D(-)$ where applicable.

All the functors in this paper are assumed to be derived, 
unless specifically mentioned otherwise. With two exceptions listed
below we suppress all the usual $\rder$'s and $\lder$'s and use 
the same notation for the derived functor as for its abelian 
category counterpart. Below we summarize basic facts about the
derived functors we make use of.

Let $X$ be a ringed space. The derived tensor product functor exists as a functor 
$$ (-) \otimes (-)\colon 
D(\mathcal{O}_X\text{-}\modd) \times D(\mathcal{O}_X\text{-}\modd)
\rightarrow
D(\mathcal{O}_X\text{-}\modd). 
$$
and always restricts to a functor $D_\text{qc}(X) \times
D_\text{qc}(X) \rightarrow D_\text{qc}(X)$
\cite[\S2.5]{Lipman-NotesOnDerivedFunctorsAndGrothendieckDuality}. 
For $X$ a locally noetherian scheme and for $A \in D_{\text{perf}}(X)$ 
the functor $A \otimes -$ restricts to a functor $D(X) \rightarrow D(X)$
\cite[Prop. II.4.3]{Hartshorne-Residues-and-Duality}.
Similarly, for any $n \in \mathbb{Z}$ the derived tensor
product functor in $n$ variables
$ (-) \otimes \dots \otimes (-) $ exists as a functor 
from the product of $n$ copies of $D(\mathcal{O}_X\text{-}\modd)$ into
$D(\mathcal{O}_X\text{-}\modd)$
\cite[\S2.5.9]{Lipman-NotesOnDerivedFunctorsAndGrothendieckDuality}.

The derived functor of the functor $\homm_{X}(-,-)$ of taking 
the global $\homm$ space between two $\mathcal{O}_X$-modules
exists as a functor
$$ \rder\homm_{X}(-,-)\colon 
D(\mathcal{O}_X\text{-}\modd)^\text{opp} \times D(\mathcal{O}_X\text{-}\modd)
\rightarrow D(\Gamma(\mathcal{O}_X)\text{-}\modd), $$
see \cite[\S2.4]{Lipman-NotesOnDerivedFunctorsAndGrothendieckDuality}.
We make an exception and do not supress `$\rder$' here in order to  
differentiate the object $\rder\homm_X(A,B)$ in 
$D(\Gamma(\mathcal{O}_X)\text{-}\modd)$ from the morphism space 
$\homm_{D(X)}(A,B)$.
Similarly, the derived functor of the sheafified $\homm$ functor
$\shhomm_X(-,-)$ exists as a functor   
$$ \rder\shhomm_{X}(-,-)\colon 
D(\mathcal{O}_X\text{-}\modd)^\text{opp} \times D(\mathcal{O}_X\text{-}\modd)
\rightarrow D(\mathcal{O}_X\text{-}\modd)$$
We do not suppress `$\rder$' here to emphasize the
relation with $\rder\homm_X$. If $X$ is a locally noetherian scheme, 
then for any $A \in D(X)$ the functor $\rder\shhomm_X(A, -)$ 
restricts to a functor $D^+_{qc}(X) \rightarrow D_{qc}(X)$ 
\cite[Prop. II.3.3]{Hartshorne-Residues-and-Duality}. Here
$D^+_{qc}(X)$ is the subcategory of $D_{qc}(X)$ consisting of
complexes with bounded below cohomology. If $X$ is a noetherian scheme
and $A$ is perfect the functor $\rder\shhomm_X(A, -)$ restricts
to a functor $D_\text{qc}(X) \rightarrow D_\text{qc}(X)$ 
and then to a functor $D(X) \rightarrow D(X)$ 
\cite[Lemma 1.4.6]{AvramovIyengarLipman-ReflexivityAndRigidityForComplexesIISchemes}.

Let now $Y$ be another ringed space, and let $f\colon X \rightarrow Y$
be a map of ringed spaces.

The derived direct image functor exists as a functor 
$$ f_*(-)\colon 
D(\mathcal{O}_X\text{-}\modd) \rightarrow
D(\mathcal{O}_Y\text{-}\modd),$$
cf. \cite[\S3.1]{Lipman-NotesOnDerivedFunctorsAndGrothendieckDuality}.
When $f$ is a concentrated 
map of schemes $f_*$ restricts to a functor 
$D_{\text{qc}}(X) \rightarrow D_{\text{qc}}(Y)$
\cite[Prop. 3.9.2]{Lipman-NotesOnDerivedFunctorsAndGrothendieckDuality}. 
If $X$ and $Y$ are noetherian and $f$ is proper\footnote{In 
a non-noetherian world one can work with a more general notion of 
\em a quasi-proper scheme map, \rm cf.
\cite[\S4.3]{Lipman-NotesOnDerivedFunctorsAndGrothendieckDuality}. } 
then $f_*$ restricts to a functor $D(X) \rightarrow D(Y)$
\cite[Th\'eor\`eme 2.2.1]{IllusieConditionsDeFinitudeRelative}. 

The derived inverse image functor exists as a functor
$$ f^*(-)\colon 
D(\mathcal{O}_Y\text{-}\modd) \rightarrow
D(\mathcal{O}_X\text{-}\modd),$$
cf. \cite[\S3.1]{Lipman-NotesOnDerivedFunctorsAndGrothendieckDuality}.
When $f$ is a concentrated map of schemes $f^*$ restricts to a
functor $D_{\text{qc}}(Y) \rightarrow D_{\text{qc}}(X)$
\cite[Prop. 3.9.1]{Lipman-NotesOnDerivedFunctorsAndGrothendieckDuality}. 
If $X$ and $Y$ are locally noetherian and $f$ is perfect, 
then $f^*$ restricts to a functor $D(Y) \rightarrow D(X)$
\cite[Prop. II4.4]{Hartshorne-Residues-and-Duality}.  

\subsection{Adjunctions and dualities for derived functors}
\label{section-adjunctions-and-dualities-for-derived-functors}

Let $X$ be a ringed space. For any $A \in
D(\mathcal{O}_X\text{-}\modd)$ the functor 
$$ A \otimes (-) \colon  D(\mathcal{O}_X\text{-}\modd) \rightarrow
D(\mathcal{O}_X\text{-}\modd) $$
is left adjoint to functor 
$$ \rder\shhomm_X(A, -) \colon  D(\mathcal{O}_X\text{-}\modd) \rightarrow
D(\mathcal{O}_X\text{-}\modd), $$
cf. \cite[Prop. 2.6.1]{Lipman-NotesOnDerivedFunctorsAndGrothendieckDuality}.

For any $A \in D(\mathcal{O}_X\text{-}\modd)$ denote by $A^\vee$ the
object $\rder\shhomm_X(A,\mathcal{O}_X) \in D(\mathcal{O}_X\text{-}\modd)$.
There is a natural morphism $A \rightarrow A^{\vee\vee}$ which 
is an isomorphism for any $A \in D_{\text{perf}}(X)$
\cite[Prop. 7.2]{IllusieGeneralitesSurLesConditionsDeFinitudeDansLesCategoriesDerivees}.
So $(-)^\vee$ restricts to 
a self-inverse category equivalence $D_{\text{perf}}(X) \rightarrow
D_{\text{perf}}(X)^{opp}$ giving us \em the duality functor for perfect complexes\rm.  

For any $A \in D_{\text{perf}}(X)$ there is a canonical isomorphism 
$A^\vee \otimes (-) \simeq \rder\shhomm_X(A, -)$, see 
\S\ref{section-standard-relations-between-derived-functors}\eqref{item-perfect-objects-and-rhom}, so  
$$ A \otimes (-)\colon  D(\mathcal{O}_X\text{-}\modd) \rightarrow
D(\mathcal{O}_X\text{-}\modd) $$
is both the left and the right adjoint of functor
$$ A^\vee \otimes (-)\colon  D(\mathcal{O}_X\text{-}\modd) \rightarrow
D(\mathcal{O}_X\text{-}\modd).$$

Let now $Y$ be another ringed space and let $f\colon X \rightarrow Y$
be a map of ringed spaces. Then functor 
$$ f^*(-) \colon  D(\mathcal{O}_Y\text{-}\modd) \rightarrow
D(\mathcal{O}_X\text{-}\modd) $$
is left adjoint to functor 
$$ f_*(-) \colon  D(\mathcal{O}_X\text{-}\modd) \rightarrow
D(\mathcal{O}_Y\text{-}\modd), $$
cf. \cite[Prop. 3.2.1]{Lipman-NotesOnDerivedFunctorsAndGrothendieckDuality}.

Suppose now that $X$ and $Y$ are concentrated schemes and let $f\colon
X \rightarrow Y$ be a scheme map. Then the functor 
$$ f_*(-) \colon  D_\text{qc}(X) \rightarrow D_{\text{qc}}(Y) $$
has a right adjoint which we denote as
$$ f^\times (-) \colon D_\text{qc}(Y) \rightarrow D_{qc}(X), $$
cf. \cite[Theorem 4.1]{Lipman-NotesOnDerivedFunctorsAndGrothendieckDuality}
or 
\cite[\S4]{Neeman-TheGrothendieckDualityTheoremViaBousfieldsTechniquesAndBrownsRepresentability}. 

To state the rest of the Grothendieck duality results in their full
presently known generality we would have to introduce a number of 
notions (pseudo-coherence, quasi-properness, etc.) which are only
meaningfully different from well-established ones in non-noetherian context.
Since the main bulk of this paper deals with schemes of finite type over 
a field, we prefer to state these results for noetherian schemes only 
and refer the reader to 
\cite[\S4]{Lipman-NotesOnDerivedFunctorsAndGrothendieckDuality} 
for a more general story. 

So let $X$ and $Y$ be noetherian schemes and let $f\colon X \rightarrow Y$ 
be a separated scheme map of finite type. 
Adjunction $(f_*, f^\times)$ induces a natural morphism
$\delta_f\colon f_* \rder\shhomm_X(A, f^\times B) \rightarrow
\rder\shhomm_Y(f_*A,B),$
see \S\ref{section-standard-relations-between-derived-functors}\eqref{item-sheafified-grothendieck-duality-morphism},
often referred to as \em the sheafified Grothendieck duality morphism\rm. 
For $\delta_f$ to be an isomorphism we need $f^\times$ 
to commute with restriction to open sets of $Y$
\cite[\S4.6]{Lipman-NotesOnDerivedFunctorsAndGrothendieckDuality}.
When $f$ is proper $f^\times$ commutes with $\tor$-independent base change
for all objects in $D^+_{qc}(Y)$ 
and so $\delta_f$ 
is an isomorphism for all $A \in D_{qc}(X)$ and $B \in D^+_{qc}(Y)$
\cite[\S4.4]{Lipman-NotesOnDerivedFunctorsAndGrothendieckDuality}. 
If $f$ is also perfect, then $f^\times$
commutes with $\tor$-independent base change for all of $D_{qc}(Y)$ and so
$\delta_f$ is an isomorphism for all $A, B \in D_{qc}(Y)$
\cite[Theorem
4.7.4]{Lipman-NotesOnDerivedFunctorsAndGrothendieckDuality}.
Moreover, the natural map
$\chi_f\colon f^\times(A) \otimes f^*(B) \xrightarrow{\sim} f^\times(A
\otimes B)$, cf. 
\S\ref{section-standard-relations-between-derived-functors}\eqref{item-f^!-f!xf^*-morphism}, 
is an isomorphism for all $A, B \in D_{qc}(X)$
\cite[\S5]{Neeman-TheGrothendieckDualityTheoremViaBousfieldsTechniquesAndBrownsRepresentability}.

By a result of Nagata any separated map of finite type
between noetherian schemes decomposes as an open immersion followed by
a proper map (\cite{Nagata-ImbeddingOfAnAbstractVarietyInACompleteVariety}, or 
\cite{Vojta-NagatasEmbeddingTheorem} for a more modern exposition).
So to make $(-)^\times$ commute with flat base change
we can try and modify its behaviour over open immersions. Indeed,
there is a unique way to paste $(-)^\times$ over proper maps with 
$(-)^*$ over \'etale maps in a way compatible with \'etale base change
of $(-)^\times$ (see \cite{Lipman-NotesOnDerivedFunctorsAndGrothendieckDuality}, Theorem
4.8.1 for more detail). The result is the pseudo-functor $(-)^!$, \em 
Deligne's twisted inverse image pseudo-functor, \rm which associates 
to any finite-type separated map $f\colon X \rightarrow Y$ of noetherian schemes 
a functor $f^!\colon D^+_{qc}(Y) \rightarrow D^+_{qc}(X)$ 
with a number of nice properties:
\begin{enumerate}
\item $f^! = f^\times|_{D^+_{qc}}$ when $f$ is proper and 
$f^!= f^*|_{D^+_{qc}}$ when $f$ is \'etale.
\item For any $f$ functor $f^!$ commutes with $\tor$-independent
base change \cite[Theorem 4.8.3]{Lipman-NotesOnDerivedFunctorsAndGrothendieckDuality}.
\item For perfect $f$ functor $f^!$ restricts to a functor $D(Y)
\rightarrow D(X)$
\cite[Remark 2.1.5]{AvramovIyengarLipman-ReflexivityAndRigidityForComplexesIISchemes}.  
\item There exists, as explained
in \cite[\S4.9.1]{Lipman-NotesOnDerivedFunctorsAndGrothendieckDuality}, 
for all $A \in D^+_{qc}(X)$ a natural morphism  
\begin{align} \label{eqn-f^!-to-f^*O_Y-f^*}
f^!(\mathcal{O}_Y) \otimes f^* (A) \rightarrow f^!(A).
\end{align}
If $f$ is perfect then \eqref{eqn-f^!-to-f^*O_Y-f^*} 
is an isomorphism \cite[Theorem 4.9.4]
{Lipman-NotesOnDerivedFunctorsAndGrothendieckDuality}
and the morphism 
\begin{align} \label{eqn-f^*-to-rhom-f^!O_Y-f^!} 
f^* (A) \rightarrow \rder\shhomm_X\left(f^!(\mathcal{O}_Y),
f^!(A)\right) 
\end{align}
right adjoint to \eqref{eqn-f^!-to-f^*O_Y-f^*} with respect to 
$f^!(\mathcal{O}_Y) \otimes (-)$ 
is also an isomorphism
\cite[Lemma 2.1.10]{AvramovIyengarLipman-ReflexivityAndRigidityForComplexesIISchemes}. 
\item If $f$ is a regular immersion of codimension $n$, then 
$f^!(\mathcal{O}_Y) = \omega_{X/Y}[-n]$ where 
$\omega_{X/Y}$ is the top wedge power of the normal bundle
$\mathcal{N}_{X/Y}$ \cite[Cor. III.7.3]{Hartshorne-Residues-and-Duality}. 
\item If $f$ is smooth of relative dimension $n$, then 
$f^!(\mathcal{O}_Y) = \omega_{X/Y}[n]$ where $\omega_{X/Y}$ is the
top wedge power of the sheaf $\Omega^1_{X/Y}$ of relative differentials
\cite[Theorem 3]{Verdier-BaseChangeForTwistedInverseImageOfCoherentSheaves}. 
\end{enumerate}
When $f$ is both perfect and proper, then $f^! = f^\times|_{D^+_{qc}}$
\em and \rm all the above properties of $f^!$ apply to  
the whole of $f^\times\colon D_{qc} \rightarrow D_{qc}$. We do not 
therefore distinguish between $f^!$ and $f^\times$ when $f$ is perfect 
and proper.  

If $f$ is proper the RHS of \eqref{eqn-f^*-to-rhom-f^!O_Y-f^!}, as
a functor in $A$, has left adjoint 
$f_*\left(f^!  \mathcal{O}_Y \otimes (-)\right)$.  
If $f$ is also perfect we denote this functor by $f_!$ and the fact that 
$\eqref{eqn-f^*-to-rhom-f^!O_Y-f^!}$
is an isomorphism implies immediately that $f_!\colon D_{qc}(X)
\rightarrow D_{qc}(Y)$ is the left adjoint of 
$f^*\colon D_{qc}(Y) \rightarrow D_{qc}(X)$ and the 
adjunction counit $f_! f^* \rightarrow \id$ is the composition 
$$ f_! f^*(\text{-}) = f_*(f^!(\mathcal{O}_Y) \otimes f^*(\text{-})) 
\xrightarrow{\eqref{eqn-f^!-to-f^*O_Y-f^*}} f_* f^!(\text{-}) 
\xrightarrow{\text{adj. counit}} \id. $$ 

Finally, let $X$ be a separated scheme of finite type over 
a field $k$ and let $\pi_k\colon X \rightarrow k$ be the structure
morphism. The functor $\rder\shhomm_X\left(-,\pi^!_k k\right)$
restricts to a self-inverse category equivalence $D(X) \rightarrow
D(X)^{opp}$, \em the global \footnote{ I.e. over a point. One can obtain 
duality theories on $X$ relative to any separated, finite-type map
$\pi_S\colon X \rightarrow S$ with $S$ noetherian, but only after 
restricting to objects of $D(X)$ perfect over $S$ (see \cite{IllusieConditionsDeFinitudeRelative}, Cor. 4.9.2 etc.). 
Since the objects perfect over a point are precisely the complexes with
bounded and coherent cohomologies, the global duality works for all of $D(X)$.} 
Grothendieck duality functor $D_{X/k}$\rm. For any separated
finite-type map $f\colon X \rightarrow Y$ between two schemes of
finite-type over $k$, the duality $D_{\bullet/k}$ interchanges $f^*$ and $f^!$
\cite[Prop. 4.10.1]{Lipman-NotesOnDerivedFunctorsAndGrothendieckDuality}. 
For proper $f$ the dual of $f_*$ under $D_{\bullet/k}$ is
$f_*$ itself - this is precisely the sheafified Grothendieck duality 
isomorphism. 

\subsection{Standard relations between derived functors}
\label{section-standard-relations-between-derived-functors}

There exists a number of well-known morphisms 
and isomorphisms between compositions of the derived functors 
listed in Sections 
\ref{section-derived-functors-and-derived-categories} and 
\ref{section-adjunctions-and-dualities-for-derived-functors}. 
Here we compile for the convenience of the reader 
a list of such elementary relations employed throughout this paper. 

For a number of these morphisms of derived functors 
we say below that they are compatible with the corresponding 
natural morphisms for sheaves. For full detail on this the reader
should consult the reference we quote for each result, but 
roughly we mean the following. A natural 
transformation of compositions of derived functors 
\begin{align}
\label{eqn-derived-functor-transformation}
\rder f_1 \circ \dots \circ \rder f_n \rightarrow \rder g_1 \circ \dots \circ \rder g_m  
\end{align}
is said to be compatible with a natural transformation of 
compositions of the underlying abelian category functors 
\begin{align}
\label{eqn-abelian-functor-transformation}
f_1 \circ \dots \circ f_n \rightarrow g_1 \circ \dots \circ g_m 
\end{align}
if the following diagram commutes
\begin{align}
\xymatrix{
Q \circ f_1 \circ \dots \circ f_n 
\ar[r]^{\eqref{eqn-abelian-functor-transformation}}
\ar[d] &
Q \circ g_1 \circ \dots \circ g_m 
\ar[d] \\
\rder f_1 \circ \dots \circ \rder f_n \circ Q
\ar[r]^{\eqref{eqn-derived-functor-transformation}} &
\rder g_1 \circ \dots \circ \rder g_m \circ Q
}
\end{align}
where $Q$ denotes localisation functor from each chain homotopy category to the 
corresponding derived category and the vertical arrows are composed from
the natural transformations $Q \circ f_i \rightarrow \rder f_i \circ Q$ and 
$Q \circ g_i \rightarrow \rder g_i \circ Q$ that $\rder f_i$ and $\rder g_i$ 
come equipped with by the definition of a right derived functor.
Compositions of left-derived functors are treated analogously. 

\begin{enumerate}

\item \em Commutativity and associativity of tensor product\rm. Let $X$ be a ringed space. Then for any $A,B,C \in
D(\mathcal{O}_X\text{-}\modd)$ there exist unique natural isomorphisms
\begin{align}
A \otimes B \xrightarrow{\sim} B \otimes A
\end{align}
and
\begin{align}
\label{eqn-associativity-of-tensor-product-definition}
\rho\colon 
(A \otimes B) \otimes C 
\xrightarrow{\sim} 
A \otimes B \otimes C 
\xrightarrow{\sim} 
A \otimes (B \otimes C)
\end{align}
which are functorial in $A$, $B$ and $C$ and which are
compatible with the corresponding natural isomorphisms for sheaves 
\cite[\S2.5.7 and \S2.5.9]{Lipman-NotesOnDerivedFunctorsAndGrothendieckDuality}.

\item \em Sheafified ($A \otimes (-)$, $\rder\shhomm(A,-)$) adjunction\rm.  
Let $X$ be a ringed space. Then for any $A,B,C \in
D(\mathcal{O}_X\text{-}\modd)$ there exist unique natural isomorphism
\begin{align}
\label{eqn-sheafified-tensor-product-adjunction}
\upsilon_A\colon
\rder\shhomm_X\left(A \otimes B, C \right)
\xrightarrow{\sim}
\rder\shhomm_X\left(B, \rder\shhomm_X\left(A,C\right) \right)
\end{align}
compatible with the corresponding natural isomorphism for sheaves
\cite[Prop. 2.6.1]{Lipman-NotesOnDerivedFunctorsAndGrothendieckDuality}.

Applying the derived global sections functor to 
\eqref{eqn-sheafified-tensor-product-adjunction} produces 
the adjunction isomorphism for the pair $\left(A \otimes -, 
\rder\shhomm_X(A, -)\right)$. We call its counit 
the \em evaluation map \rm of $A$ and denote it by
\begin{align}
\label{eqn-evaluation-map} 
\evmap_A\colon A \otimes \rder\shhomm_X(A, -) \rightarrow \id. 
\end{align}
An important instance is the morphism 
$A \otimes A^\vee \xrightarrow{\evmap_A} \mathcal{O}_X$ obtained
by applying $\evmap_A$ to $\mathcal{O}_X$. 

\item \em Perfect objects and $\rder\shhomm$. \rm 
\label{item-perfect-objects-and-rhom}
Let $X$ be a ringed space. For any $A,B,C \in 
D(\mathcal{O}_X\text{-}\modd)$ define 
\begin{align}
\label{eqn-bringing-an-object-into-rhom-bracket}
\xi\colon
\rder\shhomm_X\left(A,B\right) \otimes C
\longrightarrow 
\rder\shhomm_X\left(A,B \otimes C \right)
\end{align}
to be the right adjoint with respect to $A \otimes (-)$ of the
composition
\begin{align}
\label{eqn-bringing-an-object-into-rhom-bracket-adjoint}
A \otimes \left(\rder\shhomm_X\left(A,B\right) \otimes C\right)
\xrightarrow{\rho^{-1}}
\left(A \otimes \rder\shhomm_X\left(A,B\right) \right) \otimes C
\xrightarrow{\evmap_A}
B \otimes C.
\end{align}
If either of $C$ or $A$ belong to $D_\text{perf}(X)$, then 
$\xi$ is an isomorphism \cite[Lemma
1.4.6]{AvramovIyengarLipman-ReflexivityAndRigidityForComplexesIISchemes}.
In particular, for any $E \in  D_\text{perf}(X)$ we have an isomorphism 
\begin{align}
\label{eqn-perfect-project-rhom-is-tensor-with-dual}
\xi_E \colon
E^\vee \otimes (-) \xrightarrow{\sim} \rder \shhomm_X(E, -) 
\end{align}
of functors 
$D(\mathcal{O}_X\text{-}\modd) \rightarrow D(\mathcal{O}_X\text{-}\modd)$. 

The adjunction $\left(E \otimes -, \rder\shhomm_X(E, -)\right)$  
induces via $\xi_E$
an adjunction $(E \otimes -, E^\vee \otimes -)$
whose adjunction co-unit we also denote by $\ev_E$:
\begin{align}
\label{eqn-E-x-E^vee-x-adjunction-counit} 
E \otimes (E^\vee \otimes -) 
\xrightarrow{\xi_E}
E \otimes \rder \shhomm_X(E, -) 
\xrightarrow{\ev_E} \id.
\end{align}

\item \em $\mathcal{O}_X$-reflexivity for perfect objects\rm.  
Let $X$ be a ringed space. For any $A,B \in 
D(\mathcal{O}_X\text{-}\modd)$ define 
\begin{align}
\label{eqn-reflexivity-morphism}
\theta_{A,B} \colon
A \longrightarrow 
\rder\shhomm_X\left(\rder\shhomm_X\left(A,B\right),B \right)
\end{align}
to be the right adjoint with respect to
$\rder\shhomm_X\left(A,B\right) \otimes (-)$ of 
$$A \otimes \rder\shhomm_X\left(A,B\right)
\xrightarrow{\ev_A} B.$$ 
If $B = \mathcal{O}_X$ the resulting morphism 
\begin{align}
\label{eqn-dualizing-morphism}
\theta_A \colon A \rightarrow A^{\vee\vee}
\end{align}
an isomorphism for all $A \in D_\text{perf}(X)$
\cite[Prop 1.4.4]{AvramovIyengarLipman-ReflexivityAndRigidityForComplexesIISchemes}. 

Let $E \in D_\text{perf}$. 
The adjunction $\left(E^\vee \otimes -, E^{\vee\vee}\otimes -\right)$  
induces via the isomorphism 
$E \xrightarrow{\theta_{E}}E^{\vee\vee}$  
an adjunction $(E^\vee \otimes -, E \otimes -)$
whose adjunction co-unit we denote by $\ev_{E^\vee}$:
\begin{align}
\label{eqn-E^vee-x-E-adjunction-counit} 
E^\vee \otimes (E \otimes -) 
\xrightarrow{\theta_{E}}
E^\vee \otimes (E^{\vee\vee} \otimes -) 
\xrightarrow{\ev_{E^\vee}} \id.
\end{align}

\item \em Pseudofunctoriality of direct and inverse image\rm.  Let
$X$, $Y$, $Z$ be ringed spaces and 
$X \xrightarrow{f} Y \xrightarrow{g} Z$ be maps of 
ringed spaces. There exist unique isomorphisms
\begin{align}
\label{eqn-pseudofunctoriality-direct-image}
\eta_{g,f}\colon (g \circ f)_* \xrightarrow{\sim} g_*f_*
\quad \text{of functors }
D(\mathcal{O}_{X}\text{-}\modd) \rightarrow D(\mathcal{O}_{Z}\text{-}\modd)
\end{align}
and 
\begin{align}
\label{eqn-pseudofunctoriality-inverse-image}
\zeta_{g,f}\colon f^* g^* \xrightarrow{\sim} (g \circ f)^* 
\quad \text{of functors }
D(\mathcal{O}_{Z}\text{-}\modd) \rightarrow D(\mathcal{O}_{X}\text{-}\modd)
\end{align}
which are compatible with the corresponding natural isomorphisms for 
sheaves. These isomorphisms give $(-)_*$ and $(-)^*$ the structures
of a covariant and a contravariant pseudofunctor over the category
of ringed spaces 
\cite[\S3.6]{Lipman-NotesOnDerivedFunctorsAndGrothendieckDuality}. 
Specifically, for any map $X \xrightarrow{f} Y$ of ringed spaces we have
\begin{align}
\label{eqn-pseudofunctoriality-relations-identity}
\eta_{\id,f} = \eta_{f,\id} = \id \quad \text{ and } \quad 
\zeta_{\id,f} = \zeta_{f,\id} = \id
\end{align}
and for any maps $X \xrightarrow{f} Y  \xrightarrow{g} Z
\xrightarrow{h} W$ of ringed spaces the following diagrams commute
\begin{align}
\label{eqn-pseudofunctoriality-relations-associativity}
\xymatrix{
(h \circ g \circ f)_* 
\ar[r]^{\eta_{h \circ g,f}}
\ar[d]_{\eta_{h, g \circ f}}
&
(h \circ g)_* f_* 
\ar[d]^{\eta_{h,g}}
\\
h_* (g \circ f)_*
\ar[r]_{h_* \eta_{g,f}}
&
h_* g_* f_*
}
\quad\text{ and }\quad
\xymatrix{
f^* g^* h^*  
\ar[r]^{f^* \zeta_{h,g}}
\ar[d]_{\zeta_{g, f}}
&
f^* (h \circ g)^*
\ar[d]^{\zeta_{h \circ g, f}}
\\
(g \circ f)^* h^*  
\ar[r]_{\eta_{h, g \circ f}}
&
(h \circ g \circ f)^*
}.
\end{align}
We write $\eta_{h,g,f}$ for the corresponding isomorphism 
$(h \circ g \circ f)_* \xrightarrow{\sim} h_* g_* f_*$ and
$\zeta_{h,g,f}$ for the corresponding isomorphism 
$ f^* g^* h^* \xrightarrow{\sim} (h \circ g \circ f)^*$.

\item \label{item-sheafified-f^*-f_*-adjunction} \em Sheafified $(f^*,f_*)$ adjunction\rm. Let $X,Y$ be ringed
spaces and let $X \xrightarrow{f} Y$ be a map of ringed spaces. 
For any $A \in D(\mathcal{O}_Y\text{-}\modd)$ and $B \in
D(\mathcal{O}_X\text{-}\modd)$ there exists a unique bifunctorial isomorphism
\begin{align}
\label{eqn-definition-of-sheafified-f^*-f_*-adjunction}
\tau_f\colon 
f_* \rder\shhomm_X(f^*A,B) \xrightarrow{\sim} 
\rder\shhomm_Y(A, f_*B)
\end{align}
compatible with the corresponding
natural isomorphism for sheaves 
\cite[Prop. 3.2.3]{Lipman-NotesOnDerivedFunctorsAndGrothendieckDuality}.

Applying the derived global sections functor to 
\eqref{eqn-definition-of-sheafified-f^*-f_*-adjunction}
produces an adjunction isomorphism for the pair
$\left(f^*, f_*\right)$. We denote its unit and counit by
\begin{align}
\label{eqn-adjunction-unit-and-counit-for-f^*-f_*}
\beta_f \colon \id \rightarrow f_* f^* 
\quad &\text{ and }\quad 
\gamma_f \colon f^* f_* \rightarrow \id.
\end{align}

The adjunction $\left(f^*, f_*\right)$ is compatible with
pseudofunctoriality in the following sense. Let $X \xrightarrow{f} Y$
and $Y \xrightarrow{g} Z$ be maps of ringed spaces, then 
the following diagrams commute:
\begin{align}
\label{eqn-f^*-f_*-adjunction-and-pseudofunctoriality}
\vcenter{
\xymatrix{
\id
\ar[r]^{\beta_g}
\ar[drr]_{\beta_{g \circ f}}
&
g_* g^* 
\ar[r]^{g_* \beta_f}
&
g_* f_* f^*g^* 
\ar[d]^{\eta^{-1}_{g,f} \circ (g_* f_* \zeta_{g,f})} 
\\
&
&
(g\circ f)_* (g \circ f)^* 
}
} 
\; \text{ and } \;
\vcenter{
\xymatrix{
f^*g^* g_* f_* 
\ar[d]_{\zeta_{g,f} \circ (f^* g^* \eta^{-1}_{g,f})} 
\ar[r]^{f^* \gamma_g}
&
f^* f_*
\ar[r]^{\gamma_f}
&
\id,
\\
(g \circ f)^* (g\circ f)_*  
\ar[urr]_{\gamma_{g \circ f}}
&
&
}
}
\end{align}
see \cite[\S3.6]{Lipman-NotesOnDerivedFunctorsAndGrothendieckDuality}
for more details.

\item \em Monoidal functor structure for inverse image\rm. 
Let $X,Y$ be ringed spaces and let $X \xrightarrow{f} Y$ be a map 
of ringed spaces. 
For any $A,B \in D(\mathcal{O}_Y\text{-}\modd)$ there exists a unique isomorphism
\begin{align}
\label{eqn-definition-of-morphism-nu}
\nu_f \colon f^*(A \otimes B) \xrightarrow{\sim} f^*(A) \otimes f^* B 
\end{align}
functorial in $A$ and $B$ which is compatible with the corresponding
natural isomorphism for sheaves 
\cite[Prop. 3.2.4(i)]{Lipman-NotesOnDerivedFunctorsAndGrothendieckDuality}. 
It is worth noting that as a natural transformation 
of functors in $B$ isomorphism $\nu_f$ is conjugate 
to $\tau_f$ in sense of 
\cite[\S IV.7]{MacLane-CategoriesfortheWorkingMathematician}.

Map $\nu_f$ is compatible with the associativity of the tensor 
product in the following sense. Let $X \xrightarrow{f} Y$ be a map 
of ringed spaces. Then the following diagram 
\begin{align}
\label{eqn-nu-is-compatible-with-the-assotiatity-of-the-tensor-product}
\xymatrix{
f^* \left(\left(A \otimes B\right) \otimes C\right)
\ar[d]_{f^* \rho}
\ar[r]^{\nu_f}
&
f^* \left(A \otimes B\right) \otimes f^*C
\ar[r]^{\nu_f \otimes \id}
&
\left(f^* A \otimes f^*B \right) \otimes f^*C
\ar[d]_{\rho}
\\
f^* \left(A \otimes \left(B\otimes C\right)\right)
\ar[r]_{\nu_f}
& 
f^* A \otimes f^* \left(B\otimes C\right)
\ar[r]_{\id \otimes \nu_f}
&
f^* A \otimes \left(f^* B\otimes f^*C\right)
}
\end{align}
commutes for any $A,B,C \in D(\mathcal{O}_Y\text{-}\modd)$
\cite[\S3.4]{Lipman-NotesOnDerivedFunctorsAndGrothendieckDuality}.

Map $\nu_f$ is compatible with pseudofunctoriality in the following
sense. Let $X \xrightarrow{f} Y$ and $Y \xrightarrow{g} Z$ be maps of
ringed spaces. Then the following diagram commutes
\begin{align}
\label{eqn-map-nu-is-compatible-with-pseudofunctoriality}
\xymatrix{
f^* g^*\left(A \otimes B\right) 
\ar[r]^{f^* \nu_g} 
\ar[d]_{\zeta_{g,f}}
&
f^* \left(g^*A \otimes g^*B\right)
\ar[r]^{\nu_f}
&
f^* g^*A \otimes f^* g^* B
\ar[d]^{\zeta_{g,f} \otimes \zeta_{g,f}}
\\
(g \circ f)^* \left(A \otimes B\right) 
\ar[rr]_{\nu_{g \circ f}}
& &
(g \circ f)^* A \otimes (g \circ f)^* B 
}
\end{align}
for all $A, B \in D(\mathcal{O}_Z\text{-}\modd)$ 
\cite[\S3.6]{Lipman-NotesOnDerivedFunctorsAndGrothendieckDuality}.

\item \em Monoidal functor structure for direct image\rm. 
\label{item-symmetric-monoidal-functor-struct-for-direct-image}
Let $X,Y$ be ringed
spaces and let $X \xrightarrow{f} Y$ be a map of ringed spaces. 
For any $A,B \in D(\mathcal{O}_X\text{-}\modd)$ define morphism 
\begin{align}
\label{eqn-symmetric-monoidal-functor-struct-for-direct-image}
\kappa_f\colon f_*A \otimes f_*B \rightarrow f_*(A \otimes B),
\end{align}
functorial in $A$ and $B$, to be the right adjoint with respect to $f^*$
of the composition 
$$
f^*(f_*A \otimes f_*B) \xrightarrow{\nu_f} f^*f_* A \otimes f^*f_* B
\xrightarrow{\gamma_f \otimes \gamma_f} A \otimes B. $$

Map $\kappa_f$ is compatible with the associativity of the tensor
product and with pseudofunctoriality in a way analogous to map $\nu_f$
\cite[\S3.4 and \S3.6]{Lipman-NotesOnDerivedFunctorsAndGrothendieckDuality}. 
 
\item \em Projection formula\rm. Let $X,Y$ be ringed
spaces and let $X \xrightarrow{f} Y$ be a map of ringed spaces. 
For any $A \in D(\mathcal{O}_X\text{-}\modd)$ and $B \in
D(\mathcal{O}_Y\text{-}\modd)$ define the projection formula morphism  
\begin{align}
\label{eqn-definition-of-projection-formula-morphism}
\alpha_f\colon f_* A \otimes B \rightarrow f_* (A \otimes f^*B) 
\end{align}
to be the right adjoint with respect to $f^*$ of the composition 
$$ f^*(f_* A \otimes B)
\xrightarrow{\nu_f} f^*f_*A \otimes f^* B 
\xrightarrow{\gamma_f \otimes \id} A \otimes f^*B.$$
If $X$ and $Y$ are concentrated schemes, then $\alpha_f$ is
an isomorphism for all $A \in D_{qc}(X)$ and $B \in D_{qc}(Y)$
\cite[Prop. 3.9.4]{Lipman-NotesOnDerivedFunctorsAndGrothendieckDuality}. 

The projection formula is compatible with pseudofunctoriality
in the following sense. Let $X \xrightarrow{f} Y$ and $Y \xrightarrow{g} Z$ 
be maps of ringed spaces. Then the following diagram
\begin{align}
\label{eqn-projection-formula-is-compatible-with-pseudofunctoriality}
\xymatrix{
A \otimes g_* f_* B 
\ar[r]^{\alpha_g} 
&
g_* \left(g^* A \otimes f_*B\right)
\ar[r]^{g_* \alpha_f} 
&
g_* f_*\left(f^* g^* A \otimes B\right) 
\ar[d]^{g_* f_* \left(\zeta_{f,g} \otimes \id\right)}_{\simeq}
\\
A \otimes (g \circ f)_* B
\ar[u]^{\id \otimes \eta_{f,g}}_{\simeq}
\ar[r]_<<<<<{\alpha_{g \circ f}} 
&
(g \circ f)_* \left((g \circ f)^* A \otimes  B \right)
\ar[r]_{\eta_{g \circ f}}^{\simeq} 
& 
g_* f_* \left((g \circ f)^* A \otimes  B \right)
}
\end{align}
commutes for any $A \in
D(\mathcal{O}_Z\text{-}\modd)$ and $B \in
D(\mathcal{O}_X\text{-}\modd)$
\cite[Prop. 3.7.1]{Lipman-NotesOnDerivedFunctorsAndGrothendieckDuality}.

\item \em The sheafified Grothendieck duality morphism. \rm
\label{item-sheafified-grothendieck-duality-morphism}
Let $X$ and $Y$ be concentrated schemes and let 
$X \xrightarrow{f} Y$ be a map of schemes. Denote the 
unit and counit of the $\left(f_*, f^\times\right)$
adjunction by
\begin{align}
\label{eqn-adjunction-unit-and-counit-for-f_*-f^x}
\lambda_f \colon \id \rightarrow f^\times f_*
\quad &\text{ and }\quad
\epsilon_f \colon f_* f^\times \rightarrow \id.
\end{align}

The $(f_*, f^\times)$ adjunction is compatible with
pseudofunctoriality, in the sense that the analogues
of diagrams \eqref{eqn-f^*-f_*-adjunction-and-pseudofunctoriality} 
for $\delta_f$ and $\lambda_f$ also commute, see 
\cite[Cor. 4.1.2]{Lipman-NotesOnDerivedFunctorsAndGrothendieckDuality}
for more details.

Define for any $A \in D_{qc}(X)$ and $B \in D_{qc}(Y)$
\em the sheafified Grothendieck duality morphism \rm  
\begin{align}
\label{eqn-sheafified-Grothendieck-duality-morphism}
\delta_f\colon f_* \rder\shhomm_X(A, f^\times B) \rightarrow
\rder\shhomm_Y(f_*A,B) 
\end{align}
to be the composition  
$$
f_* \rder\shhomm_X(A, f^\times B) 
\xrightarrow{\gamma_f}
f_* \rder \shhomm_X\left(f^*f_* A, f^\times B \right)  
\xrightarrow{\tau_f}
\rder \shhomm\left(f_* A, f_* f^\times B \right)
\xrightarrow{\epsilon_f}
\rder \shhomm\left(f_* A, B \right).
$$
When $X$ and $Y$ are Noetherian and $f$ is proper $\delta_f$ is 
an isomorphism for all $A \in D_\text{qc}(X)$ and $B \in
D^+_\text{qc}(Y)$
\cite[Theorem. 4.4.1]{Lipman-NotesOnDerivedFunctorsAndGrothendieckDuality}. 
If, in addition to the above, $f$ is perfect, $\delta_f$ 
is an isomorphism for all 
$A \in D_\text{qc}(X)$ and $B \in D_\text{qc}(X)$ 
\cite[Theorem 4.7.4]{Lipman-NotesOnDerivedFunctorsAndGrothendieckDuality}.

\item 
\label{item-f^!-f!xf^*-morphism}
Let $X,Y$ be concentrated schemes and let $X \xrightarrow{f} Y$
be a map of schemes. For any $A \in D_{qc}(X)$ and $B \in D_{qc}(Y)$
define morphism 
\begin{align}
\label{eqn-definition-of-f^!-f!xf^*-morphism}
\chi_f\colon f^\times A \otimes f^* B \rightarrow f^\times(A \otimes B) 
\end{align}
functorial in $A$ and $B$ to be the right adjoint with 
respect to $f_*$ of the composition
$$ f_*(f^\times A \otimes f^* B) \xrightarrow{\alpha_f^{-1}}
f_* f^\times A \otimes B \xrightarrow{\epsilon_f \otimes \id}
A \otimes B$$ 
where $\alpha_f^{-1}$ is the inverse of the projection formula
isomorphism. When $f$ is proper and perfect $\chi_f$ is an 
isomorphism 
\cite[Exercise 4.7.3.4(a)]{Lipman-NotesOnDerivedFunctorsAndGrothendieckDuality}. 

\item \em Base change\rm. Let $\sigma$ be a commutative square 
\begin{align}
\vcenter{
\xymatrix{
X' \ar[r]^{g'} \ar[d]_{f'} &
X \ar[d]^f \\
Y' \ar[r]_{g} &
Y
} 
}
\end{align}
of ringed spaces. We define the base change morphism
\begin{align}
\label{eqn-definition-of-base-change-morphism}
\mu_{\sigma}\colon g^* f_* \rightarrow f'_* g'^* 
\end{align}
to be the right adjoint with respect to $f'^*$ of the composition
\begin{align*}
f'^* g^* f_* \xrightarrow{\zeta_{g,f'}} (g \circ f')^* f_* = (f \circ g')^*
f_* \xrightarrow{\zeta^{-1}_{f,g'}} g'^* f^* f_* \xrightarrow{\gamma_f} g'^*
\end{align*}
or, equivalently
\cite[Prop. 3.7.2]{Lipman-NotesOnDerivedFunctorsAndGrothendieckDuality},
the left adjoint with respect to $g_*$ of the composition
\begin{align*}
f_* \xrightarrow{\beta_{g'}} f_* g'_* g'^*
\xrightarrow{\eta^{-1}_{f,g'}} (f \circ g')_* g'^* = (g \circ f')_* g'^*
\xrightarrow{\eta_{g,f'}} g_* f'_* g'^*.
\end{align*}
This defines $\mu_\sigma$ as a morphism of functors 
$D(\mathcal{O}_X\text{-}\modd) \rightarrow
D(\mathcal{O}_Y'\text{-}\modd)$.
When $\sigma$ is a square of concentrated schemes the base change
map restricts to a morphism of functors $D_{qc}(X) \rightarrow D_{qc}(Y')$. 

We use $\sigma^T$ to denote the transposed square 
\begin{align}
\vcenter{
\xymatrix{
X' \ar[r]^{f'} \ar[d]_{g'} &
Y'\ar[d]^g \\
X \ar[r]_{f} &
Y.
}
}
\end{align}
In particular, we denote by $\mu_{\sigma^T}$ the base change map 
$f^* g_* \rightarrow f'^* g'_*$ for $\sigma^T$.

If the restriction of $\mu_{\sigma}$ to complexes with quasi-coherent
cohomology is an isomorphism, then $\sigma$ is said to be 
\em independent\rm. A fiber-square of concentrated schemes is independent 
if and only if it is \em$\tor$-independent\rm, i.e. for 
any $x \in X$ and $y' \in Y'$ such that $f(x) = g(y') = y \in Y$ we have
\begin{align}
\tor^i_{\mathcal{O}_{Y,y}}(\mathcal{O}_{X,x},\mathcal{O}_{Y',y'}) = 0
\quad \text{ for all } i > 0,
\end{align}
cf. \cite[Theorem 3.10.3]{Lipman-NotesOnDerivedFunctorsAndGrothendieckDuality}.
In particular, a fiber-square of concentrated schemes is independent 
if $f$ or $g$ are flat. Another good reference for the above material is 
\cite[\S2.4]{Kuznetsov-HyperplaneSectionsAndDerivedCategories}, 
where the proofs are carried out via computations with 
the underlying Fourier-Mukai kernels. 
\end{enumerate}

\subsection{Further relations}

To prove our main results in Section \ref{section-fmt-adjunction-morphism} 
we need three technical results which we could not find in the literature. 
The first two state that the projection formula commutes with certain 
adjunction units and counits of the direct image functor. 

\begin{lemma} \label{lemma-projection-formula-commutes-with-adjunction-1}
Let $X \xrightarrow{g} Y \xrightarrow{f} Z$ be maps of ringed spaces.
Let $A \in D(\mathcal{O}_{Y}\text{-}\modd)$ and $B \in
D(\mathcal{O}_{Z}\text{-}\modd)$. 
Then the following diagram commutes:
\begin{align} \label{eqn-proj-form-commutes-with-adjunction}
\vcenter{
\xymatrix{ 
f_* A \otimes B \ar[rr]^{f_* \beta_g \otimes \id}
\ar[d]_{\alpha_f} & & f_* g_* g^* A \otimes B \ar[d]^{f_*\alpha_{g} \circ
\alpha_{f}} \\
f_* \left( A \otimes f^* B \right) 
\ar[r]_>>>>>{f_* \beta_g} &
f_* g_* g^* \left(A \otimes  f^* B \right) 
\ar[r]_{f_* g_* \nu_g} &
f_* g_* \left( g^* A \otimes g^* f^* B \right) .
}
}
\end{align}
\end{lemma}
\begin{proof}
By functoriality of $\alpha_f$ it suffices to show that the square
\begin{align*}
\xymatrix{
f_* \left( A \otimes f^* B \right) 
\ar[d]_{f_* \beta_g}
\ar[rr]^{f_*(\beta_g \otimes f^*\id)}
&  &
f_* \left( g_* g^* A \otimes f^* B \right)
\ar[d]^{f_*\alpha_{g}} 
\\
f_* g_* g^* \left(A \otimes  f^* B \right) 
\ar[rr]_{f_* g_* \nu_g} & &
f_* g_* \left( g^* A \otimes g^* f^* B \right)
} 
\end{align*}
commutes. This square is the image under $f_*$ of the square
\begin{align}
\label{eqn-proj-form-commutes-with-adjunction-after-alpha-f}
\vcenter{
\xymatrix{
A \otimes f^* B 
\ar[d]_{\beta_g}
\ar[rr]^{\beta_g \otimes \id}
&  &
g_* g^* A \otimes f^* B 
\ar[d]^{\alpha_{g}} 
\\
g_* g^* \left(A \otimes  f^* B \right) 
\ar[rr]_{g_* \nu_g} & &
g_* \left( g^* A \otimes g^* f^* B \right).
}
}
\end{align}

To show that \eqref{eqn-proj-form-commutes-with-adjunction-after-alpha-f} 
commutes we show that its left adjoint 
with respect to $g_*$ commutes.
By definition of $\alpha_g$ its left adjoint 
with respect to $g_*$ is $(\gamma_g \otimes \id) \circ \nu_g$. 
So the left adjoint with respect to $g_*$ of 
\eqref{eqn-proj-form-commutes-with-adjunction-after-alpha-f} is 
\begin{align*}
\xymatrix{
g^*\left(A \otimes f^* B\right)
\ar[rr]^{g^*(\beta_g \otimes \id)}
\ar[rrrd]_{\nu_g}
& & 
g^*\left(g_* g^* A \otimes f^* B\right) 
\ar[r]^{\nu_g} &
g^* g_* g^* A \otimes g^* f^* B
\ar[d]^{\gamma_g \otimes \id}
\\
& & &
g^* A \otimes g^* f^* B 
}
\end{align*}
and by functoriality of $\nu_g$ it suffices to show that
the following composition is the identity morphism:
$$ 
g^* A \otimes g^* f^* B 
\xrightarrow{g^*\beta_g \otimes \id}
g^*g_* g^* A \otimes g^* f^* B 
\xrightarrow{\gamma_g \otimes \id}
g^* A \otimes g^* f^* B. $$
Rewrite it as $(g^* \beta_g \circ \gamma_g) \otimes \id$. 
Since $\beta_g$ and $\gamma_g$ are the unit and the counit 
of the adjunction $(g^*,g_*)$, the morphism 
$g^*A \xrightarrow{g^* \beta_g \circ \gamma_g} g^* A$ is 
the identity morphism. The result follows.
\end{proof}

\begin{lemma} \label{lemma-projection-formula-commutes-with-adjunction-2}
Let $X$, $Y$, $Z$ be concentrated schemes and $X \xrightarrow{g} Y 
\xrightarrow{f} Z$ be scheme maps. Let $A \in D_{\text{qc}}(Y)$ and 
$B \in D_{\text{qc}} (Z)$. Then the following diagram commutes:
\begin{align} \label{eqn-proj-form-commutes-with-adjunction-2}
\vcenter{
\xymatrix{ 
f_* g_* g^\times A \otimes B \ar[rr]^{f_* \epsilon_g \otimes \id}
\ar[d]_{f_* \alpha_{g} \circ \alpha_f} & &
f_* A \otimes B \ar[d]^{\alpha_f} \\
 f_* g_* \left( g^\times A \otimes g^* f^* B \right) \ar[r]_{f_* g_* \chi_g} &
f_* g_* g^\times \left(A \otimes f^* B \right) \ar[r]_>>>>>{f_* \epsilon_g} &
f_* \left(A \otimes f^* B \right) .
}
}
\end{align}
\end{lemma}
\begin{proof}
The proof is analogous to that of Lemma
\ref{lemma-projection-formula-commutes-with-adjunction-1}. 
By functoriality of $\alpha_f$ it suffices to show 
that the image under $f_*$ of 
\begin{align*} 
\xymatrix{ 
g_* g^\times A \otimes f^*B 
\ar[rrd]^{\epsilon_g \otimes \id}
\ar[d]_{\alpha_{g}} & &
\\
g_* \left( g^\times A \otimes g^* f^* B \right) \ar[r]_{g_* \chi_g} &
g_* g^\times \left(A \otimes f^* B \right) \ar[r]_>>>>>{\epsilon_g} &
A \otimes f^* B.
}
\end{align*}
commutes. Since $\alpha_g$ is an isomorphism, this is equivalent to the diagram
\begin{align*} 
\xymatrix{ 
g_* g^\times A \otimes f^*B 
\ar[rrd]^{\epsilon_g \otimes \id}
\\
\ar[u]^{\alpha^{-1}_{g}} 
g_* \left( g^\times A \otimes g^* f^* B \right) \ar[r]_{g_* \chi_g} &
g_* g^\times \left(A \otimes f^* B \right) \ar[r]_>>>>>{\epsilon_g} &
A \otimes f^* B.
}
\end{align*}
commuting. But as $\epsilon_g$ is the adjunction counit, 
the composition $\epsilon_g \circ g_* \chi_g$ is the left adjoint
of $\chi_g$ with respect to $g^\times$. By the definition  
of $\chi_g$ this left adjoint is precisely 
$(\epsilon_g \otimes \id) \circ \alpha^{-1}_g$. The result follows.  
\end{proof}

The third result shows that for a perfect object $E$
the adjunction co-units for $E \otimes (-)$ commute with
the associativity of the tensor product:
\begin{lemma}
\label{lemma-tensor-product-adjunction-and-associativity} 
Let $X$ be a ringed space. Then for any $A \in
D(\mathcal{O}_X\text{-}\modd)$ and $E \in D_{perf}(X)$ the following
diagrams commute
\begin{align}
\vcenter{\xymatrix{
E \otimes \left(E^\vee \otimes A \right) 
\ar[rr]^<<<<<<<<<<<{\ev_E}
\ar[d]_{\rho^{-1}}
& \quad &
A 
\ar[d]^{\id}
\\
\left(E \otimes E^\vee \right) \otimes A
\ar[rr]_<<<<<<<<<<{\ev_E(\mathcal{O}_X) \otimes \id} 
& \quad &
A
}}
\quad \text{ and } \quad
\vcenter{\xymatrix{
\left(E^\vee \otimes E \right) \otimes A 
\ar[rr]^<<<<<<<<<<{\ev_E(\mathcal{O}_X) \otimes \id} 
\ar[d]_{\rho}
& \quad &
A 
\ar[d]^{\id}
\\
E^\vee \otimes \left(E \otimes A \right)
\ar[rr]_<<<<<<<<<<<{\ev_{E^\vee}}
& \quad &
A.
}}
\end{align}
\end{lemma}
\begin{proof}
The adjunction counit 
$E \otimes (E^\vee \otimes A) \xrightarrow{\ev_E} A$ 
was defined as the composition 
\begin{align*}
\xymatrix{
E \otimes (E^\vee \otimes A)
\ar[rr]_>>>>>>>>>>>>>>{\sim}^>>>>>>>>>>>>>>{\id \otimes \xi_E}
& \quad &
E \otimes \rder \shhomm(E,A) 
\ar[rr]^>>>>>>>>>>>>>{\ev_E}
& \quad \quad &
A
}.
\end{align*}
Therefore its right adjoint with respect to $E \otimes (-)$ is
isomorphism $\xi_E$. 
But isomorphism $\xi_E$ was defined to be the right adjoint 
with respect to $E \otimes (-)$ of the composition 
$$
E \otimes \left(E^\vee \otimes A\right)
\xrightarrow{\rho^{-1}}
\left(E \otimes E^\vee \right) \otimes A
\xrightarrow{\ev_E(\mathcal{O}_X) \otimes \id} A .$$
Therefore the left diagram commutes. 

For the right diagram, recall that by its definition the adjunction co-unit 
$E^\vee \otimes (E \otimes A) \xrightarrow{\evmap_{E^\vee}} \id$ 
is
\begin{align*}
\xymatrix{
E^\vee \otimes (E \otimes A) 
\ar[rr]_>>>>>>>>>>>>>>{\sim}^>>>>>>>>>>>>>>{(\id \otimes \theta_E) \otimes \id}
& \quad &
E^\vee \otimes (E^{\vee\vee} \otimes A) 
\ar[rr]^>>>>>>>>>>>>>{\evmap_{E^\vee}}
& \quad \quad &
A
}.
\end{align*}
Since the left diagram commutes and $\rho$ is functorial, 
we can rewrite the composition above as 
\begin{align*}
\xymatrix{
E^\vee \otimes \left(E \otimes A\right)
\ar[r]^{\rho^{-1}}
& 
\left(E^\vee \otimes E\right) \otimes A
\ar[rr]_>>>>>>>>>>>>>>{\sim}^>>>>>>>>>>>>>>{\theta_E}
& \quad &
\left(E^\vee \otimes E^{\vee\vee}\right) \otimes A
\ar[rr]^>>>>>>>>>>>>>{\evmap_{E^\vee}(\mathcal{O}_X)}
& \quad \quad &
A
}.
\end{align*}
To show that the right diagram commutes it now remains only to show
that
\begin{align*}
\xymatrix{
E^\vee \otimes E
\ar[rr]_>>>>>>>>>>>>>>{\sim}^>>>>>>>>>>>>>>{\id \otimes \theta_E}
& \quad &
E^\vee \otimes E^{\vee \vee}
\ar[rr]^>>>>>>>>>>>>>{\evmap_{E^\vee}(\mathcal{O}_X)}
& \quad \quad &
\mathcal{O}_X
}
\end{align*}
is the map 
$E^\vee \otimes E \xrightarrow{\evmap_E(\mathcal{O}_X)} \mathcal{O}_X$. 
The right adjoint of the composition above with respect to 
$E^\vee \otimes (-)$ is just the map 
$E \xrightarrow{\theta_E} E^{\vee \vee}$. 
But $\theta_E$ was defined as the right adjoint with respect to 
$E^\vee \otimes (-)$ of 
$E^\vee \otimes E \xrightarrow{\evmap_E(\mathcal{O}_X)} \mathcal{O}_X$. 
The claim follows.
\end{proof}

Define a morphism 
\begin{align}
\label{eqn-universal-adj-counit-for-tensor-product}
\evmap_E\colon E^\vee \otimes E \otimes (-) \longrightarrow \id
\end{align}
to be the composition 
$$
E^\vee \otimes E \otimes (-) \;
\underset{\sim}{\overset{\rho}{\longrightarrow}} \;
(E^\vee \otimes E) \otimes (-) \;
\xrightarrow{\evmap_E(\mathcal{O}_X)\otimes \id}\; \id.$$
By Lemma \ref{lemma-tensor-product-adjunction-and-associativity}
the canonical isomorphisms identifying $E^\vee \otimes E \otimes -$
with $E^\vee \otimes \left( E \otimes -\right)$
and $E \otimes \left( E^\vee \otimes -\right)$
identify \eqref{eqn-universal-adj-counit-for-tensor-product} 
with the adjunction counits for the adjunctions $(E^\vee \otimes -, E
\otimes -)$ and $(E \otimes -, E^\vee \otimes -)$, respectively. 
We thus abuse notation by speaking of 
\eqref{eqn-universal-adj-counit-for-tensor-product}
as ``the adjunction counit'' for these two adjunctions.

\section{Adjunction morphisms for Fourier-Mukai transforms} 
\label{section-fmt-adjunction-morphism}

 \subsection{Compact case}
\label{section-compact-case}

Let $X_1$ and $X_2$ be a pair of separable schemes of finite 
type over an algebraically closed field $k$ of characteristic $0$ 
with $X_2$ proper. We have the following commutative diagram 
\begin{align} \label{eqn-big-projection-tree}
\vcenter{
\xymatrix{
& & X_1 \times X_2 \times X_1 \ar[ld]_{\pi_{12}} \ar[d]^{\pi_{13}} \ar[rd]^{\pi_{23}} & & \\
& X_1 \times X_2 \ar[ld]_{\pi_1} \ar[rd]^>>>>>>>{\pi_2} & X_1 \times X_1
\ar[lld]_>>>>>>>>>>>>>>{\tilde{\pi}_1}
\ar[rrd]^>>>>>>>>>>>>>>{\tilde{\pi}_2} & X_2 \times X_1 \ar[ld]_>>>>>>>{\pi_2} \ar[rd]^{\pi_1} & \\
X_1 & & X_2 & & X_1
}
}
\end{align}
All the morphisms in it are separated and of finite-type. They are 
also flat, and therefore perfect.  
Moreover, morphisms $\pi_1$ and $\pi_{13}$ are proper.  

\begin{defn}
Let $E$ be a perfect object of $D(X_1 \times X_2)$. 
The \tt Fourier-Mukai transform $\Phi_E$ from $X_1$ to
$X_2$ with kernel $E$ \rm is the functor 
$D_{qc}(X_1) \rightarrow  D_{qc}(X_2)$ given by
$$\Phi_E(-) =  \pi_{2 *} \left(E \otimes \pi^*_1
\left(-\right)\right).$$
\end{defn}

By the adjunctions described in Section 
\ref{section-adjunctions-and-dualities-for-derived-functors}
functor $\Phi_E$ has both left and right adjoints. The left
adjoint $\Phi^{\text{ladj}}_E$ is 
isomorphic to the Fourier-Mukai transform from $X_2$ to $X_1$ with kernel 
$E^\vee \otimes \pi^!_1(\mathcal{O}_{X_1})$. The composition $\Phi^{\text{ladj}}_E \Phi_E$ is 
then isomorphic \cite[Prop 1.3]{Muk81} 
to the Fourer-Mukai transform from $X_1$ to $X_1$
with kernel
$$ Q = \pi_{13 *}\left(\pi_{12}^* E \otimes \pi_{23}^* E^\vee
\otimes \pi_{23}^* \pi^!_1(\mathcal{O}_{X_1}) \right).$$

Let now $\Delta$ denote the diagonal inclusion $X_1 \hookrightarrow X_1
\times X_1$ and, by abuse of notation, let it also denote the induced
inclusion $X_1 \times X_2 \hookrightarrow X_1 \times X_2 \times X_1$, 
so that there is the following fiber square: 
\begin{align} \label{eqn-diagonal-morphism-square}
\vcenter{
\xymatrix{
X_1 \times X_2 \ar@{^{(}->}[r]^>>>>>{\Delta} \ar[d]_{\pi_1} & X_1 \times X_2 \times X_1
\ar[d]^{\pi_{13}} \\
X_1 \ar@{^{(}->}[r]_>>>>>>>>>{\Delta} & X_1 \times X_1
}
}
\end{align}
The identity functor $\id$ is isomorphic to the
Fourier-Mukai transform from $X_1$ to $X_1$ with kernel $\Delta_*
\mathcal{O}_{X_1}$. We now state the main result of this section:

\begin{theorem} \label{theorem-left-adjunction-counit-morphism}
Let $X_1$ and $X_2$ be two separable schemes of finite type over
$k$ with $X_2$ proper. Let $E$ be a perfect object of $D(X_1 \times X_2)$ 
and $\Phi_E$ be a Fourier-Mukai transform 
from $D_{qc}(X_1)$ to $D_{qc}(X_2)$ defined by $E$. 

The adjunction counit $\gamma_E\colon \Phi^{\text{ladj}}_E \Phi_E \rightarrow
\id$ is isomorphic to the morphism of
Fourier-Mukai transforms $D_{qc}(X_1) \rightarrow D_{qc}(X_1)$
induced by the following morphism of their kernels:
\begin{align} \label{eqn-derived-restriction-morphism}
Q = \pi_{13 *}\left(\pi_{12}^* E \otimes \pi_{23}^* E^\vee
\otimes \pi_{23}^* \pi^!_1(\mathcal{O}_{X_1})\right) 
\xrightarrow{\pi_{13 *} \beta_\Delta}
\pi_{13 *} \Delta_* \Delta^* 
\left(\pi_{12}^* E \otimes \pi_{23}^* E^\vee \otimes \pi_{23}^*
\pi^!_1(\mathcal{O}_{X_1})\right) \\
\label{eqn-two-squares-isomorphism}
\pi_{13 *} \Delta_* \Delta^* 
\left(\pi_{12}^* E \otimes \pi_{23}^* E^\vee \otimes \pi_{23}^*
\pi^!_1(\mathcal{O}_{X_1})\right)
\quad \simeq \quad
\Delta_* \pi_{1 *}  
\left(E \otimes E^\vee \otimes \pi^!_1(\mathcal{O}_{X_1})\right)
\\ 
\label{eqn-E-E-dual-trace-morphism}
\Delta_* \pi_{1 *}  
\left(E \otimes E^\vee \otimes \pi^!_1(\mathcal{O}_{X_1})\right) 
\xrightarrow{\Delta_* \pi_{1 *} \ev_E}
\Delta_* \pi_{1 *} \left(\pi^!_1(\mathcal{O}_{X_1})\right)  
\\
\label{eqn-rpi_1-trace-morphism}
\Delta_* \pi_{1 *} \pi^!_1(\mathcal{O}_{X_1})
\xrightarrow{\Delta_* \epsilon_{\pi_1}} \Delta_* \mathcal{O}_{X_1}
\end{align}
where \eqref{eqn-two-squares-isomorphism} is composed of 
isomorphism $\nu_\Delta\colon \Delta^*\left(- \otimes - \right) 
\xrightarrow{\sim}
\Delta^*\left( - \right) \otimes \Delta^* \left( - \right)$ and 
of pseudofunctoriality 
isomorphisms corresponding to the identities 
$\pi_{13} \circ \Delta = \Delta \circ \pi_1$  and $\pi_{12}
\circ \Delta = \pi_{23} \circ \Delta = \id$.
\end{theorem}

We first need the following crucial lemma:

\begin{lemma}
\label{lemma-change-of-base-turns-adjunct-into-diag-restrict}
Let $\sigma$ be the  fiber square
\begin{align}
\vcenter{
\xymatrix{
X_{1} \times X_{2} \times X_{1} 
\ar[r]^{\pi_{12}}
\ar[d]_{\pi_{23}} 
& 
X_{1} \times X_{2}
\ar[d]^{\pi_{2}}
\\
X_{1} \times X_{2}
\ar[r]_{\pi_{2}}
&
X_{2}.
}
}
\end{align}
Then the following diagram of functors commutes:
\begin{align} 
\vcenter{
\xymatrix{
\pi^*_2 \pi_{2 *} \ar[d]_{\mu_\sigma}^{\simeq}
\ar[r]^>>>>>>>>>>{\gamma_{\pi_2}} &
\id  
\ar[d]_{\simeq}^{\eta_{\pi_{23}, \Delta} \circ \zeta^{-1}_{\pi_{12}, \Delta}}
\\
\pi_{23 *} \pi_{12}^* \ar[r]_<<<<<{\pi_{23 *} \beta_\Delta} &
\pi_{23 *} \Delta_* \Delta^* \pi_{12}^*.
}
}
\end{align}
\end{lemma}
\begin{proof}
It suffices to show that the right adjoints with respect to 
$\pi^{*}_{2}$ of the composition
\begin{align}
\label{eqn-base-change-compos-adjunction-unit-for-Delta}
\pi^*_2 \pi_{2 *} 
\xrightarrow{\mu} 
\pi_{23 *} \pi_{12}^*
\xrightarrow{\pi_{23 *} \beta_\Delta} 
\pi_{23 *} \Delta_* \Delta^* \pi_{12}^*
\end{align}
and of the composition
\begin{align}
\label{eqn-adjunction-counit-for-pi2-compos-pseudofunc}
\pi^*_2 \pi_{2 *} 
\xrightarrow{\gamma_{\pi_2}}
\id
\xrightarrow{\eta_{\pi_{23}, \Delta} \circ \zeta^{-1}_{\pi_{12}, \Delta}} 
\pi_{23 *} \Delta_* \Delta^* \pi_{12}^*
\end{align}
coincide. By the definition of morphism $\mu_\sigma$ 
the right adjoint with respect to $\pi^*_2$ 
of \eqref{eqn-base-change-compos-adjunction-unit-for-Delta} is
\begin{align*}
\pi_{2 *} 
\xrightarrow{\pi_{2 *} \beta_{\pi_{12}}}
\pi_{2 *} \pi_{12 *} \pi^*_{12}
\xrightarrow{\eta_{\pi_2, \pi_{23}} \circ \eta^{-1}_{\pi_{2},\pi_{12}}}
\pi_{2 *} \pi_{23*} \pi^*_{12}
\xrightarrow{\pi_{2 *} \pi_{23*} \beta_\Delta}
\pi_{2 *} \pi_{23*} \Delta_* \Delta^* \pi^*_{12}
\end{align*}
which by functoriality of 
$\eta_{\pi_2, \pi_{23}} \circ \eta^{-1}_{\pi_{2},\pi_{12}}$ 
is the same as
\begin{align}
\label{eqn-right-adjoint-of-base-change-compos-adjunction-unit-for-Delta}
\pi_{2 *} 
\xrightarrow{\pi_{2 *} \beta_{\pi_{12}}}
\pi_{2 *} \pi_{12 *} \pi^*_{12}
\xrightarrow{\pi_{2 *} \pi_{12*} \beta_\Delta}
\pi_{2 *} \pi_{12*} \Delta_* \Delta^* \pi^*_{12}
\xrightarrow{\eta_{\pi_2, \pi_{23}} \circ \eta^{-1}_{\pi_{2},\pi_{12}}}
\pi_{2 *} \pi_{23*} \Delta_* \Delta^* \pi^*_{12}.
\end{align}
By pseudofunctoriality of the direct image, 
cf. \eqref{eqn-pseudofunctoriality-relations-associativity},
the morphism of functors
$$
\pi_{2 *} \pi_{12*} \Delta_*
\xrightarrow{\eta_{\pi_2, \pi_{23}} \circ \eta^{-1}_{\pi_{2},\pi_{12}} }
\pi_{2 *} \pi_{23*} \Delta_* $$
is the same as the morphism of functors
$$  
\pi_{2 *} \pi_{12*} \Delta_*
\xrightarrow{\pi_{2_*}(\eta_{\pi_{23}, \Delta} \circ 
\eta^{-1}_{\pi_{12},\Delta})}
\pi_{2 *} \pi_{23*} \Delta_* $$
and we can therefore rewrite
\eqref{eqn-right-adjoint-of-base-change-compos-adjunction-unit-for-Delta}
as
\begin{align}
\label{eqn-right-adjoint-of-base-change-compos-adjunction-unit-for-Delta-2}
\pi_{2 *} \left(
\id
\xrightarrow{\beta_{\pi_{12}}}
\pi_{12 *} \pi^*_{12}
\xrightarrow{\pi_{12*} \beta_\Delta}
\pi_{12*} \Delta_* \Delta^* \pi^*_{12}
\xrightarrow{\eta_{\pi_{23}, \Delta} \circ \eta^{-1}_{\pi_{12},\pi_{\Delta}}}
\pi_{23*} \Delta_* \Delta^* \pi^*_{12} 
\right).
\end{align}

By the compatibility of $\beta$ with pseudofunctoriality as per diagram 
\eqref{eqn-f^*-f_*-adjunction-and-pseudofunctoriality}
we can rewrite
\eqref{eqn-right-adjoint-of-base-change-compos-adjunction-unit-for-Delta-2}
as
\begin{align*}
\pi_{2 *} \left(
\id
\xrightarrow{\beta_{\pi_{12} \circ \Delta}}
(\pi_{12} \circ \Delta)_* (\pi_{12} \circ \Delta)^* 
\xrightarrow{\eta_{\pi_{12}, \Delta} \circ \zeta^{-1}_{\pi_{12},\Delta}}
\pi_{12*} \Delta_* \Delta^* \pi^*_{12}
\xrightarrow{\eta_{\pi_{23}, \Delta} \circ \eta^{-1}_{\pi_{12}, \Delta}}
\pi_{23*} \Delta_* \Delta^* \pi^*_{12}
\right).
\end{align*}
Cancelling out $\eta^{-1}_{\pi_{12}, \Delta} \circ \eta_{\pi_{12},
\Delta}$ and noting that $\beta_{\pi_{12} \circ \Delta} = \id$ 
since $\pi_{12} \circ \Delta = \id$ yields
\begin{align*}
\pi_{2 *} \left(
\id
\xrightarrow{\eta_{\pi_{23}, \Delta} \circ \zeta^{-1}_{\pi_{12},\Delta}}
\pi_{23*} \Delta_* \Delta^* \pi^*_{12}
\right).
\end{align*}
which is clearly the right adjoint of  
\eqref{eqn-adjunction-counit-for-pi2-compos-pseudofunc} with respect
to $\pi^*_2$. The result follows. 
\end{proof}

\begin{proof}[Proof of Theorem \ref{theorem-left-adjunction-counit-morphism}]
Set
$$ Q' = \pi^*_{23} \left(\pi^!_1\mathcal{O}_{X_1} \otimes E^\vee
\right) \otimes \pi^*_{12} E $$ so that $Q = \pi_{13 *} Q'$.  
Since $\pi_{12} \circ \Delta = \pi_{23} \circ \Delta = \id$ we 
have a natural isomorphism 
\begin{align}
\label{eqn-Delta-Q'-iso}
\Delta^* Q' 
\xrightarrow{\nu_\Delta}
\Delta^* \pi^*_{23} \left(\pi^!_1\mathcal{O}_{X_1} \otimes E^\vee
\right) \otimes \Delta^* \pi^*_{12} E 
\xrightarrow{\zeta_{\pi_{23},\Delta} \otimes \zeta_{\pi_{12},\Delta}}
\pi^!_1\mathcal{O}_{X_1} \otimes E^\vee \otimes E.
\end{align}
We therefore define a morphism 
\begin{align}
\label{eqn-Delta-Q'-pi^!_1-O_X1-morphism}
\Delta^* Q' \xrightarrow{\eqref{eqn-Delta-Q'-iso}}
\pi^!_1\mathcal{O}_{X_1} \otimes E^\vee \otimes E
\xrightarrow{\evmap_E} 
\pi^!_1\mathcal{O}_{X_1}.
\end{align}

Let us write the morphism of functors induced by the morphism 
$
Q 
\xrightarrow{  
\eqref{eqn-derived-restriction-morphism}-\eqref{eqn-rpi_1-trace-morphism}
}
\Delta_* \mathcal{O}_{X_1}
$
of FM-kernels as:
\begin{align} \label{eqn-transform-derived-restriction}
\tilde{\pi}_{2 *} \left( \pi_{13 *} Q' \otimes \tilde{\pi}^*_1(-)\right)  
\xrightarrow{\beta_\Delta} 
\tilde{\pi}_{2 *} \left( \pi_{13 *}\Delta_*\Delta^* Q' \otimes \tilde{\pi}^*_1(-)\right)  
\\ 
\label{eqn-transform-pi_13-Delta-to-Delta-pi_1-relabelling}
\tilde{\pi}_{2 *} \left( \pi_{13 *}\Delta_*\Delta^* Q' \otimes \tilde{\pi}^*_1(-)\right)  
\xrightarrow{\eta_{\Delta,\pi_1} \circ \eta^{-1}_{\pi_{13},\Delta}}
\tilde{\pi}_{2 *} \left( \Delta_*\pi_{1 *} \Delta^* Q' \otimes \tilde{\pi}^*_1(-)\right)  
\\ \label{eqn-transform-trace-map}
\tilde{\pi}_{2 *} \left( \Delta_*\pi_{1 *} \Delta^* Q' \otimes \tilde{\pi}^*_1(-)\right)  
\xrightarrow{\eqref{eqn-Delta-Q'-pi^!_1-O_X1-morphism}}
\tilde{\pi}_{2 *} \left( \Delta_*\pi_{1 *} \pi^{!}_1 \mathcal{O}_{X_1} \otimes \tilde{\pi}^*_1(-)\right)  
\\ \label{eqn-transform-adjunction}
\tilde{\pi}_{2 *} \left( \Delta_*\pi_{1 *} \pi^{!}_1 \mathcal{O}_{X_1} \otimes \tilde{\pi}^*_1(-)\right)  
\xrightarrow{\epsilon_{\pi_1}}
\tilde{\pi}_{2 *} \left( \Delta_* \mathcal{O}_{X_1} \otimes \tilde{\pi}^*_1(-)\right)  
\end{align}
On the other hand, $\Phi_E$ is the composition of functors $\pi^*_{1}$,
$E \otimes (-)$ and $\pi_{2 *}$. Each of these functors has a left
adjoint, these adjoints are $\pi_{1 *}(\pi_1^!\mathcal{O}_{X_1} \otimes -)$, 
$E^\vee \otimes (-)$ and $\pi^*_{2}$, respectively. 
Therefore, the adjunction counit 
$\Phi^{\text{ladj}}_E \Phi_E \rightarrow \id$
is the composition of the three corresponding adjunction counits: 
\begin{align} \label{eqn-pi1-pi2-morphism}
\pi_{1 *} \left(\pi^!_1\mathcal{O}_{X_1} \otimes E^\vee \otimes
\pi^*_2 \pi_{2 *} \left(E \otimes \pi^*_1 \left(-\right)\right) \right)
\xrightarrow{\gamma_{\pi_2}}
\pi_{1 *} \left(\pi^!_1\mathcal{O}_{X_1} \otimes E^\vee \otimes
E \otimes \pi^*_1 \left(-\right)\right) \\
\label{eqn-pi1-trace-morphism}
\pi_{1 *} \left(\pi^!_1\mathcal{O}_{X_1} \otimes E^\vee \otimes
E \otimes \pi^*_1 \left(-\right)\right) 
\xrightarrow{\evmap_E} 
\pi_{1 *} \left(\pi^!_1\mathcal{O}_{X_1} \otimes \pi^*_1 \left(-\right)\right) \\
\label{eqn-pi1-trace-adjunction}
\pi_{1 *} \left(\pi^!_1\mathcal{O}_{X_1} \otimes \pi^*_1 \left(-\right)\right) 
\xrightarrow{\epsilon_{\pi_1} \circ \chi_{\pi_1}}
\id
\end{align}
The claim of the theorem is that the composition 
\eqref{eqn-pi1-pi2-morphism}-\eqref{eqn-pi1-trace-adjunction} 
is isomorphic to the composition
\eqref{eqn-transform-derived-restriction}-\eqref{eqn-transform-adjunction}.

Let us clarify some terminology. We say that 
two morphisms of functors $f \rightarrow g$ and $f' \rightarrow g'$ 
are isomorphic if there exist 
connecting isomorphisms $f \xrightarrow{\sim} f'$ and 
$g \xrightarrow{\sim} g'$ such that the diagram 
\begin{align}
\label{eqn-f-g-iso-to-f'-g'}
\xymatrix{
f \ar[r] \ar[d]_{\sim} &  
g \ar[d]^{\sim} \\
f' \ar[r] &
g'
}
\end{align}
commutes. Clearly it is an equivalence relation on the set of all
morphisms between all functors between two given categories. In
particular, it is transitive.  

If we further have a morphism of functors $g \rightarrow h$
which is isomorphic to a morphism of functors $g'' \rightarrow h''$
then $f \rightarrow g \rightarrow h$ is isomorphic
to $f' \rightarrow g' \xrightarrow{\sim} g'' \rightarrow h''$, 
where the connecting isomorphism $g' \xrightarrow{\sim} g''$
is the composition of the inverse of the connecting isomorphism
$g \xrightarrow{\sim} g'$ with the connecting isomorphism 
$g \xrightarrow{\sim} g''$.

Our strategy therefore is to 
consecutively replace the morphisms which compose
\eqref{eqn-pi1-pi2-morphism}-\eqref{eqn-pi1-trace-adjunction} 
by isomorphic ones until we obtain 
\eqref{eqn-transform-derived-restriction}-\eqref{eqn-transform-adjunction}.
However, every time we replace a composant 
by an isomorphic one, we introduce a new connecting
isomorphism. In the end we have to compose a long chain of 
these isomorphisms (each composed of natural isomorphisms detailed in 
\S \ref{section-standard-relations-between-derived-functors})
and simplify the result. It is a mechanical exercise in
pseudofunctoriality of direct and inverse image and the associativity 
of tensor product. To present it in full detail
would be very tedious, the end result being always obvious from 
the start. This had long been lamented in 
the literature, cf. \cite[\S II.6]{Hartshorne-Residues-and-Duality}.
To keep the focus on the substance of a proof
we only state the final result of each such computation 
of a connecting isomorphism, unless something non-trivial is involved.
For our most meticulous readers (and our most inquisitive
referees) we have included in the Appendix an unabbreviated proof, 
where all such computations are carried out in full detail. 

We begin with morphism \eqref{eqn-pi1-pi2-morphism}. By
Lemma \ref{lemma-change-of-base-turns-adjunct-into-diag-restrict} it
is isomorphic to
\begin{footnotesize}
\begin{align}
\label{eqn-pi1-pi2-morphism-transformed-1}
\pi_{1 *} \left(E^\vee \otimes \pi^!_1\mathcal{O}_{X_1}
\otimes \pi_{23*} \pi^*_{12} \left(E \otimes \pi^*_1
\left(-\right)\right)\right) 
\xrightarrow{ \quad\beta_\Delta \quad}
\pi_{1 *} \left(E^\vee \otimes \pi^!_1\mathcal{O}_{X_1}
\otimes \pi_{23*} \Delta_{*} \Delta^* \pi^*_{12} \left( E \otimes \pi^*_1
\left(-\right)\right)\right). 
\end{align}
\end{footnotesize}
By Lemma \ref{lemma-projection-formula-commutes-with-adjunction-1}
morphism \eqref{eqn-pi1-pi2-morphism-transformed-1} is further
isomorphic to 
\begin{align}
\label{eqn-pi1-pi2-morphism-transformed-3}
\pi_{1 *}\pi_{23_*}\left(
Q' \otimes \pi^*_{12} \pi^*_1 \left(-\right)
\right)
\xrightarrow{\quad \nu_\Delta \circ \beta_\Delta \quad}
\pi_{1 *} \pi_{23_*} \Delta_* \left(
\Delta^* Q' \otimes \Delta^* \pi^*_{12} \pi^*_1 \left(-\right) \right).
\end{align}
Finally, since $\pi_1 \circ \pi_{23} = \tilde{\pi}_{2} \circ \pi_{13}$
and $\pi_{1} \circ \pi_{12} = \tilde{\pi}_{1} \circ \pi_{13}$,
see diagram \eqref{eqn-big-projection-tree}, the corresponding
pseudofunctoriality isomorphisms imply
that \eqref{eqn-pi1-pi2-morphism-transformed-3} is isomorphic to 
\begin{align} \label{eqn-pi_2-pi_13-derived-restriction}
\tilde{\pi}_{2 *} \pi_{13 *} 
\left( Q' \otimes \pi^*_{13} \tilde{\pi}^*_1 (-) \right) 
\xrightarrow{\nu_\Delta \circ \beta_\Delta}
\tilde{\pi}_{2 *} \pi_{13 *} \Delta_* 
\left( \Delta^* Q' \otimes \Delta^* \pi^*_{13} \tilde{\pi}^*_1 (-) \right)   
\end{align}

We proceed to morphism \eqref{eqn-pi1-trace-morphism} which is 
induced by 
the adjunction counit $\pi^!_1\mathcal{O}_{X_1} \otimes E^\vee \otimes E
\rightarrow \pi^!_1\mathcal{O}_{X_1}$. By its definition 
morphism \eqref{eqn-Delta-Q'-pi^!_1-O_X1-morphism}
is isomorphic to this adjunction counit, and 
so \eqref{eqn-pi1-trace-morphism} is isomorphic to 
\begin{align}
\label{eqn-transforming-trace-map-1}
\pi_{1 *} \left( \Delta^* Q'\otimes \pi^*_1 \left(-\right) \right)  
\xrightarrow{\quad \eqref{eqn-Delta-Q'-pi^!_1-O_X1-morphism} \quad}
\pi_{1 *} \left(\pi^!_1(\mathcal{O}_{X_1}) \otimes \pi^*_1 \left(-\right) \right)  
\end{align}
As $\tilde{\pi}_2 \circ \Delta = \tilde{\pi}_1 \circ \Delta =
\id$ by pseudofunctoriality 
\eqref{eqn-transforming-trace-map-1} is isomorphic to 
\begin{align} 
\label{eqn-pi_2-Delta-pi_1-trace-map}
\tilde{\pi}_{2 *} \Delta_*
\pi_{1 *} \left( \Delta^* Q'\otimes  
\pi^*_1 \Delta^* \tilde{\pi}^*_{1}  \left(-\right) \right)  
\xrightarrow{\eqref{eqn-Delta-Q'-pi^!_1-O_X1-morphism}}
\tilde{\pi}_{2 *} \Delta_*
\pi_{1 *} \left( \pi^!_1\mathcal{O}_{X_1} \otimes  
\pi^*_1 \Delta^* \tilde{\pi}^*_{1} \left(-\right) \right).
\end{align}

Finally, the same pseudofunctoriality isomorphisms imply that 
\eqref{eqn-pi1-trace-adjunction} is isomorphic to 
\begin{align} 
\label{eqn-pi_2-Delta-pi_1-adjunction}
\tilde{\pi}_{2 *} \Delta_*
\pi_{1 *} \left( \pi^!_1\mathcal{O}_{X_1}\otimes  
\pi^*_1 \Delta_* \tilde{\pi}^*_{1} \left(-\right) \right)  
\xrightarrow{\epsilon_{\pi_1} \circ \chi_{\pi_1}}
\tilde{\pi}_{2 *} \Delta_*
\Delta^* \tilde{\pi}^*_{1}  \left(-\right). 
\end{align}

We have now shown that \eqref{eqn-pi1-pi2-morphism}, 
\eqref{eqn-pi1-trace-morphism} and 
\eqref{eqn-pi1-trace-adjunction} are isomorphic to 
\eqref{eqn-pi_2-pi_13-derived-restriction}, 
\eqref{eqn-pi_2-Delta-pi_1-trace-map} and 
\eqref{eqn-pi_2-Delta-pi_1-adjunction}, respectively. 
Next, we compute the connecting isomorphisms. The isomorphism
connecting \eqref{eqn-pi_2-pi_13-derived-restriction} to
\eqref{eqn-pi_2-Delta-pi_1-trace-map} works out 
to be the pseudofunctoriality isomorphism  
\begin{align}
\label{eqn-pi_2-pi_13-Delta-to-Delta-pi_1-relabelling}
\tilde{\pi}_{2 *} \pi_{13 *} \Delta_* 
\left( \Delta^* Q' \otimes \Delta^* \pi^*_{13} \tilde{\pi}^*_1 (-) \right)   
\xrightarrow{\eta_{\Delta,\pi_1} \circ \eta^{-1}_{\pi_{13},\Delta}}
\tilde{\pi}_{2 *} \Delta_*
\pi_{1 *} \left( \Delta^* Q'\otimes  
\pi^*_1 \Delta^* \tilde{\pi}^*_{1}  \left(-\right) \right).
\end{align}
and the isomorphism connecting \eqref{eqn-pi_2-Delta-pi_1-trace-map} 
to \eqref{eqn-pi_2-Delta-pi_1-adjunction} works out to be the identity.

We can now conclude that the adjunction counit 
$\Phi^{\text{ladj}}_E \Phi_E \rightarrow \id$, 
being the composition of \eqref{eqn-pi1-pi2-morphism},
\eqref{eqn-pi1-trace-morphism} and
\eqref{eqn-pi1-trace-adjunction},
is isomorphic 
to the composition of \eqref{eqn-pi_2-pi_13-derived-restriction}, 
\eqref{eqn-pi_2-pi_13-Delta-to-Delta-pi_1-relabelling},
\eqref{eqn-pi_2-Delta-pi_1-trace-map} and
\eqref{eqn-pi_2-Delta-pi_1-adjunction}.
The claim of the theorem then follows from the fact that the 
following diagram commutes:
\begin{align} \label{eqn-big-commutative-ladder}
\xymatrix{
\tilde{\pi}_{2 *} \left(\pi_{13 *} Q' \otimes \tilde{\pi}^*_1(-) \right)
\ar[r]^{\sim} \ar[d]_{\eqref{eqn-transform-derived-restriction}} &
\tilde{\pi}_{2 *} \pi_{13 *} \left(Q' \otimes \pi^*_{13} \tilde{\pi}^*_1(-) \right)
\ar[d]^{\eqref{eqn-pi_2-pi_13-derived-restriction}} \\
\tilde{\pi}_{2 *} \left(\pi_{13 *} \Delta_* \Delta^* Q' \otimes \tilde{\pi}^*_1(-) \right)
\ar[r]^{\sim}
\ar[d]_{\eqref{eqn-transform-pi_13-Delta-to-Delta-pi_1-relabelling}} &
\tilde{\pi}_{2 *} \pi_{13 *} \Delta_* \left(\Delta^* Q' \otimes \Delta^* \pi^*_{13} \tilde{\pi}^*_1(-) \right)
\ar[d]^{\eqref{eqn-pi_2-pi_13-Delta-to-Delta-pi_1-relabelling}} \\
\tilde{\pi}_{2 *} \left(\Delta_* \pi_{1 *} \Delta^* Q' \otimes \tilde{\pi}^*_1(-) \right)
\ar[r]^{\sim} \ar[d]_{\eqref{eqn-transform-trace-map}} & 
\tilde{\pi}_{2 *} \Delta_* \pi_{1 *} \left(\Delta^* Q' \otimes \pi^*_{1} \Delta^* \tilde{\pi}^*_1(-) \right)
\ar[d]^{\eqref{eqn-pi_2-Delta-pi_1-trace-map}} \\
\tilde{\pi}_{2 *} \left(\Delta_* \pi_{1 *} \pi^{!}_1\mathcal{O}_{X_1} \otimes \tilde{\pi}^*_1(-) \right)
\ar[r]^{\sim} \ar[d]_{\eqref{eqn-transform-adjunction}} & 
\tilde{\pi}_{2 *} \Delta_* \pi_{1 *} \left(\pi^{!}_1\mathcal{O}_{X_1} \otimes \pi^*_{1} \Delta^* \tilde{\pi}^*_1(-) \right)
\ar[d]^{\eqref{eqn-pi_2-Delta-pi_1-adjunction}} \\
\tilde{\pi}_{2 *} \left(\Delta_* \mathcal{O}_{X_1} \otimes \tilde{\pi}^*_1(-) \right)
\ar[r]^{\sim} &
\tilde{\pi}_{2 *} \Delta_* \Delta^* \tilde{\pi}^*_1(-)
} 
\end{align}
where the horizontal isomorphisms are all due to the projection formula.
To see that diagram \eqref{eqn-big-commutative-ladder} indeed commutes, 
observe that its topmost square commutes by Lemma
\ref{lemma-projection-formula-commutes-with-adjunction-1}, the middle 
two commute by functoriality and the lowermost square commutes by Lemma 
\ref{lemma-projection-formula-commutes-with-adjunction-2}.
\end{proof}

An identical proof yields an analogous result 
for the right adjunction counit:
\begin{theorem}\label{theorem-right-adjunction-counit-morphism}
Let $X_1$ and $X_2$ be two separable schemes of finite type over
$k$ with $X_2$ proper. Let $E$ be a perfect object of 
$D(X_1 \times X_2)$ and $\Psi_E$ be a Fourier-Mukai transform 
from $D_{qc}(X_2)$ to $D_{qc}(X_1)$ defined by $E$. 

The adjunction counit $\gamma'_E\colon \Psi_E \Psi^{\text{radj}}_E \rightarrow
\id$ is isomorphic to the morphism of
Fourier-Mukai transforms $D_{qc}(X_1) \rightarrow D_{qc}(X_1)$
induced by the following morphism of objects of $D(X_1 \times X_1)$:
\begin{align} \label{eqn-derived-restriction-morphism-radj-version}
\tilde{Q} = \pi_{13 *}\left(\pi_{12}^* E^\vee \otimes \pi_{23}^* E
\otimes \pi_{12}^* \pi^!_1(\mathcal{O}_{X_1})\right) 
\xrightarrow{\pi_{13 *} \beta_{\Delta}}
\pi_{13 *} \Delta_* \Delta^* 
\left(\pi_{12}^* E^\vee \otimes \pi_{23}^* E \otimes \pi_{12}^*
\pi^!_1(\mathcal{O}_{X_1})\right) \\
\label{eqn-two-squares-isomorphism-radj-version}
\pi_{13 *} \Delta_* \Delta^* 
\left(\pi_{12}^* E^\vee \otimes \pi_{23}^* E \otimes \pi_{12}^*
\pi^!_1(\mathcal{O}_{X_1})\right)
\quad \simeq \quad
\Delta_* \pi_{1 *}  
\left(E \otimes E^\vee \otimes \pi^!_1(\mathcal{O}_{X_1})\right)
\\ 
\label{eqn-E-E-dual-trace-morphism-radj-version}
\Delta_* \pi_{1 *}  
\left(E \otimes E^\vee \otimes \pi^!_1(\mathcal{O}_{X_1})\right) 
\xrightarrow{\Delta_* \pi_{1 *} \ev_E}
\Delta_* \pi_{1 *} \left(\pi^!_1(\mathcal{O}_{X_1})\right)  
\\
\label{eqn-rpi_1-trace-morphism-radj-version}
\Delta_* \pi_{1 *} \pi^!_1(\mathcal{O}_{X_1})
\xrightarrow{\Delta_* \epsilon_{\pi_1}}
\Delta_* \mathcal{O}_{X_1}. 
\end{align}
where \eqref{eqn-two-squares-isomorphism-radj-version} is composed of 
isomorphism $\nu_\Delta\colon \Delta^*\left(- \otimes - \right) 
\xrightarrow{\sim}
\Delta^*\left( - \right) \otimes \Delta^* \left( - \right)$ and 
of pseudofunctoriality 
isomorphisms corresponding to the identities 
$\pi_{13} \circ \Delta = \Delta \circ \pi_1$  and $\pi_{12}
\circ \Delta = \pi_{23} \circ \Delta = \id$.
\end{theorem}

\subsection{Non-compact case}
\label{section-non-compact-case}

In practice, one often has to deal with cases when neither $X_1$ 
nor $X_2$ are proper. A common way to deal with such situations is to
restrict to the full subcategories of objects with proper support. However, 
with a bit of care it is still possible to work in full generality.

So let $X_1$ and $X_2$ be any two separable schemes of finite type over 
$k$, not necessarily proper, and let $E$ be a perfect object 
of $D(X_1 \times X_2)$. We want to write down the left adjoint 
$\Phi^{ladj}_E$ of $\Phi_E = \pi_{2*} \left(E \otimes \pi_1^*(-)\right)$, 
but since $\pi_1$ is not necessarily a proper morphism, 
the left adjoint to $\pi_1^*$ does not necessarily exist.

To construct $\Phi^{ladj}_E$, we compactify $X_2$ - that is, 
we choose an open immersion $j\colon X_2 \hookrightarrow \bar{X}_2$ with
$\bar{X}_2$ proper over $k$, cf. \cite{Nagata-ImbeddingOfAnAbstractVarietyInACompleteVariety}, or 
\cite{Vojta-NagatasEmbeddingTheorem} for a more modern exposition. 
We shall abuse the notation by using $j$ to also denote 
immersions $X_1 \times X_2 \rightarrow
X_1 \times \bar{X}_2$ and $X_1 \times X_2 \times X_1 \rightarrow X_1
\times \bar{X}_2 \times X_1$ where it causes no confusion. For any such 
compactified product space we shall denote by $\bar{\pi}_{i}$ and
$\bar{\pi}_{ij}$ projections onto corresponding factors. Also, write 
$\bar{E}$ for $j_* E$. 

We have the following commutative diagram:
\begin{align} \label{eqn-compactification-diagram}
\vcenter{\xymatrix{
& 
X_1 \times \bar{X}_2 
\ar[dr]_{\bar{\pi}_{2}} \ar[dl]_{\bar{\pi}_1} &
X_1 \times X_2 \ar[dr]+/u 2ex/^{\pi_{2}}
\ar@{_{(}->}[l]_{j}
& \\
X_1 &
&
\bar{X}_2 &
X_2 \ar@{_{(}->}[l]^{j}.
}}
\end{align}

\begin{lemma}  
\label{lemma-Phi-E-and-Phi-E-bar}
Let $E \in D_{perf}(X_1 \times X_2)$. 
There is an isomorphism of functors 
$D_{qc}(X_1) \rightarrow D_{qc}(\bar{X}_2)$
\begin{align}
\label{eqn-Phi-bar-E-equals-iota^*-Phi-E}
\Phi_{\bar{E}} \xrightarrow{\sim} j_* \Phi_{E}.
\end{align}
Its left adjoint with respect to $j_*$ is an isomorphism 
of functors $D_{qc}(X_1) \rightarrow D_{qc}(X_2)$
\begin{align}
\label{eqn-iota^*-Phi-bar-E-equals-Phi-E}
j^* \Phi_{\bar{E}} \xrightarrow{\sim} \Phi_{E}.
\end{align}
\end{lemma}
 
\begin{proof}
For the first claim, we set
\eqref{eqn-Phi-bar-E-equals-iota^*-Phi-E} to be 
\begin{align*}
\Phi_{\bar{E}} = 
\bar{\pi}_{2 *} \left( j_* E \otimes \bar{\pi}^*_1(-) \right)
\xrightarrow{\alpha_j} 
\bar{\pi}_{2 *} j_* \left( E \otimes j^* \bar{\pi}^*_1(-) \right)
\xrightarrow{\eta_{j, \pi_2} \circ \eta^{-1}_{\bar{\pi}_2,j} \circ \zeta_{\bar{\pi}_1, j}} 
j_* \pi_{2 *} \left( E \otimes \pi^*_1(-) \right) 
= j_* \Phi_E.
\end{align*}

For the second claim: \eqref{eqn-iota^*-Phi-bar-E-equals-Phi-E}
is the composition of the image of 
$\eqref{eqn-Phi-bar-E-equals-iota^*-Phi-E}$ under
$j^*$ with the adjunction counit $\gamma_j \colon 
j^* j_* \Phi_{E} \rightarrow \Phi_{E}$. 
And $\gamma_j$ is an isomorphism since $j$ is an open immersion 
\cite[\it Prop. 9.4.2\rm]{Grothendieck-EGA-I}. 
\end{proof}

We now need the following key lemma: 
\begin{lemma}
\label{lemma-open-immersion-i^*-and-i_*-equivalences}
Let $X$ be a concentrated scheme and let $U \xrightarrow{j} X$ 
be an open immersion. Let $D_{qc}^{j}(X)$ be the full subcategory 
of $D_{qc}(X)$ formed by the images of the objects of $D_{qc}(U)$ 
under $j_{*}$. Let $D^{j}(X)$ and $D^{j}_{perf}(X)$ be the full
subcategories of $D^{j}_{qc}(X)$ consisting of complexes with bounded and coherent cohomology and of perfect complexes. Then:
\begin{enumerate}
\item \label{item-i^*-and-i_*-mutually-inverse-autoequivalences}
Functors $j_*$ and $j^*$ restrict to mutually inverse
equivalences between $D_{qc}^{j}(X)$ and $D_{qc}(U)$.
\item \label{item-i^*-and-i_*-preserve-tor-and-rhom-A-*}
For any $A \in D_{qc}(X)$ functors 
$A \otimes (-)$ and $\rder\shhomm_X(A,-)$ 
restrict to functors $D_{qc}^{j}(X) \rightarrow D_{qc}^{j}(X)$ and 
are identified by $j^*$ with  
$j^* A \otimes (-)$ and $\rder\shhomm_U(j^* A, -)$. 
\item  \label{item-i^*-and-i_*-base-change}
Let $X' \xrightarrow{f} X$ be a concentrated map and consider 
the following base change diagram: 
\begin{align}
\sigma\colon
\vcenter{
\xymatrix{
U'
\ar[r]^{j'} 
\ar[d]_{g}
&
X'
\ar[d]^{f}
\\
U 
\ar[r]^{j} 
&
X
}
}
\end{align}
The functors $f_*$ and $f^*$ restrict to functors between 
$D_{qc}^{j'}(X')$ and $D_{qc}^{j}(X)$ and are identified 
by the equivalences $j^*$ and $j'^{*}$ 
with functors $g_{*}$ and $g^{*}$.
\item \label{item-i^*-and-i_*-closed-supports} Let $X$ be Noetherian. 
The equivalence $j^*$ identifies $D^{j}(X)$ and $D^{j}_{perf}(X)$
with  $D^{cls}(U)$ and $D^{cls}_{perf}(U)$, the full subcategories 
of $D(U)$ and $D_{perf}(U)$ consisting of objects whose support is closed in $X$.
\item \label{item-i^*-and-i_*-preserve-rhom-*-A}
Let $X$ be Noetherian. For any $A \in D^{+}_{qc}(X)$ functor 
$\rder\shhomm_{X}(-,A)$ restricts to a functor 
$D^{j}(X) \rightarrow D^{j}_{qc}(X)$ and the equivalence $j^{*}$
identifies it with $\rder\shhomm_{U}(-,j^{*}A)$.
\end{enumerate}
\end{lemma}
\begin{proof}
Since $j$ is an open immersion, the adjunction co-unit $j^* j_*
\xrightarrow{\gamma_j} \id$ is an isomorphism. It follows that 
$j_*$ is fully faithful, so its restriction to a functor 
$D_{qc}(U) \rightarrow D_{j}(X)$ is 
tautologically an equivalence. It also follows that
$j^*$ is the inverse equivalence to $j_*$. This settles
claim \eqref{item-i^*-and-i_*-mutually-inverse-autoequivalences}.

For claim \eqref{item-i^*-and-i_*-preserve-tor-and-rhom-A-*}, 
it follows from the projection formula isomorphism 
$$ A \otimes j_*(-) \xrightarrow{\alpha_j}
j_* (j^*A \otimes -) $$ that $A \otimes (-)$ 
restricts to a functor $D^j(X) \rightarrow D_j(X)$
and that this restriction is identified by $j^*$ with 
$$j^*A \otimes (-) \colon D_{qc}(U) \rightarrow D_{qc}(U).$$
The assertion for the functor $\rder\shhomm_X\left(A,-\right)$ 
follows similarly from the sheafified adjunction isomorphism 
$$ j_* \rder\shhomm_U\left(j^*A, -\right)
\xrightarrow{\tau_j} \rder\shhomm_X\left(A, j_*(-)\right). $$

The claim \eqref{item-i^*-and-i_*-base-change} follows in
the same way from the pseudo-functoriality isomorphism
$f_* j'_* \xrightarrow{\eta_{j,g}\circ\eta^{-1}_{f,j'}}
j_* g_*$ and the flat base change isomorphism 
$f^* j_* \xrightarrow{\mu_{\sigma}} j'_* g^*$. 

For claim \eqref{item-i^*-and-i_*-closed-supports}, first note that $j$ is 
an open immersion
of Notherian schemes and thus perfect. Now let $A$ be any object 
of $D^{j}(X)$ and let $B = j^{*} A$ so that $A = j_* B$. 
Since $j$ is perfect, $B$ lies in $D(U)$. We have
$\supp_{U} B = (\supp_{X} A) \cap U$ and we need to check 
that this set is closed in $X$. Since $A \in D(X)$, we know
that $\supp_{X} A$ is closed in $X$ and 
any point $p \in X$ lies in $\supp_{X} A$ if and only if 
$\iota_{p}^{*} A \neq 0$, where $\iota_{p}$ is the inclusion map. 
On the other hand, for any $p \in X \setminus U$ we have 
$\iota_{p}^{*} A = \iota_{p}^{*} j_{*} B = 0$ 
by the base change formula. Hence $\supp_{X}A \subset U$, 
so $\supp_{U} B =\supp_{X} A$ and hence closed in $X$. 
We conclude that $B \in D^{cls}(U)$ as required.
  
Conversely, let $B \in D^{cls}(U)$ and let $A = j_{*} B$. Since $B \in D(U)$ 
we can find a fat enough closed subscheme  $Z \xrightarrow{k} U$ with 
the underlying set $\supp_{U} B$ to ensure that $B \simeq k_{*} C$ for 
some $C \in D(Z)$. Since $\supp_{U} B$ is closed in $X$, the composite 
map $Z \xrightarrow{j \circ k} X$ is a closed immersion. We conclude that
$A = j_{*} B \simeq j_{*} k_{*} C \simeq (j \circ k)_{*} C$ lies in $D(X)$, 
as required. 

We have now shown that $j^{*}$ identifies $D^{j}(X)$ with $D^{cls}(U)$. Finally, any inverse image functor takes perfect complexes to perfect complexes \cite[Cor. 4.19.1]{IllusieGeneralitesSurLesConditionsDeFinitudeDansLesCategoriesDerivees}, therefore $j^*$ takes $D^{j}_{perf}(X)$ to $D^{cls}_{perf}(U)$. Conversely, let $A$ be a perfect object in 
$D^{cls}(U)$, then it is, in particular, of finite $\tor$-dimension. As $j$ is perfect, $j_{*} A$ is also of finite $\tor$-dimension \cite[Cor. 3.7.2]{IllusieConditionsDeFinitudeRelative}. Since we already know that $j_{*}A \in D(X)$, we conclude that $j_{*} A$ is perfect. Thus 
$j^{*}$ identifies $D_{perf}^{j}(X)$ with $D_{perf}^{cls}(U)$. This settles claim
\eqref{item-i^*-and-i_*-closed-supports}.

For claim \eqref{item-i^*-and-i_*-preserve-rhom-*-A}, take any $B \in D^{j}(X)$. Then, as before, we can find a closed immersion $Z \xrightarrow{k} U$ and 
an object $C \in D(Z)$ such that $B = (j \circ k)_{*} C$. We then have 
a functorial isomorphism
$$ \rder\shhomm_{X}\left( (j \circ k)_{*} C, A\right) 
\xrightarrow{\eta_{j, k}\circ \delta_{j \circ k}}
j_{*}k_{*}\rder\shhomm_{Z}\left( C, (j \circ k)^! A\right) $$
which shows that functor $\rder\shhomm_{X}\left(-, A\right)$ 
restricts to a functor $D^{j}(X) \rightarrow D^{j}_{qc}(X)$. Finally, 
this restriction is identified by $j^*$ 
with $\rder\shhomm_{U}\left(-, j^{*}A\right)$
because $j$ is an open immersion and hence 
the natural morphism
$$ j^{*} \rder\shhomm_{X}\left(B, A\right) \rightarrow 
\rder\shhomm_{X}\left(j^{*}B, j^{*}A\right)
$$
is an isomorphism \cite[Lemma 2.1.7]{AvramovIyengarLipman-ReflexivityAndRigidityForComplexesIISchemes}.
\end{proof}

\begin{cor} 
The Fourier-Mukai transform 
$$\Phi_{E}\colon \quad
D_{qc}(X_{1}) \rightarrow D_{qc}(X_{2})$$ 
has a left adjoint $\Phi^{ladj}_{E}$, and 
this adjoint is isomorphic to the Fourier-Mukai transform
$$ \Psi_{E^{\vee}\otimes \pi^{!}_{1}(\mathcal{O}_{X_1})}\colon\quad
D_{qc}(X_2) \rightarrow D_{qc}(X_1).$$  
If $\supp_{X_1 \times X_2} E$ is proper over $X_{1}$ and $X_{2}$, then $\Phi_{E}$ and
$\Phi^{ladj}_{E}$ restrict to functors between $D(X_{1})$ and $D(X_{2})$.
\end{cor}
\begin{proof}
We only prove the first claim, as the assertion about 
the restriction to $D(X_1)$ and $D(X_2)$ is standard. 
By Lemma \ref{lemma-Phi-E-and-Phi-E-bar} functor $\Phi_{\bar{E}}$ 
is isomorphic to $j_* \Phi_{E}$. Hence it restricts 
to a functor $D_{qc}(X_1) \rightarrow D^j_{qc}(\bar{X}_2)$. 
Thus, by the same lemma, $\Phi_E$ is isomorphic to the composition
$$ D_{qc}(X_1) \xrightarrow{\Phi_{\bar{E}}} D^j_{qc}(\bar{X}_2)
\xrightarrow{j^*} D_{qc}(X_2).$$
By Lemma
\ref{lemma-open-immersion-i^*-and-i_*-equivalences}\eqref{item-i^*-and-i_*-mutually-inverse-autoequivalences}
the functor $D^j_{qc}(\bar{X}_2) \xrightarrow{j^*} D_{qc}(X_2)$ is an
equivalence whose inverse is the functor $j_*$. 
Therefore $\Phi_E$ has a left adjoint $\Phi^{ladj}_E$ isomorphic 
to $\Phi^{ladj}_{\bar{E}} j_*$, that is to 
$$\bar{\pi}_{1 *} \left(\bar{E}^\vee \otimes
\bar{\pi}^!_1(\mathcal{O}_{X_1}) \otimes \bar{\pi}_2^*j_*(-)\right).$$
By Lemma \ref{lemma-open-immersion-i^*-and-i_*-equivalences}\eqref{item-i^*-and-i_*-preserve-tor-and-rhom-A-*}-\eqref{item-i^*-and-i_*-preserve-rhom-*-A} 
this is further isomorphic to 
$$\bar{\pi}_{1 *} j_* \left(E^\vee \otimes
j^* \bar{\pi}^!_1(\mathcal{O}_{X_1}) \otimes \pi_2^*(-)\right).$$
Since $\pi_1 = \bar{\pi}_1 \circ j$, 
the claim now follows by the pseudofunctoriality of $(-)_*$ and $(-)^!$.  
\end{proof}

The isomorphism 
$\Phi_{\bar{E}} \xrightarrow{\eqref{eqn-Phi-bar-E-equals-iota^*-Phi-E}} j_* \Phi_{E}$
of functors induces the unique isomorphism 
\begin{align}
\label{eqn-Phi-bar-E-equals-iota^*-Phi-E-ladj-version}
\Phi^{ladj}_{\bar{E}} \xrightarrow{\sim}
\Phi^{ladj}_E j^*
\end{align}
of their left adjoints $D_{qc}(\bar{X}_2) \rightarrow D_{qc}(X_1)$
which makes the diagram 
\begin{align}
\xymatrix{ 
\Phi^{ladj}_{\bar{E}} \Phi_{\bar{E}} 
\ar[d]^\sim_{
\eqref{eqn-Phi-bar-E-equals-iota^*-Phi-E-ladj-version}
\circ 
\eqref{eqn-Phi-bar-E-equals-iota^*-Phi-E}
} 
\ar[drrr]^{\quad \gamma_{\bar{E}}}
& & & \\
\Phi^{ladj}_{E} j^* j_* \Phi_{E} \ar[rr]^{\sim}_{\gamma_j}
& &
\Phi^{ladj}_{E} \Phi_{E} \ar[r]_{\quad\gamma_E} &
\id
}
\end{align}
of functors $D_{qc}(X_1) \rightarrow D_{qc}(X_1)$ commute.
Therefore 
the adjunction co-unit $\Phi^{ladj}_E \Phi_E \xrightarrow{\gamma_E} \id$
is isomorphic to the adjunction co-unit
$\Phi^{ladj}_{\bar{E}} \Phi_{\bar{E}} \xrightarrow{\gamma_{\bar{E}}} \id$. 
The standard Fourier-Mukai kernel of $\Phi^{ladj}_{\bar{E}} \Phi_{\bar{E}}$ is 
$$ \bar{Q} = \bar{\pi}_{13 *}\left(\bar{\pi}_{12}^* \bar{E} \otimes
\bar{\pi}_{23}^* \bar{E}^\vee \otimes \bar{\pi}_{23}^* 
\bar{\pi}^!_1(\mathcal{O}_{X_1}) \right)$$
and Theorem \ref{theorem-left-adjunction-counit-morphism} 
supplies us with the morphism $\bar{Q} \rightarrow \Delta_*
\mathcal{O}_{X_1}$ which induces 
$\Phi^{ladj}_{\bar{E}} \Phi_{\bar{E}} \xrightarrow{\gamma_{\bar{E}}} \id$. 
We obtain:
\begin{prps} 
The adjunction counit
$\gamma_E\colon \Phi^{ladj}_E \Phi_E \rightarrow \id$
is isomorphic to the morphism of Fourier-Mukai transforms
$D_{qc}(X_1) \rightarrow D_{qc}(X_1)$ induced by the morphism 
$\bar{Q} \rightarrow \Delta_* \mathcal{O}_{X_1}$ of 
Theorem \ref{theorem-left-adjunction-counit-morphism}.
\end{prps}

As a non-essential aside, the standard Fourier-Mukai kernel 
of $\Phi^{\text{ladj}}_E \Phi_E$ itself is 
$$ Q = \pi_{13 *}\left(\pi_{12}^* E \otimes \pi_{23}^* E^\vee
\otimes \pi_{23}^* \pi^!_1(\mathcal{O}_{X_1}) \right) $$
The functors $\Phi^{\text{ladj}}_E \Phi_E$ and 
$\Phi^{\text{ladj}}_{\bar{E}} \Phi_{\bar{E}}$ are isomorphic, but it 
does not a priori mean that $Q$ and $\bar{Q}$ are isomorphic. However, 
it is easy to check that they are --- we leave the details as 
an exercise for the reader. 

\section{An alternative description for the pushforward kernels} 
\label{section-the-pushforward-case}

Whenever $E$ is direct image of an object from the derived category
of some subscheme of $X_1 \times X_2$ the decomposition
of the morphism $Q \rightarrow \mathcal{O}_\Delta$ given in 
Theorem \ref{theorem-left-adjunction-counit-morphism} is usually 
very poorly suited for computing cones. We first illustrate this in Section 
\ref{section-the-global-intersection-example}
with an example where $E$ is the structure sheaf of a global complete
intersection subscheme and so everything can be worked out 
explicitly using Koszul-type resolutions. For a general 
closed subscheme of $X_1 \times X_2$ such a resolution does not exist
and a different approach is needed. But with an insight obtained 
from Section \ref{section-the-global-intersection-example}
we set up some general machinery in Sections
\ref{section-a-decomposition-of-the-evaluation-map} and 
\ref{section-kunneth-maps-and-base-change}
which we then apply in Section
\ref{section-adjunction-counit-for-the-pushforward-kernels}
to obtain a better description of the morphism 
$Q \rightarrow \mathcal{O}_\Delta$ for $E$ being a pushforward 
from an arbitrary closed subscheme.

\subsection{The global complete intersection example}
\label{section-the-global-intersection-example}
Let $X_1$ and $X_2$ be a pair of smooth varieties over $k$ with $X_2$ proper. 
Let $\mathcal{N}$ be a vector bundle of rank $d$ on $X_1 \times X_2$ and let 
$s$ be a regular global section of $\mathcal{N}$. 
Let $Z$ be the zero-locus of $s$ in $X_1 \times X_2$, 
it is a closed subscheme of codimension $d$ and 
normal bundle $\mathcal{N}|_Z$. Let $Z \times X_1$ and  $X_1 \times Z$ 
be $\tor$-independent in $X_1 \times X_2 \times X_1$, 
i.e. the derived tensor product 
$\mathcal{O}_{Z \times X_1} \otimes \mathcal{O}_{X_1 \times Z}$
is $\mathcal{O}_{Z'}$ where $Z' = (Z \times X_1) \cap (X_1 \times Z)$. 
We can rewrite the first two morphisms in the decomposition of
Theorem~\ref{theorem-left-adjunction-counit-morphism} for $E = \mathcal{O}_Z$
as the images under 
$\pi_{13 *}\left( - \otimes \pi^*_{23} \pi^!_2 \mathcal{O}_{X_1}\right)$ of
the following morphism in $D_{qc}(X_1 \times X_2 \times X_1)$:
\begin{align} \label{eqn-wasteful-decomposition}
\pi_{12}^* \mathcal{O}_Z \otimes \pi_{23}^* \mathcal{O}_Z^\vee 
\xrightarrow{\beta_\Delta} 
\Delta_* \left(\mathcal{O}_Z \otimes \mathcal{O}_Z^\vee \right)
\xrightarrow{\Delta_* \ev_{\mathcal{O}_Z}}  
\Delta_* \mathcal{O}_{X_1 \times X_2}.
\end{align}
Note that by the flat base change for the twisted inverse image
pseudofunctor (see 
\S\ref{section-adjunctions-and-dualities-for-derived-functors})
the object $\pi_1^! \mathcal{O}_{X_1}$ is just 
the shifted line bundle $\pi_2^* \omega_{X_2}[\dim X_2]$. 

The structure sheaf $\mathcal{O}_Z$ has a global Koszul
resolution on $X_1 \times X_2$
\begin{align} \label{eqn-global-koszul-resolution}
\wedge^d \mathcal{N}^\vee \rightarrow \wedge^{d-1} \mathcal{N}^\vee 
\rightarrow \dots \rightarrow \mathcal{N}^\vee \rightarrow \mathcal{O}_{X_1
\times X_2}
\end{align}
whose differential maps are defined in the usual way by valuations
at $s$. In particular, they all vanish along $Z$. Dualizing the Koszul 
complex, we see immediately that $\mathcal{O}_Z^\vee$ is isomorphic
to $\mathcal{O}_Z \otimes \wedge^d \mathcal{N}[-d]$ in $D(X_1 \times X_2)$. 

We have $\pi_{12}^{-1}(Z) = Z \times X_1$ and 
$\pi_{23}^{-1}(Z) = X_1 \times Z$. So $\pi_{12}^* \mathcal{O}_Z
\simeq O_{Z \times X_1}$, while $\pi_{23}^* \mathcal{O}_Z^\vee
= \mathcal{O}_{X_1 \times Z} \otimes \pi^*_{23} (\wedge^d \mathcal{N})
[-d]$. Thus 
$\pi_{12}^* \mathcal{O}_Z \otimes \pi_{23}^* \mathcal{O}_Z^\vee$, 
the first term in \eqref{eqn-wasteful-decomposition}, equals
$\mathcal{O}_{Z \times X_1} \otimes \mathcal{O}_{X_1 \times Z} \otimes 
\pi^*_{23} \left(\wedge^d \mathcal{N}\right)[-d]$. 
By assumption $Z \times X_1$ and $X_1 \times Z$ are 
$\tor$-independent, and $\pi^*_{23} \wedge^d \mathcal{N}[-d]$ 
is a line bundle, so we conclude that 
the first term in \eqref{eqn-wasteful-decomposition}
equals $(\pi^*_{23} \wedge^d \mathcal{N})|_{Z'}[-d]$. 

On the other hand, 
$\Delta_* \left( \mathcal{O}_Z \otimes \mathcal{O}_Z^\vee\right)$, 
the second term in \eqref{eqn-wasteful-decomposition},  
is isomorphic to the image under $\Delta_*$ of 
the restriction of the dual of the complex 
\eqref{eqn-global-koszul-resolution} to $Z$. Since all the 
differentials vanish along $Z$, this equals 
\begin{align} 
\Delta_* \left( \mathcal{O}_Z \xrightarrow{0} 
\mathcal{N}|_Z \xrightarrow{0}
\dots \xrightarrow{0} 
\wedge^d \mathcal{N}|_Z\; \right) 
= \bigoplus_{i = 0}^d \wedge^i \mathcal{N} |_{\Delta(Z)} [-i],
\end{align}
where $\Delta(Z)$ is the image of $Z$ under 
$X_1 \times X_2 \xrightarrow{\Delta} X_1 \times X_2 \times X_1$. 

Thus the decomposition \eqref{eqn-wasteful-decomposition} is not
practical from the point of view of computing cones.
Its first map goes from $(\pi^*_{23} \wedge^d \mathcal{N})|_{Z'}[-d]$,
a single sheaf sitting in the degree $d$, to 
$\bigoplus_{i = 0}^d \wedge^i \mathcal{N} |_{\Delta(Z)} [-i]$,
a huge complex with non-zero cohomologies in all degrees 
from $0$ to $d$. Its second map goes from this huge complex
to $\mathcal{O}_{\Delta(X_1 \times X_2)}$, a single sheaf sitting in 
the degree $0$. We get two huge cones with non-zero cohomologies 
in all degrees from $0$ to $d$ which almost entirely annihilate each 
other when we take the cone of the map between them.

In the rest of this section we prove, in a much more general setting, 
that there exists a more economical decomposition than 
\eqref{eqn-wasteful-decomposition}. Applied to the case at hand, 
our result will tell us that the decomposition 
\eqref{eqn-wasteful-decomposition} filters
through the summand $\wedge^d \mathcal{N}|_{\Delta(Z)}[-d]$
of $\bigoplus_{i = 0}^d \wedge^i \mathcal{N} |_{\Delta(Z)} [-i]$, 
and can be written simply as: 
\begin{align} \label{eqn-gci-alternative-decomposition}
(\pi^*_{23} \wedge^d \mathcal{N})|_{Z'}[-d] 
\xrightarrow{Z' \rightarrow  \Delta(Z)}  
\wedge^d \mathcal{N}|_{\Delta(Z)}[-d]
\simeq \Delta_* \mathcal{O}_Z^\vee
\xrightarrow{ \Delta_* \left(\mathcal{O}_{X_1 \times X_2}  \rightarrow
\mathcal{O}_Z\right)^\vee}
\Delta_* \mathcal{O}_{X_1 \times X_2}. 
\end{align}
The cones of these two maps are small compared to those 
in \eqref{eqn-wasteful-decomposition} and easy to compute. 

\subsection{A decomposition of the evaluation map}
\label{section-a-decomposition-of-the-evaluation-map}
Let $Y \xrightarrow{f} X$ be a map of concentrated schemes.  

\begin{prps} \label{prps-decomposition-of-the-evaluation-map}
For any $E \in D(\mathcal{O}_Y\text{-}\modd)$ 
the following diagram commutes
\begin{align} \label{eqn-decomposition-of-the-evaluation-map}
\vcenter{
\xymatrix{
f_* E \otimes \rder \shhomm\left(f_* E, \mathcal{O}_X\right)
\ar[rddd]_{ev_{f_* E}}  &
f_* E \otimes f_* \rder \shhomm\left(E, f^\times \mathcal{O}_X\right). 
\ar[d]^{\kappa_f} 
\ar[l]_>>>>>>{\id \otimes \delta_f} 
\\ 
& f_* \left(E \otimes \rder \shhomm\left(E, f^\times \mathcal{O}_X\right)\right)
\ar[d]^{ev_E}
\\
& f_* f^\times \mathcal{O}_X \ar[d]^{\epsilon_f} 
\\
& \mathcal{O}_X 
}
}
\end{align}
\end{prps}
\begin{proof} 
Let us show that
the right adjoint of \eqref{eqn-decomposition-of-the-evaluation-map}
with respect to $f_*E \otimes (-)$ commutes. The result in
\cite[Prop. 3.2.4(ii)]{Lipman-NotesOnDerivedFunctorsAndGrothendieckDuality}
tells what is the right adjoint of 
$f_* E \otimes f_* (-) \xrightarrow{\kappa_f} f_* (E \otimes -)$
with respect to $f_*E \otimes (-)$. It follows immediately 
that the right adjoint of  
$$
f_* E \otimes f_* \rder \shhomm\left(E, - \right)
\xrightarrow{\kappa_f} 
f_* \left(E \otimes \rder \shhomm\left(E, - \right)\right)
\xrightarrow{\evmap_E}
f_* \left(-\right)
$$
with respect to $f_*E \otimes (-)$ is
$$ f_* \rder \shhomm\left(E, - \right)  
\xrightarrow{\gamma_f}
f_* \rder \shhomm\left(f^*f_* E, - \right)  
\xrightarrow{\tau_f}
\rder \shhomm\left(f_* E, f_* - \right). $$
Therefore the right adjoint of the composition 
$\epsilon_f \circ \ev_E \circ \kappa_f$ in 
\eqref{eqn-decomposition-of-the-evaluation-map} is 
$$ f_* \rder \shhomm\left(E, f^\times \mathcal{O}_X \right)  
\xrightarrow{\gamma_f}
f_* \rder \shhomm\left(f^*f_* E, f^\times \mathcal{O}_X \right)  
\xrightarrow{\tau_f}
\rder \shhomm\left(f_* E, f_* f^\times \mathcal{O}_X \right) 
\xrightarrow{\epsilon_f}
\rder \shhomm\left(f_* E, \mathcal{O}_X \right).
$$
By definition this is just the sheafified Grothendieck duality morphism 
$$ f_* \rder \shhomm\left(E, f^\times \mathcal{O}_X \right)  
\xrightarrow{\delta_f}
\rder \shhomm\left(f_* E, \mathcal{O}_X \right).
$$
So is clearly the right adjoint of the composition 
$\ev_{f_* E} \circ \left(\id \otimes \delta_f\right)$ in
\eqref{eqn-decomposition-of-the-evaluation-map}.
The claim follows.
\end{proof}

\subsection{K{\"u}nneth maps and the base change}
\label{section-kunneth-maps-and-base-change}
Let $Y \xrightarrow{f} X$ be a map of concentrated schemes.
Morphism $\kappa_f\colon f_*(-) \otimes f_*(-) \rightarrow f_*\left( - \otimes -\right)$ can be interpreted as the K{\"u}nneth map 
of the commutative square: 
\begin{align}
\vcenter{
\xymatrix{
Y \ar[d]_{\iden} \ar[r]^{\iden} & Y \ar[d]^f \\
Y \ar[r]_f & X
}
}
\end{align}

We recall the basics on the K{\"u}nneth map, cf. 
\cite[\S3.10]{Lipman-NotesOnDerivedFunctorsAndGrothendieckDuality}:
\begin{defn}
Let  
\begin{align}
\sigma\colon\quad 
\vcenter{
\xymatrix{
Z \ar[d]_{g_1} \ar[r]^{g_2} & Y_2 \ar[d]^{f_2} \\
Y_1 \ar[r]_{f_1} & X
}}
\end{align}
be a commutative square of concentrated schemes. 
Setting $h = f_1 \circ g_1 = f_2 \circ g_2$ define the 
\em K{\"u}nneth map \rm to be the bifunctorial morphism 
\begin{align}
\label{eqn-kunneth-map-statement}
\kappa_{\sigma} \colon\; 
f_{1 *} (A_1) \otimes f_{2 *} (A_2) \rightarrow 
h_*\left(g_1^*(A_1) \otimes g_2^*(A_2)\right) \quad\quad A_i \in D(Y_i)
\end{align}
which is the composition 
\begin{align} \label{eqn-kunneth-map-definition}
f_{1 *} (A_1) \otimes f_{2 *} (A_2) \xrightarrow{\beta_h}
& 
h_* h^* \left( f_{1 *} (A_1) \otimes f_{2 *} (A_2) \right) 
\xrightarrow{\nu_h }
h_* \left(h^* f_{1 *} (A_1) \otimes h^*f_{2 *} (A_2) \right) 
\xrightarrow{\zeta^{-1}_{f_1, g_1} \otimes \zeta^{-1}_{f_2, g_2}}
\\ \notag
\xrightarrow{\zeta^{-1}_{f_1, g_1} \otimes \zeta^{-1}_{f_2, g_2}}
&
h_* \left(g_1^* f_1^* f_{1 *} (A_1) \otimes g_2^* f_2^* f_{2 *} (A_2) \right) 
\xrightarrow{\gamma_{f_1} \otimes \gamma_{f_2}} 
h_* \left(g_1^* (A_1) \otimes g_2^* (A_2) \right) 
\end{align}
with $\beta_h$ being the adjunction unit $\iden_X \rightarrow h_*
h^*$ and $\gamma_{f_i}$ being the adjunction counits $f_i^* f_{i *}
\rightarrow \iden_{Y_i}$. 
\end{defn}

A commutative square is called \em K{\"u}nneth-independent \rm 
if its K{\"u}nneth map is a bifunctorial isomorphism. For fiber squares of
concentrated schemes this notion of independence is equivalent to 
several others:
\begin{prps}[\cite{Lipman-NotesOnDerivedFunctorsAndGrothendieckDuality},
Theorem 3.10.3]
\label{prps-notions-of-independency-for-commutative-squares}
Let 
\begin{align}
\sigma\colon\;
\vcenter{
\xymatrix{
Z = Y_1 \times_X Y_2 \ar[d]_{g_1} \ar[r]^>>>>>{g_2} & Y_2 \ar[d]^{f_2} \\
Y_1 \ar[r]_{f_1} & X
}
}
\end{align}
be a fiber square of concentrated schemes, then 
the following are equivalent:
\begin{enumerate}
\item $\sigma$ is independent, i.e. the base change map $
\mu_\sigma\colon f_1^* f_{2 *} \rightarrow  g_{1 *} g_2^* $ 
is a functorial isomorphism. 
\item $\sigma$ is K{\"u}nneth-independent. 
\item $\sigma$ is $\tor$-independent, i.e. for any pair of points
$y_1 \in Y_1$ and $y_2 \in Y_2$ with $f_1(y_1) = f_2(y_2) = x \in
X$ we have
\begin{align}
\tor^i_{\mathcal{O}_{X,x}}\left(\mathcal{O}_{Y_1, y_1}, \mathcal{O}_{Y_2,
y_2}\right) = 0 \text{ for all } i > 0 .
\end{align}
\end{enumerate}
\end{prps}

What we saw in Section \ref{section-the-global-intersection-example}
is a special case of a very general base change statement for 
K{\"u}nneth maps:
\begin{prps}\label{prps-base-change-for-kunneth-maps}
Let 
\begin{align}
\sigma\colon\;
\vcenter{
\xymatrix{
Z \ar[d]_{g_1} \ar[r]^{g_2} & Y_2 \ar[d]^{f_2} \\
Y_1 \ar[r]_{f_1} & X
}
}
\end{align}
be a commutative square of concentrated schemes and
set $h = f_1 \circ g_1 = f_2 \circ g_2$. 
Let $u\colon X' \rightarrow X$ be any morphism and let $\sigma'$
be the fiber product of $\sigma$ with $X'$ over $X$, that is -
the outer square $(Z', Y'_1, Y'_2, X')$ in the commutative diagram
\begin{align}
\vcenter{
\xymatrix{
Z' \ar[ddd]_{g'_1} \ar[rrr]^{g'_2} \ar[dr]^{u} & & &
Y'_2 \ar[ddd]^{f'_2} \ar[dl]_{u} \\
& Z \ar[d]_{g_1} \ar[r]^{g_2} & Y_2 \ar[d]^{f_2} & \\
& Y_1 \ar[r]_{f_1} & X & \\
Y'_1 \ar[rrr]_{f'_1} \ar[ur]_{u} & & & X' \ar[ul]^{u} \\
}
}
\end{align}
where $Y'_i = Y_i \times_{X, f_i,u} X'$, 
$Z' = Z \times_{X,h,u} X' = Z \times_{Y_i,g_i,u} Y'_i$
and the four squares between $\sigma'$ and $\sigma$ are 
the corresponding fiber squares.  Let also 
$h' = f'_1 \circ g'_1 = f'_2 \circ g'_2$. Finally,
to unburden the notation, write 
\begin{itemize}
\item $\eta_{f_1}$ for the pseudofunctoriality isomorphism
$f_{1 *} u_* \xrightarrow{\eta_{u, f'_1} \circ \eta^{-1}_{f_1, u}}
u_* f'_{1 *}$. 
\item $\zeta_{f_1}$ for the pseudofunctoriality isomorphism 
$u^* f_1^* \xrightarrow{\zeta^{-1}_{u, f'_1} \circ \zeta_{f_1, u}} 
f'^*_1 u^*$   
\item $\mu_{f_1}$ for the base change map 
$u^* f_{1 *} \rightarrow f'_{1 *} u^*$ 
of the corresponding fiber square.    
\end{itemize}
and analogously for $f_2, g_1, g_2$ and $h$. 

Then for any objects $A_1 \in D(Y_1)$ and $A_2 \in D(Y_2)$:
\begin{enumerate} 
\item \label{item-the-base-change-for-kunneth-map}
The following diagram commutes in $D(X')$: 
\begin{align} \label{eqn-the-base-change-for-kunneth-map}
\vcenter{
\xymatrix{
u^* \left(f_{1 *}(A_1) \otimes f_{2 *}(A_2)\right)
\ar[r]^{u^* \kappa_\sigma} 
\ar[d]_{\left(\mu_{f_1} \otimes \mu_{f_2}\right) \circ \nu_u} &
u^* h_* \left(g^*_{1}(A_1) \otimes g^*_{2}(A_2)\right)
\ar[d]^{h'_*\left(\left(\zeta_{g1} \otimes \zeta_{g2}\right) \circ \nu_{u}\right) \circ \mu_{h}}
\\
f'_{1 *}(u^* A_1) \otimes f'_{2 *}(u^* A_2)
\ar[r]_{\kappa_{\sigma'}}
 &
h'_* \left(g'^*_1(u^* A_1) \otimes g'^*_2(u^* A_2)\right)
} 
}
\end{align}
\item \label{item-the-adjoint-base-change-for-kunneth-map}
The following diagram commutes in $D(X)$:
\begin{align}
\label{eqn-the-adjoint-base-change-for-kunneth-map}
\vcenter{
\xymatrix{
f_{1 *}(A_1) \otimes f_{2 *}(A_2)
\ar[r]^{\kappa_\sigma} \ar[d]_{\beta_u} & 
h_* \left(g^*_{1}(A_1) \otimes g^*_{2}(A_2)\right)
\ar[d]^{h_* \beta_u }
\\
u_* u^* (f_{1 *}(A_1) \otimes f_{2 *}(A_2))
\ar[d]_{u_* \left( \left( \mu_{f_{1}} \otimes \mu_{f_{2}} \right)
\circ \nu_u \right)} & 
h_* u_* u^* \left(g^*_{1}(A_1) \otimes g^*_{2}(A_2)\right)
\ar[d]^{u_* h'_*\left(\left(\zeta_{g1} \otimes \zeta_{g2}\right) \circ
\nu_{u}\right) \circ \eta_{h}} 
\\
u_* (f'_{1 *}(u^* A_1) \otimes f'_{2 *}(u^* A_2))
\ar[r]_{u_* \kappa_{\sigma'}}
&
u_* h'_* \left(g'^*_1(u^* A_1) \otimes g'^*_2(u^* A_2)\right)
}
}
\end{align}
\end{enumerate}
\end{prps}
\begin{proof}
By definition, the right adjoint of the base change map
$u^* h_* \xrightarrow{\mu_h} h'_* u^*$ with respect to 
$u^*$ is the composition $h_* \xrightarrow{h_* \beta_u}
h_* u_* u^* \xrightarrow{\eta_h} u_* h'_* u^*$. 
It follows that the diagram 
$\eqref{eqn-the-adjoint-base-change-for-kunneth-map}$
is the right adjoint of the diagram 
\eqref{eqn-the-base-change-for-kunneth-map}
with respect to $u_*$, so we only need to prove that
\eqref{eqn-the-adjoint-base-change-for-kunneth-map} commutes.

Let $B \xrightarrow{m} h_* u_* C$ be any morphism between
some $B \in D(X)$ and some $C \in D(Z')$. Let $l$ be the left adjoint 
$u^* h^* B \rightarrow C$ of $m$ with respect to $h_* u_*$. 
By compatibility of the inverse image/direct image adjunction 
with pseudofunctoriality, the left adjoint with respect to $u_* h'_*$
of the composition $B \xrightarrow{m} h_* u_* C
\xrightarrow{\eta_h} u_* h'_* C$ is 
the composition $h'^* u^* B \xrightarrow{\zeta^{-1}_h} u^* h^*B
\xrightarrow{l} C$. Hence the left adjoint with 
respect to $u_* h'_*$ of the upper-right half 
$$  f_{1 *}A_1 \otimes f_{2 *}A_2 \xrightarrow{\kappa_\sigma} 
h_* \left(g^*_{1}A_1 \otimes g^*_{2}A_2\right) 
\xrightarrow{\nu_u \circ h_* \beta_u} 
h_* u_* \left(u^* g_1^* A_1 \otimes u^* g^*_2 A_2\right)
\xrightarrow{\eta_h \circ h_* u_* \left(\zeta_{g_1} \otimes
\zeta_{g_2}\right)} 
u_* h'_* \left(g'^*_1 u^* A_1 \otimes g'^*_2 u^* A_2\right) $$
of $\eqref{eqn-the-adjoint-base-change-for-kunneth-map}$
is the composition 
of $h'^* u^* \left(f_{1 *}A_1 \otimes f_{2 *}A_2 \right)
\xrightarrow{\zeta^{-1}_h} u^* h^* 
\left(f_{1 *}A_1 \otimes f_{2 *}A_2 \right)$ with the left adjoint
of  
$$  f_{1 *}A_1 \otimes f_{2 *}A_2 \xrightarrow{\kappa_\sigma} 
h_* \left(g^*_{1}A_1 \otimes g^*_{2}A_2\right) \xrightarrow{\nu_u
\circ h_* \beta_u} 
h_* u_* \left(u^* g_1^* A_1 \otimes u^* g^*_2 A_2\right)
\xrightarrow{h_* u_* \left(\zeta_{g_1} \otimes \zeta_{g_2}\right) }
h_* u_* \left(g'^*_1 u^* A_1 \otimes g'^*_2 u^* A_2\right) $$
with respect to $h_* u_*$. Making use of the definition 
of $\kappa_\sigma$ in \eqref{eqn-kunneth-map-definition}, this
adjoint works out to be
$$ u^*h^*\left( \bigotimes_i f_i A_i \right) 
\xrightarrow{\nu_u \circ \left(\bigotimes_i u^* \zeta^{-1}_{f_i, g_i}
\right) \circ u^* \nu_h} 
\bigotimes_i u^*g^*_i f^*_i f_{i *} A_i  
\xrightarrow{\bigotimes_i u^*g^*_i \gamma_{f_i}} 
\bigotimes_i u^* g^*_{i} A_i 
\xrightarrow{\bigotimes_i \zeta_{g_i}}
\bigotimes_i g'^*_i u^* A_i $$
Composing with $h'^* u^* \left(f_{1 *}A_1 \otimes f_{2 *}A_2 \right)
\xrightarrow{\eta^{-1}_h} u^* h^* \left(f_{1 *}A_1 \otimes f_{2 *}A_2 \right)$
and simplifying we see that the left adjoint of 
$\eqref{eqn-the-adjoint-base-change-for-kunneth-map}$ with respect
to $u_* h'_*$ is 
$$ h'^* u^* \left(
\bigotimes_i f_{i *} A_i \right) 
\xrightarrow{
\left(\bigotimes_i \zeta^{-1}_{f'_i, g'_i}
\right)  \circ \nu_{h'} \circ h'^* \nu_u } 
\bigotimes_i g'^*_i f'^*_i u^* f_{i *} A_i  
\xrightarrow{\bigotimes_i g'^*_i \zeta^{-1}_{f_i}} 
\bigotimes_i g'^*_i u^* f^*_i f_{i *} A_i  
\xrightarrow{\bigotimes_i g'^*_i u^* \gamma_{f_i}} 
\bigotimes_i g'^*_i u^* A_i. $$

Similarly, the left adjoint of the lower-left 
half $$ f_{1 *}A_1 \otimes f_{2 *}A_2
\xrightarrow{\beta_u}
u_* u^* \left(f_{1 *} A_1 \otimes f_{2 *} A_2\right)
\xrightarrow{u_* \left(\bigotimes_i \mu_{f_i} \right) \circ u_* \nu_u} 
u_* \left(f'_{1 *} u^* A_1 \otimes f'_{2 *} u^* A_2\right)
\xrightarrow{u_* \kappa_{\sigma'}}
u_* h'_* \left(g'^*_{1} u^* A_1 \otimes g'^*_{2} u^* A_2\right) $$
of \eqref{eqn-the-adjoint-base-change-for-kunneth-map}
with respect to $u_* h'_*$ works out as 
$$ h'^* u^* \left(\bigotimes_i f_{i *} A_i \right)
\xrightarrow{ 
\left(\bigotimes_i \zeta^{-1}_{f'_i, g'_i}
\right)  \circ \nu_{h'} \circ h'^* \nu_u 
}
\bigotimes_i g'^*_i f'^*_i u^* f_{i *}A_i 
\xrightarrow{\bigotimes_i g'^*_i f'^*_i \mu_{f_i}} 
\bigotimes_i g'^*_i f'^*_i f'_{i *} u^* A_i 
\xrightarrow{\otimes_i g'^*_i \gamma_{f'_i}}
\bigotimes_i g'^*_{i} u^* A_i$$

It therefore suffices to show that the following diagram commutes
for $i = 1, 2$ and for all $A_i \in D(Y_i)$
\begin{align} \label{eqn-the-base-change-for-kunneth-map-simplified}
\vcenter{
\xymatrix{
f'^*_i u^* f_{i *} A_i \ar[r]^{g'^*_i \zeta^{-1}_{f_i}} 
\ar[d]_{f'^*_i \mu_{f_i}} &
u^* f^*_i f_{i *} A_i \ar[d]^{u^* \gamma_{f_i}} \\
f'^*_i f'_{i *} u^* A_i \ar[r]_{\gamma_{f'_i}} & u^* A_i.
} 
}
\end{align}
By definition of $\mu_{f_i}$ in \eqref{eqn-definition-of-base-change-morphism} 
the right adjoint with respect to $f'^*_i$ 
of $f'^*_i u^* f_{i *} \xrightarrow{\zeta^{-1}_{f_i}} 
u^* f^*_i f_{i *} \xrightarrow{u^* \gamma_{f_i}} u^*$
is precisely $u^* f_{i *} \xrightarrow{\mu_{f_i}}  f'_{i *} u^*$. 
So the right adjoint with respect to $f'^*_i$ of 
\eqref{eqn-the-base-change-for-kunneth-map-simplified}
is the diagram
\begin{align*} 
\xymatrix{
u^* f_{i *} A_i 
\ar[rd]^{\mu_{f_i}}
\ar[d]_{\mu_{f_i}} &
\\
f'_{i *} u^* A_i \ar[r]_{\iden} & f'_{i *} u^* A_i
} 
\end{align*}
which clearly commutes. 
\end{proof}

\subsection{The adjunction counits for the pushforward Fourier-Mukai
kernels}
\label{section-adjunction-counit-for-the-pushforward-kernels}

We can now apply the generalities of the previous two sections 
to obtain an alternative decomposition to that in Theorem
\ref{theorem-left-adjunction-counit-morphism} of the morphism of 
Fourier-Mukai kernels which induces the canonical adjunction 
morphism $\Phi^{ladj}_E \Phi_E \rightarrow \iden$ in case
where $E$ is a pushforward of an object on some 
$Z \hookrightarrow X_1 \times X_2$.

Let $X_1$ and $X_2$ be a pair of separable schemes of finite 
type over $k$. Let $Z \xrightarrow{\iota_Z} X_1 \times X_2$
be a closed immersion proper over both $X_1$ and $X_2$. Denote
by $\pi_{Z1}$ the composition $Z \xrightarrow{\iota_Z} X_1 \times X_2
\xrightarrow{\pi_1} X_1$. 
Consider the following fiber squares: 
\begin{align*}
\sigma_{12}\colon
\vcenter{
\xymatrix{
Z \times X_1 \ar[r]^<<<<<{\iota_{Z12}} \ar[d]_{\pi_{Z12}}
& X_1 \times X_2 \times X_1 \ar[d]^{\pi_{12}} \\ 
Z \ar[r]_{\iota_Z} & X_1 \times X_2
}
}
\quad \text{ and } \quad
\sigma_{23}\colon
\vcenter{
\xymatrix{
X_1 \times Z \ar[r]^<<<<<{\iota_{Z23}} \ar[d]_{\pi_{Z23}}
& X_1 \times X_2 \times X_1 \ar[d]^{\pi_{23}} \\ 
Z \ar[r]_{\iota_Z} & X_1 \times X_2.
}
}
\end{align*}

Then $Z' = \left(Z \times X_1\right) \cap \left( X_1 \times Z\right) 
\xrightarrow{\iota_{Z'}} X_1 \times X_2 \times X_1$ 
fits into the fiber square
\begin{align}\label{eqn-fiber-square-for-z'}
\sigma\colon 
\vcenter{
\xymatrix{
Z' \ar[dr]^{\iota_{Z'}} \ar[d]_{\iota_{'12}} \ar[r]^{\iota_{'23}} & 
X_1 \times Z 
\ar[d]^{\iota_{Z23}}\\
Z \times X_1 \ar[r]_>>>>{\iota_{Z12}} & X_1 \times X_2 \times X_1.
}
}.
\end{align}
Let $\sigma_\Delta$ denote the square 
obtained from 
\eqref{eqn-fiber-square-for-z'}
by base change 
$X_1 \times X_2 \xrightarrow{\Delta}  X_1 \times X_2 \times X_1$:
\begin{align}\label{eqn-base-change-diagram-for-the-fiber-square-of-z'}
\vcenter{
\xymatrix{
Z \ar[ddd]_{\iden} \ar[rrr]^{\iden} \ar[dr]^{\Delta} & & 
&
Z \ar[ddd]^{\iota_Z} \ar[dl]_{\Delta} \\
& Z'\ar[dr]_{\iota_{Z'}} \ar[d]_{\iota_{'12}} \ar[r]^{\iota_{'23}} & 
X_1 \times Z
\save[]-<-1.5cm,0.7cm>*++[o][F-]\txt{\footnotesize $\sigma_{Z23}$} \restore
\ar[d]^{\iota_{Z23}} & \\
& Z \times X_1 \ar[r]_>>>>{\iota_{Z12}} & X_1 \times X_2 \times X_1 & \\
Z \ar[rrr]_{\iota_Z} \ar[ur]_{\Delta} &
\save[]+<1.5cm,0.5cm>*++[o][F-]\txt{\footnotesize $\sigma_{Z12}$} \restore
& & 
X_1 \times X_2 \ar[ul]^{\Delta} 
}
}
\end{align}
Observe that:
\begin{itemize}
\item Composition
$Z \xrightarrow{\Delta} Z' \xrightarrow{\iota_{Z'}} X_1 \times X_1
\times X_1$ equals $Z \xrightarrow{\pi_{Z1}} X_1
\xrightarrow{\Delta} X_1 \times X_1$. 
\item Compositions
$Z \xrightarrow{\Delta} Z' \xrightarrow{\iota'_{12}} Z \times X_1
\xrightarrow{\pi_{Z12}} Z$
and $Z \xrightarrow{\Delta} Z \times X_1 \xrightarrow{\pi_{Z12}} Z$
are the identity map. 
\item  
Compositions
$Z \xrightarrow{\Delta} Z' \xrightarrow{\iota'_{23}} X_1 \times Z
\xrightarrow{\pi_{Z23}} Z$
and $Z \xrightarrow{\Delta} X_1 \times Z \xrightarrow{\pi_{Z23}} Z$
are the identity map.
\end{itemize}

\begin{theorem} 
\label{theorem-pushfwd-left-adjunction-counit}
Let $E_Z \in D(Z)$ be such that 
$E = \iota_{Z *}(E_Z)$ is perfect in $D(X_1 \times X_2)$.
Let $\Phi_E$ be the Fourier-Mukai transform $D(X_1) \rightarrow
D(X_2)$ with kernel $E$. 
The adjunction counit $\Phi^{ladj}_{E} \Phi_{E} \rightarrow \iden_{X_1}$ 
is isomorphic to the 
morphism of Fourier-Mukai transforms induced by the composition:
\begin{align}
\label{eqn-pushfwd-adj-morphism-left-counit}
\vcenter{
\xymatrix{
Q_Z = \pi_{13 *} \left(\iota_{Z12*} \pi^*_{Z12} E_{Z} \otimes
\iota_{Z23 *} \pi_{Z23}^* 
\rder\shhomm\left(E_{Z}, \pi_{Z1}^!(\mathcal{O}_{X_1})\right)\right)
\ar[d]^{\pi_{13 * } \kappa_\sigma} 
\\
\pi_{13 *} \iota_{Z' *} \left(\iota'^*_{12} \pi^*_{Z12} E_{Z} 
\otimes \iota'^*_{23} \pi_{Z23}^* 
\rder\shhomm\left(E_{Z}, \pi_{Z1}^!(\mathcal{O}_{X_1})\right) \right) 
\ar[d]^{\pi_{13 *} \iota_{Z' *} \beta_{\Delta}}
\\
\pi_{13 *} \iota_{Z' *} \Delta_* \Delta^* \left(\iota'^*_{12}
\pi^*_{Z12} E_{Z} \otimes \iota'^*_{23} \pi_{Z23}^* 
\rder\shhomm\left(E_{Z}, \pi_{Z1}^!(\mathcal{O}_{X_1})\right) \right) 
\ar[d]_{\simeq}^{
\Delta_* \pi_{Z1 *}
\left(\left(\zeta_{\pi_{Z12}, \iota'_{12}, \Delta} \otimes 
\zeta_{\pi_{Z23}, \iota'_{23}, \Delta} \right) \circ
\nu_\Delta \right) \circ \eta_{\Delta, \pi_{Z1}} \circ \eta^{-1}_{\pi_{13}, \iota_{Z'}, \Delta} }
\\
\Delta_* \pi_{Z1 *} \left(E_Z \otimes 
\rder\shhomm\left(E_{Z}, \pi_{Z1}^!(\mathcal{O}_{X_1})\right)\right) 
\ar[d]^{\Delta_* \pi_{Z1 *} \left( \ev_{E_Z} \right)}
\\
\Delta_* \pi_{Z1 *} \pi_{Z1}^!(\mathcal{O}_{X_1}) 
\ar[d]^{\Delta_* \epsilon_{\pi_{Z1}}}
\\
\Delta_* \mathcal{O}_{X_1}
}
}
\end{align}
\end{theorem} 
\begin{proof}
Assume first that $X_2$ is proper. 
By Theorem \ref{theorem-left-adjunction-counit-morphism}
the adjunction counit $\Phi^{ladj}_{E} \Phi_{E} \rightarrow \iden_{X_1}$ 
is induced by the morphism of Fourier-Mukai kernels which 
we reproduce here for the convenience of our readers: 
\begin{align}
\label{eqn-pushfwd-global-restriction-to-diagonal}
\pi_{13 *}\left(\pi_{12}^* E \otimes \pi_{23}^*
\left( E^\vee \otimes \pi^!_1(\mathcal{O}_{X_1})\right)\right)
\xrightarrow{\beta_\Delta}
\pi_{13 *} \Delta_* \Delta^* 
\left(\pi_{12}^* E \otimes \pi_{23}^* \left( E^\vee 
\otimes \pi^!_1(\mathcal{O}_{X_1})\right)\right) \\
\label{eqn-pushfwd-rearrange}
\pi_{13 *} \Delta_* \Delta^* 
\left(\pi_{12}^* E \otimes \pi_{23}^* \left( 
E^\vee \otimes \pi^!_1(\mathcal{O}_{X_1})\right)\right)
\simeq 
\Delta_* \pi_{1 *}  
\left(E \otimes E^\vee 
\otimes \pi^!_1(\mathcal{O}_{X_1})\right)
\\ 
\label{eqn-pushfwd-eval-map}
\Delta_* \pi_{1 *}  
\left(E \otimes E^\vee 
\otimes \pi^!_1(\mathcal{O}_{X_1})\right) 
\xrightarrow{\Delta_* \pi_{1 *} \left( \ev_E \otimes \id \right)}
\Delta_* \pi_{1 *} \left(\pi^!_1(\mathcal{O}_{X_1})\right)  
\\
\label{eqn-pushfwd-trace-morphism}
\Delta_* \pi_{1 *} \left(\pi^!_1(\mathcal{O}_{X_1})\right) 
\xrightarrow{\Delta_* \epsilon_{\pi_1}}
\Delta_* \mathcal{O}_{X_1}. 
\end{align}
where the connecting isomorphism
\eqref{eqn-pushfwd-rearrange}
is $ \Delta_* \pi_{1 *} 
\left(
\left(
\zeta_{\pi_{12}, \Delta} 
\otimes 
\zeta_{\pi_{23}, \Delta} 
\right)
\circ 
\nu_\Delta 
\right) 
\circ
\eta_{\Delta,\pi_1}
\circ
\eta^{-1}_{\pi_{13},\Delta}$. 

We have $E = \iota_{Z *} E_Z$ and 
\begin{align}
\label{eqn-E^vee-to-E_Z^vee}
E^\vee \otimes \pi^!_{1} \mathcal{O}_{X_1}
& = 
\left(\iota_{Z *} E_Z\right)^\vee \otimes \pi^!_{1} \mathcal{O}_{X_1}
\xrightarrow{\eqref{eqn-bringing-an-object-into-rhom-bracket}}
\rder \shhomm\left( \iota_{Z *} E_Z, \pi^!_{1} \mathcal{O}_{X_1} \right)
\xrightarrow{\delta^{-1}_{\iota_{Z}}} 
\iota_{Z *} \rder\shhomm\left(E_Z, \pi^!_{Z1}\mathcal{O}_{X_1}\right). 
\end{align}
Using the isomorphisms
$\pi^*_{12} \iota_{Z *} \xrightarrow{\mu_{\sigma^T_{12}}} 
\iota_{Z12 *} \pi^*_{Z12}$ and 
$\pi^*_{23} \iota_{Z *} 
\xrightarrow{\mu_{\sigma^T_{23}}}
\iota_{Z23 *} \pi^*_{Z23} 
$
and functoriality of $\beta_\Delta$, we see that
\eqref{eqn-pushfwd-global-restriction-to-diagonal}
is isomorphic to
\begin{align}
\label{eqn-transform-of-beta_Delta}
\vcenter{
\xymatrix{
\pi_{13 *} \left(\iota_{Z12*} \pi^*_{Z12} E_{Z} \otimes
\iota_{Z23 *} \pi_{Z23}^* 
\rder\shhomm\left(E_Z, \pi^!_{Z1}\mathcal{O}_{X_1}\right) \right)
\ar[d]^{\pi_{13 *} \beta_\Delta} 
\\
\pi_{13 *} \Delta_* \Delta^* \left(\iota_{Z12*} \pi^*_{Z12} E_{Z} \otimes
\iota_{Z23 *} \pi_{Z23}^*\rder\shhomm\left(E_Z, \pi^!_{Z1}\mathcal{O}_{X_1}\right) \right).
}
}
\end{align}

By Prop. \ref{prps-decomposition-of-the-evaluation-map}
it also follows that
\eqref{eqn-pushfwd-eval-map}-\eqref{eqn-pushfwd-trace-morphism}
is isomorphic to the composition 
\begin{align} 
\label{eqn-transform-of-ev_E}
\vcenter{
\xymatrix{
\Delta_* \pi_{1 *}
\left(\iota_{Z *}E_Z \otimes \iota_{Z *} \rder\shhomm\left(E_Z, \pi^!_{Z1}\mathcal{O}_{X_1}\right)\right)
\ar[d]^{\Delta_* \pi_{1 *} \kappa_{\sigma_\Delta}}
\\
\Delta_* \pi_{1 *}\iota_{Z *}
\left(E_Z \otimes 
\rder\shhomm\left(E_Z, \pi^!_{Z1}\mathcal{O}_{X_1}\right)\right)
\ar[d]^{\Delta_* \pi_{1 *}\iota_{Z *}  \ev_{E_Z} }
\\
\Delta_* \pi_{1 *}\iota_{Z *}
\iota_Z^! \pi^!_{1}(\mathcal{O}_{X_1})
\ar[d]^{\Delta_* \pi_{1 *} \epsilon_{\iota_Z}}
\\
\Delta_* \pi_{1 *} \pi^!_{1}(\mathcal{O}_{X_1})
\ar[d]^{\Delta_* \epsilon_{\pi_1} }
\\
\Delta_* (\mathcal{O}_{X_1}).
}
}
\end{align}

The connecting isomorphism from 
\eqref{eqn-transform-of-beta_Delta}
to \eqref{eqn-transform-of-ev_E}  works out to be 
\begin{align}
\label{eqn-transform-connecting-iso}
\vcenter{
\xymatrix{
\pi_{13 *} \Delta_* \Delta^* \left(\iota_{Z12*} \pi^*_{Z12} E_{Z} \otimes
\iota_{Z23 *} \pi_{Z23}^* \rder\shhomm\left(E_Z, \pi^!_{Z1}\mathcal{O}_{X_1}\right) \right)
\ar[d]_{\simeq}^{
\left( \mu_{\sigma_{Z12}} \otimes \mu_{\sigma_{Z23}}\right)\circ \nu_{\Delta} 
}
\\
\pi_{13 *} \Delta_*  \left(\iota_{Z *} \Delta^* \pi^*_{Z12} E_{Z} \otimes
\iota_{Z *} \Delta^* \pi_{Z23}^* \rder\shhomm\left(E_Z, \pi^!_{Z1}\mathcal{O}_{X_1}\right) \right)
\ar[d]_{\simeq}^{
\Delta_* \pi_{1 *} \left( \iota_{Z *} \zeta_{\pi_{Z12},\Delta} \otimes
\iota_{Z *} \zeta_{\pi_{Z23},\Delta} \right)
\circ
\eta_{\Delta,\pi_1}
\circ \eta^{-1}_{\pi_{13}, \Delta} 
}
\\
\Delta_* \pi_{1 *}
\left(\iota_{Z *}E_Z \otimes \iota_{Z *}
\rder\shhomm\left(E_Z, \pi^!_{Z1}\mathcal{O}_{X_1}\right)\right)
}
}
\end{align}

By functoriality the bottom isomorphism 
of \eqref{eqn-transform-connecting-iso}
commutes with the top morphism 
of \eqref{eqn-transform-of-ev_E}, so we conclude that 
\eqref{eqn-pushfwd-global-restriction-to-diagonal}-\eqref{eqn-pushfwd-trace-morphism} 
is isomorphic to the composition of 
\begin{align}
\label{eqn-pushfwd-global-restriction-to-diagonal-final-part1}
\vcenter{
\xymatrix{
\pi_{13 *} \left(\iota_{Z12*} \pi^*_{Z12} E_{Z} \otimes
\iota_{Z23 *} \pi_{Z23}^* \rder\shhomm\left(E_Z, \pi^!_{Z1}\mathcal{O}_{X_1}\right) \right)
\ar[d]^{\pi_{13 *} \beta_\Delta} 
\\
\pi_{13 *} \Delta_* \Delta^* \left(\iota_{Z12*} \pi^*_{Z12} E_{Z} \otimes
\iota_{Z23 *} \pi_{Z23}^* \rder\shhomm\left(E_Z, \pi^!_{Z1}\mathcal{O}_{X_1}\right) \right)
\ar[d]_{\simeq}^{\pi_{13 *} \Delta_* \left(
\left( \mu_{Z12} \otimes \mu_{Z23}\right)\circ \nu_{\Delta} \right)
}
\\
\pi_{13 *} \Delta_*
\left(\iota_{Z *}\Delta^* \pi^*_{Z12} E_Z \otimes \iota_{Z *}\Delta^*
\pi^*_{Z23}\rder\shhomm\left(E_Z, \pi^!_{Z1}\mathcal{O}_{X_1}\right)\right)
\ar[d]^{\pi_{13 *} \Delta_* \kappa_{\sigma_\Delta}}
\\
\pi_{13 *} \Delta_* \iota_{Z *}
\left(\Delta^* \pi^*_{Z12} E_Z \otimes  
\Delta^* \pi^*_{Z23}\rder\shhomm\left(E_Z, \pi^!_{Z1}\mathcal{O}_{X_1}\right)\right)
}
}
\end{align}
with
\begin{align}
\label{eqn-pushfwd-global-restriction-to-diagonal-final-part2}
\vcenter{
\xymatrix{
\pi_{13 *} \Delta_* \iota_{Z *}
\left(\Delta^* \pi^*_{Z12} E_Z \otimes  
\Delta^* \pi^*_{Z23} \rder\shhomm\left(E_Z, \pi^!_{Z1}\mathcal{O}_{X_1}\right)\right)
\ar[d]_{\simeq}^{
\Delta_* \pi_{1 *} \left(\zeta_{\pi_{Z12},\Delta} \otimes
\zeta_{\pi_{Z23},\Delta} \right)
\circ
\eta_{\Delta,\pi_1}
\circ \eta^{-1}_{\pi_{13}, \Delta} }
\\
\Delta_* \pi_{1 *}\iota_{Z *}
\left(E_Z \otimes \rder\shhomm\left(E_Z, \pi^!_{Z1}\mathcal{O}_{X_1}\right)\right)
\ar[d]^{\Delta_* \pi_{1 *}\iota_{Z *} \ev_{E_Z}}
\\
\Delta_* \pi_{1 *}\iota_{Z *}
\iota_Z^! \pi^!_{1}(\mathcal{O}_{X_1})
\ar[d]^{\Delta_* \pi_{1 *} \epsilon_{\iota_Z}}
\\
\Delta_* \pi_{1 *} \pi^!_{1}(\mathcal{O}_{X_1}).
\ar[d]^{\Delta_* \epsilon_{\pi_1}}
\\
\Delta_* \mathcal{O}_{X_1}
}
}
\end{align}

The claim of the theorem follows by applying 
the base change for K{\"u}nneth maps of
Prop.
\ref{prps-base-change-for-kunneth-maps}(\ref{item-the-adjoint-base-change-for-kunneth-map})
to \eqref{eqn-pushfwd-global-restriction-to-diagonal-final-part1}
and noting that as $\pi_{Z1} = \pi_1 \circ \iota_Z$ so 
by compatibility of the $(f_*, f^\times)$ adjunction with
pseudo-functoriality, counits 
$\pi_{1 *} \epsilon_{\iota_Z}$ and $\epsilon_{\pi_{1}}$ 
at the bottom of 
\eqref{eqn-pushfwd-global-restriction-to-diagonal-final-part2}
compose to give $\epsilon_{\pi_{Z1}}$. 

Suppose now $X_2$ is not proper. Then, following Section 
\ref{section-non-compact-case}, we compactify $X_2$ by choosing
an open immersion $j\colon X_2 \rightarrow \bar{X}_2$ with 
$\bar{X}_2$ proper. Similar to the conventions in 
Section \ref{section-non-compact-case}, we use $j$ to
also denote all the compactification maps induced by 
$j\colon X_2 \rightarrow \bar{X}_2$ and we 
put a bar over various objects and morphisms
to denote their compactified versions. E.g. we denote the inclusion 
$Z \xrightarrow{\iota_Z} X_1 \times X_2 \xrightarrow{j} X_1 \times
\bar{X}_2$ by $\bar{\iota}_Z$. By the argument above 
the compactified version of the composition 
\eqref{eqn-pushfwd-adj-morphism-left-counit}
gives a morphism $\bar{Q}_Z \rightarrow \Delta_{\mathcal{O}_{X_1}}$
which induces the compactified adjunction counit 
$\Phi^{ladj}_{\bar{E}} \Phi_{\bar{E}} \rightarrow
\iden_{X_1}$. By the results of Section \ref{section-non-compact-case}
the compactified and the uncompactified adjunction counits are
naturally isomorphic, therefore to prove the claim of the theorem 
it suffices to exhibit an isomorphism 
$\bar{Q}_Z \xrightarrow{\sim} Q_Z$ which composed with the uncompactified
\eqref{eqn-pushfwd-adj-morphism-left-counit}
gives the compactified  
\eqref{eqn-pushfwd-adj-morphism-left-counit}. 

All the morphisms in \eqref{eqn-pushfwd-adj-morphism-left-counit}
except for the first one are independent of the ambient space $X_2$. 
To be more precise, we have 
$\bar{\pi}_{13} \circ \bar{\iota}_{Z'} = 
\bar{\pi}_{13} \circ j \circ \iota_{Z'} = 
\pi_{13} \circ \iota_{Z'}$,
and hence the compactified versions of last four morphisms in 
\eqref{eqn-pushfwd-adj-morphism-left-counit} are 
isomorphic to the uncompactified ones via pseudofunctoriality
isomorphisms. It therefore suffices to find an isomorphism $\bar{Q}_Z
\xrightarrow{\sim} Q_Z$ that would make the following diagram commute:
\begin{align} \label{eqn-desired-isomorphism-diagram}
\vcenter{
\xymatrix{
Q \ar[rrr]^<<<<<<<<<<<<<<{\pi_{13 *} \kappa_\sigma}
& & & 
\pi_{13 *} \iota_{Z' *} \left(\iota'^*_{12} \pi^*_{Z12} E_{Z} 
\otimes \iota'^*_{23} \pi_{Z23}^* \rder\shhomm\left(E_Z, \pi^!_{Z1}\mathcal{O}_{X_1}\right) \right)
\\
\bar{Q} \ar[u]^{\sim}  
\ar[rrr]_<<<<<<<<<<<<<<{\bar{\pi}_{13 *} \kappa_{\bar{\sigma}}}
& & &  
\bar{\pi}_{13 *} \bar{\iota}_{Z' *} \left(\iota'^*_{12} \pi^*_{Z12} E_{Z} 
\otimes \iota'^*_{23} \pi_{Z23}^* \rder\shhomm\left(E_Z, \pi^!_{Z1}\mathcal{O}_{X_1}\right) \right).
\ar[u]^{\simeq}_{ \eta_{\pi_{13}, \iota_{Z'}} \circ 
\eta^{-1}_{\bar{\pi}_{13}, \bar{\iota}_{Z'}} }
}
}
\end{align}
But $\pi_{13 *} \simeq \bar{\pi}_{13 *} j_*$ 
and square $\sigma$ is obtained from square $\bar{\sigma}$ by 
the base change $j\colon X_1 \times X_2 \times X_1
\rightarrow X_1 \times \bar{X}_2 \times X_1$. So 
the desired statement is precisely the base change 
for K{\"u}nneth maps of Prps. \ref{prps-base-change-for-kunneth-maps}.
\end{proof}

We have similarly:

\begin{theorem} 
\label{theorem-pushfwd-right-adjunction-counit}
Under the assumptions of Theorem 
\ref{theorem-pushfwd-left-adjunction-counit}
let $\Psi_E\colon D(X_2) \rightarrow D(X_1)$
be the Fourier-Mukai transform with kernel $E$.  
The adjunction counit $\Psi_E \Psi^{\text{radj}}_E \rightarrow
\id$ is isomorphic to the morphism of Fourier-Mukai transforms
induced by the composition:

\begin{align}
\label{eqn-pushfwd-adj-morphism-right-counit}
\vcenter{
\xymatrix{
Q'_Z = \pi_{13 *} \left(
\iota_{Z12*} \pi^*_{Z12}
\rder\shhomm\left(E_Z, \pi_{Z1}^! \mathcal{O}_{X_1}\right) 
\otimes 
\iota_{Z23 *} \pi_{Z23}^* E_{Z} \right)
\ar[d]^{\pi_{13 * } \kappa_\sigma} 
\\
\pi_{13 *} \iota_{Z' *} \left(\pi^*_{Z12}
\rder\shhomm\left(E_Z, \pi_{Z1}^! \mathcal{O}_{X_1}\right)
\otimes 
\pi_{Z23}^* E_{Z} \right)
\ar[d]^{\pi_{13 *} \iota_{Z' *} \beta_{\Delta}}
\\
\pi_{13 *} \iota_{Z' *} \Delta_* \Delta^* \left(\pi^*_{Z12}
\rder\shhomm\left(E_Z, \pi_{Z1}^! \mathcal{O}_{X_1}\right)
\otimes 
\pi_{Z23}^* E_{Z}^\vee \right)
\ar[d]_{\simeq}^{
\Delta_* \pi_{Z1 *}
\left(\left(\zeta_{\pi_{Z12}, \iota'_{12}, \Delta} \otimes 
\zeta_{\pi_{Z23}, \iota'_{23}, \Delta} \right) \circ
\nu_\Delta \right) \circ \eta_{\Delta, \pi_{Z1}} \circ \eta^{-1}_{\pi_{13}, \iota_{Z'}, \Delta} }
\\
\Delta_* \pi_{Z1 *} \left(
\rder\shhomm\left(E_Z, \pi_{Z1}^! \mathcal{O}_{X_1}\right)
\otimes E_Z \right) 
\ar[d]^{\Delta_* \pi_{Z1 *} \ev_{E_Z}}
\\
\Delta_* \pi_{Z1 *} \pi_{Z1}^!(\mathcal{O}_{X_1}) 
\ar[d]^{\Delta_* \epsilon_{\pi_{Z1}}}
\\
\Delta_* \mathcal{O}_{X_1}.
}
}
\end{align}
\end{theorem}

One of the main advantages of the alternative decompositions
offered by Theorems \ref{theorem-pushfwd-left-adjunction-counit}
and \ref{theorem-pushfwd-right-adjunction-counit} is that most of 
the morphisms in them can become isomorphisms under fairly reasonable
assumptions on $Z$, $X_1$ and $X_2$. We can then write down 
twists of $\Phi_E$ and $\Psi_E$ fairly easily, for example:  

\begin{cor}\label{cor-best-case-scenario}
Let $X_1$ and $X_2$ be separable schemes of finite 
type over a field $k$. Let $Z \xrightarrow{\iota_Z} X_1 \times X_2$
be a regular closed immersion proper over $X_1$ and $X_2$. Suppose 
$\pi_{Z1 *} \mathcal{O}_Z = \mathcal{O}_{X_1}$
where $\pi_{Z1}$ is 
the composition $Z \xrightarrow{\iota_Z} X_1 \times X_2
\xrightarrow{\pi_1} X_1$. 
Suppose also that $Z \times X_1$ and $X_1 \times Z$
are $\tor$-independent inside $X_1 \times X_2 \times X_1$
and denote by $Z'$ their intersection. 
Denote by $\iota_{Z'}$ the inclusion 
$Z' \hookrightarrow X_1 \times X_2 \times X_1$. 

Then the Fourier-Mukai kernel of 
the dual co-twist of $\Phi_{\mathcal{O}_Z}\colon D(X_1)
\rightarrow D(X_2)$ is 
$\pi_{13 *} \iota_{Z' *}\left(\mathcal{L} \otimes
\mathcal{I}_{\Delta'}[1]\right)$ where $\mathcal{I}_{\Delta'}$ 
is the ideal sheaf of the diagonal $Z$ in $Z'$
and $\mathcal{L}$ is the pullback of $\pi^!_{Z1}(\mathcal{O}_{X_1})$ 
via $X_1 \times Z$ to $Z'$. 
\end{cor}
\begin{proof}
The Fourier-Mukai kernel of the dual co-twist of $\Phi_E$ 
is the cone of the morphism of kernels underlying 
$\Phi^{ladj}_E \Phi_E \rightarrow \id$.
Applying Theorem \ref{theorem-pushfwd-left-adjunction-counit}, 
we note that under the assumptions of 
this corollary, all the morphisms in 
\eqref{eqn-pushfwd-adj-morphism-left-counit} become isomorphisms 
with the exception of 
\begin{align*}
\xymatrix{
\pi_{13 *} \iota_{Z' *} \left(\iota'^*_{12} \pi^*_{Z12} E_{Z} 
\otimes \iota'^*_{23} \pi_{Z23}^* \rder\shhomm\left(E_Z, \pi_{Z1}^! \mathcal{O}_{X_1}\right) \right) 
\ar[d]^{\pi_{13 *} \iota_{Z' *} \beta_{\Delta}}
\\
\pi_{13 *} \iota_{Z' *} \Delta_* \Delta^* \left(\iota'^*_{12}
\pi^*_{Z12} E_{Z} \otimes \iota'^*_{23} \pi_{Z23}^* \rder\shhomm\left(E_Z, \pi_{Z1}^! \mathcal{O}_{X_1}\right) \right).
}
\end{align*}
Since $E_Z = \mathcal{O}_Z$ the above simplifies to
the direct image under $\pi_{13 *} \iota_{Z' *}$ of
$$ 
\iota'^*_{23} \pi_{Z23}^* \pi_{Z1}^!(\mathcal{O}_{X_1}) 
\xrightarrow{\beta_{\Delta}}
\Delta_* \Delta^* \left(
\iota'^*_{23} \pi_{Z23}^* \pi_{Z1}^!(\mathcal{O}_{X_1}) \right).
$$
Write $\mathcal{L}$ for 
$\iota'^*_{23} \pi_{Z23}^* \pi_{Z1}^!(\mathcal{O}_{X_1})$. 
By Lemma \ref{lemma-projection-formula-commutes-with-adjunction-1} 
(with $f = \id$) the morphism
$\mathcal{L} \xrightarrow{\beta_\Delta} \Delta_* \Delta^* \mathcal{L}$
is isomorphic to 
$\mathcal{L} \otimes \left( \mathcal{O}_{Z'} \rightarrow \Delta_*
\Delta^* \mathcal{O}_{Z'} \right)$.  Since
$\mathcal{O}_{Z'} \xrightarrow{\beta_\Delta} 
\Delta_* \Delta^* \mathcal{O}_{Z'}$
is just the sheaf restriction $\mathcal{O}_{Z'} \rightarrow
\Delta_* \mathcal{O}_Z$, its cone is $\mathcal{I}_{\Delta'}[1]$
and the claim follows. 
\end{proof}

\section{An example} 
\label{section-an-example}

Let us give a concrete example of using the results of section 
\ref{section-the-pushforward-case}. For this example we choose
the naive derived category transform induced by the Mukai 
flop. This transform is not an equivalence - it was proved by 
Namikawa in \cite{Namikawa-MukaiFlopsAndDerivedCategories} 
by direct comparison of $\homm$ spaces. 
Below we use Cor. \ref{cor-best-case-scenario} to  
compute the kernel which defines its dual co-twist as 
the Fourier-Mukai transform. We stress that the value of this section 
lies not in the answer itself, but in demonstrating how the methods
of the paper apply to obtain it. However, the reader may observe
that the kernel we obtain is a line bundle supported on the
zero-section of the product. We shall demonstrate 
in \cite{AnnoLogvinenko-OnBraidingCriteriaForSphericalTwistsByFlatFibrations} 
that this is the reason for the braiding which occurs between natural 
spherical twists in the derived categories of the cotangent 
bundles of complete flag varieties (see
\cite{KhovanovThomas-BraidCobordismsTriangulatedCategoriesAndFlagVarieties},
\S 4).  

Let $V$ be a $3$-dimensional vector space and let 
$X_1$ be the scheme $T^*{\mathbb P}(V)$, that is - the total 
space of the cotangent bundle of ${\mathbb P}(V)$. 
Similarly, let $X_2$ be the scheme $T^*{\mathbb P}(V^\vee)$.
These schemes admit the following description:
\begin{align*}
X_1 = \left\{
\xymatrix@C=0.1em{
0 &\subset & U_1 \ar_{\alpha}@/_1em/[ll] &\subset &V
\ar_{\alpha}@/_1em/[ll]} \right\} \; := 
\left\{ 
\begin{array}{c|c}
U_1 \subset V , \; \alpha \in \eend(V) & 
\dim U_1 = 1, \alpha(V) \subseteq U_1,\; \alpha(U_1) = 0 
\end{array}
\right\}
\end{align*}

\begin{align*}
X_2 = \left\{
\xymatrix@C=0.1em{
0 &\subset & U_2 \ar_{\alpha}@/_1em/[ll] &\subset &V
\ar_{\alpha}@/_1em/[ll]} \right\} \; := 
\left\{ 
\begin{array}{c|c}
U_2 \subset V , \; \alpha \in \eend(V) & 
\dim U_2 = 2, \alpha(V) \subseteq U_2,\; \alpha(U_2) = 0 
\end{array}
\right\}
\end{align*}

We also have a variety
\begin{align*}
Z = \left\{
\xymatrix@C=0.1em{
0 &\subset & U_1 &\subset
&U_2 \ar_{\alpha}@/_1em/[llll]
\subset &V
\ar_{\alpha}@/_1em/[lll]} \right\} 
\end{align*}
with natural ``forgetful'' maps $\phi_k\colon Z
\rightarrow X_k$ which forget one of the subspaces. 
Each map $\phi_k$ is isomorphic to the blow-up of 
the zero section carved out by $\alpha = 0$ in $X_k$. Both blowups
have the same exceptional divisor $F \subset Z$ which is carved out
by $\alpha = 0$:
$$ F = \left\{ 0 \subset U_1 \subset U_2 \subset V \right\}. $$
The resulting birational transformation $X_1 \dashrightarrow X_2$ 
which transforms the zero-section $\mathbb{P}(V) \hookrightarrow X_1$ 
into the zero-section 
$\mathbb{P}(V^\vee) \hookrightarrow X_2$ 
is a local model of a four-dimensional Mukai flop. 
Note that maps $\phi_k$ are proper and, since 
each map $\phi_k$ is a blowup of $X_k$, we have
${\phi_k}_*{\mathcal O}_Z = {\mathcal O}_{X_k}$. 

Let $\Phi$ be the functor $\phi_{2 *} \phi^*_1$ from $D(X_1)$ to
$D(X_2)$ and let us compute its dual co-twist.  
The functor $\Phi$ is a Fourier-Mukai transform with the kernel 
${\iota_Z}_*{\mathcal O}_Z$, where
$\iota_Z=\phi_1\times\phi_2: Z \to X_1\times X_2$. We have:
\begin{align*}
X_1 \times X_2 \times X_1 = \left\{
\xymatrix@C=0.1em{
0 &\subset & U_1, U_2, U'_1 \ar_{\alpha_{1},\alpha_2, \alpha'_1}@/_1em/[ll] &\subset &V
\ar_{\alpha_1,\alpha_2,\alpha'_1}@/_1em/[ll]
}
\right\}
\end{align*}
\begin{align*}
Z \times X_1 = \left\{
\begin{array}{c}
\xymatrix@C=0.1em{
0 &\subset &U_1 &\subset
&U_2 \ar_{\alpha_1 = \alpha_2}@/_1em/[llll]
&\subset &V
\ar_{\alpha_1 = \alpha_2}@/_1em/[llll]}, \;
\xymatrix@C=0.1em{
0 &\subset & U'_1 \ar_{\alpha'_1}@/_1em/[ll] &\subset &V
\ar_{\alpha'_1}@/_1em/[ll]}
\end{array}
\right\}
\end{align*}
\begin{align*}
X_1 \times Z = \left\{
\begin{array}{c}
\xymatrix@C=0.1em{
0 &\subset & U_1 \ar_{\alpha_1}@/_1em/[ll] &\subset &V
\ar_{\alpha_1}@/_1em/[ll]}, \;
\xymatrix@C=0.1em{
0 &\subset &U'_1 &\subset
&U_2 \ar_{\alpha_2 = \alpha'_1}@/_1em/[llll]
&\subset &V
\ar_{\alpha_2 = \alpha'_1}@/_1em/[llll]}
\end{array}
\right\}.
\end{align*}
It follows that $Z' = \left(Z \times X_1\right) \cap \left(X_1 \times
Z\right) \subset X_1 \times X_2 \times X_1$ can be described as
\begin{align*}
Z' = \left\{
\begin{array}{c|c}
\xymatrix@C=0.1em{
0 &\subset&U_1,U'_1&\subset
&U_2 \ar_{\alpha}@/_1em/[llll]
&\subset &V
\ar_{\alpha}@/_1em/[llll]} &
\alpha(V) \subseteq U_1 \cap U'_1
\end{array}
\right\}.
\end{align*}

Observe that for any point of $Z'$ we have $U_1 = U'_1$ or 
$\alpha = 0$ (or both). Therefore $Z'$ consists of two 
irreducible components: the diagonal $\Delta Z$
and the zero section 
\begin{align*}
P = \left\{
0 \subset U_1, U'_1 \subset
U_2 \subset V \right\}.
\end{align*}
The intersection $\Delta Z \cap P$ considered as a subvariety of 
$\Delta Z$ is the exceptional divisor $F$ of the blowups 
$Z \xrightarrow{\phi_i} X_i$ described above.
On the other hand, let $P \xrightarrow{\phi_{13}} \mathbb{P}(V) \times
\mathbb{P}(V)$ be the map which forgets the subspace $U_2$. It 
is the blowup of the diagonal of $\mathbb{P}(V) \times \mathbb{P}(V)$ and its
exceptional divisor in $P$ is carved out by $U_1 = U'_1$, i.e. 
it is $F = \Delta Z \cap P$ again. 

By Cor. \ref{cor-best-case-scenario} the dual co-twist 
of $\Phi$ is the Fourier-Mukai
transform $X_1 \rightarrow X_1$ with kernel
$$K = {\pi_{13}}_* {\iota_{Z'}}_* \left(\mathcal{L} \otimes {\mathcal I}_{\Delta} [1] \right) \in D(X_1 \times X_1).$$
Here $\iota_{Z'}$ is 
the inclusion $Z' \hookrightarrow X_1 \times X_2 \times X_1$, 
${\mathcal I}_{\Delta}$ is the ideal sheaf of $\Delta Z$ in $Z'$
and $\mathcal{L}$ is the pullback of $\phi^!_1(\mathcal{O}_{X_1})$ 
to $Z'$ via $X_1 \times Z$. 

Since $Z \xrightarrow{\phi_1} X_1$ is the blow-up of 
the zero-section $\mathbb{P}(V) \hookrightarrow X_1$ whose 
codimension is $2$, we know that $\phi^!_1(\mathcal{O}_{X_1})$ is the 
line bundle $\mathcal{O}_Z(F)$ where $F$ is the exceptional divisor of the 
blow-up. On the other hand, pulling back along the projection 
$$ Z \rightarrow \mathbb{P}(V) \times \mathbb{P}(V^\vee) $$ 
induces an isomorphism 
$$ \picr Z \simeq \picr \mathbb{P}(V) \times \picr \mathbb{P}(V^\vee).$$ 
A simple calculation shows that $\mathcal{O}_Z(F)$ is the pullback of
$\mathcal{O}_{\mathbb{P}(V) \times \mathbb{P}(V^\vee)}(-1,-1)$.
Similarly
$$ \picr Z' \simeq \picr \mathbb{P}(V) \times \picr \mathbb{P}(V^\vee)
\times \picr \mathbb{P}(V)$$
and $\mathcal{L}$, being the pullback to $Z'$ of
$\phi^!_1(\mathcal{O}_{X_1})$
via $X_1 \times Z$, is then the pullback of 
$\mathcal{O}_{\mathbb{P}(V) \times \mathbb{P}(V^\vee) \times
\mathbb{P}(V)}(0,-1,-1)$.

Since $Z'$ has two irreducible 
components $\Delta Z$ and $P$, we have $\mathcal{I}_\Delta \simeq
\iota_{P *} \mathcal{O}_P(- \Delta Z \cap P)$ where $\iota_P$
is the inclusion $P \hookrightarrow Z'$. 
We therefore have $K \simeq \pi_{13 *} \iota_{Z ' *} \iota_{P *} 
\left(\iota^*_P \mathcal{L} \otimes \mathcal{O}_P(-F) [1]\right)$.
A simple computation shows that $\mathcal{O}_P(-F)$ is the pullback of 
$\mathcal{O}_{\mathbb{P}(V) \times \mathbb{P}(V^\vee) \times
\mathbb{P}(V)}(-1,1,-1)$ and  
therefore $\iota^*_P \mathcal{L} \otimes \mathcal{O}_P(-F)$
is the pullback of $\mathcal{O}_{\mathbb{P}(V) \times
\mathbb{P}(V^\vee) \times \mathbb{P}(V)}(-1,0,-2)$.
We conclude that $K \simeq 
\pi_{13 *} \iota_{Z ' *}
\iota_{P *} \phi_{13}^*\left(\mathcal{O}_{\mathbb{P}(V) \times
\mathbb{P}(V)}(-1,-2)[1]\right)$. 

Now observe that the following diagram commutes
\begin{align*}
\xymatrix{
P \ar[r]^{\iota_P} \ar[d]_{\phi_{13}} & Z' \ar[r]^>>>>>{\iota_{Z'}} & 
X_1 \times X_2 \times X_1 \ar[d]_{\pi_{13}} \\  
\mathbb{P}(V) \times \mathbb{P}(V) \ar[rr]_{\iota_{0}} & & X_1 \times X_1
}
\end{align*}
where $\iota_0$ is the zero-section inclusion of $\mathbb{P}(V) \times
\mathbb{P}(V)$ into $X_1 \times X_1$. 
We conclude that
\begin{align*}
K \simeq
\iota_{0 *} \phi_{13 *} \phi_{13}^* 
\left(\mathcal{O}_{\mathbb{P}(V) \times
\mathbb{P}(V)}(-1,-2)[1]\right)
\simeq 
\iota_{0 *} 
\left( \mathcal{O}_{\mathbb{P}(V) \times
\mathbb{P}(V)}(-1,-2)[1]\right).
\end{align*}

\appendix

\section{ The unabridged proof of Theorem
\ref{theorem-left-adjunction-counit-morphism}}
Here we give a complete version of the proof of Theorem 
\ref{theorem-left-adjunction-counit-morphism}. It contains 
explicit computations of all the connecting isomorphisms which 
we left out of the version in the main body of the paper so 
as to emphasise the meaningful part of the proof. 
The version below is for referees and 
others who relish seeing how the monoidal structure 
of the inverse image functor commutes with pseudofunctoriality 
and with the associativity of tensor product.
Lasciate ogne speranza, voi ch'intrate.

\begin{proof}
Set
$$ Q' = \pi^*_{23} \left(\pi^!_1\mathcal{O}_{X_1} \otimes E^\vee
\right) \otimes \pi^*_{12} E $$ so that $Q = \pi_{13 *} Q'$.  
Since $\pi_{12} \circ \Delta = \pi_{23} \circ \Delta = \id$ we 
have a natural isomorphism 
\begin{align}
\label{eqn-appendix-Delta-Q'-iso}
\Delta^* Q' \xrightarrow{\quad \nu_\Delta \quad} 
\Delta^* \pi^*_{23} \left(\pi^!_1\mathcal{O}_{X_1} \otimes E^\vee
\right) \otimes \Delta^* \pi^*_{12} E 
\xrightarrow{\zeta_{\pi_{23},\Delta} \otimes \zeta_{\pi_{12},\Delta}}
\left(\pi^!_1\mathcal{O}_{X_1} \otimes E^\vee \right) \otimes E.
\end{align}
We therefore define a morphism 
\begin{align}
\label{eqn-appendix-Delta-Q'-pi^!_1-O_X1-morphism}
\Delta^* Q' \xrightarrow{\eqref{eqn-appendix-Delta-Q'-iso}}
\left(\pi^!_1\mathcal{O}_{X_1} \otimes E^\vee \right) \otimes E
\xrightarrow{E \otimes \left(E^\vee \otimes (-)\right) \rightarrow \id} 
\pi^!_1\mathcal{O}_{X_1}.
\end{align}

In these terms, the morphism of FM-transforms 
$D(X) \rightarrow D(X)$ induced by
$
Q 
\xrightarrow{  
\eqref{eqn-derived-restriction-morphism}-\eqref{eqn-rpi_1-trace-morphism}
}
\Delta \mathcal{O}_X
$
is: 
\begin{align} \label{eqn-appendix-transform-derived-restriction}
\tilde{\pi}_{2 *} \left( \pi_{13 *} Q' \otimes \tilde{\pi}^*_1(-)\right)  
\xrightarrow{\id \rightarrow \Delta_* \Delta^*} 
\tilde{\pi}_{2 *} \left( \pi_{13 *}\Delta_*\Delta^* Q' \otimes \tilde{\pi}^*_1(-)\right)  
\\ \label{eqn-appendix-transform-pi_13-Delta-to-Delta-pi_1-relabelling}
\tilde{\pi}_{2 *} \left( \pi_{13 *}\Delta_*\Delta^* Q' \otimes \tilde{\pi}^*_1(-)\right)  
\xrightarrow{\eta_{\Delta,\pi_1} \circ \eta^{-1}_{\pi_{13},\Delta}}
\tilde{\pi}_{2 *} \left( \Delta_*\pi_{1 *} \Delta^* Q' \otimes \tilde{\pi}^*_1(-)\right)  
\\ \label{eqn-appendix-transform-trace-map}
\tilde{\pi}_{2 *} \left( \Delta_*\pi_{1 *} \Delta^* Q' \otimes \tilde{\pi}^*_1(-)\right)  
\xrightarrow{\eqref{eqn-appendix-Delta-Q'-pi^!_1-O_X1-morphism}}
\tilde{\pi}_{2 *} \left( \Delta_*\pi_{1 *} \pi^{!}_1 \mathcal{O}_{X_1} \otimes \tilde{\pi}^*_1(-)\right)  
\\ \label{eqn-appendix-transform-adjunction}
\tilde{\pi}_{2 *} \left( \Delta_*\pi_{1 *} \pi^{!}_1 \mathcal{O}_{X_1} \otimes \tilde{\pi}^*_1(-)\right)  
\xrightarrow{\pi_{1 *} \pi^!_{1} \rightarrow \id}
\tilde{\pi}_{2 *} \left( \Delta_* \mathcal{O}_{X_1} \otimes \tilde{\pi}^*_1(-)\right)  
\end{align}
On the other hand, $\Phi_E$ is the composition of functors $\pi^*_{1}$,
$E \otimes (-)$ and $\pi_{2 *}$. Each of these functors has a left
adjoint, these adjoints are $\pi_{1 *}(\pi_1^!\mathcal{O}_{X_1} \otimes -)$, 
$E^\vee \otimes (-)$ and $\pi^*_{2}$, respectively. 
Therefore, the adjunction counit 
$\Phi^{\text{ladj}}_E \Phi_E \rightarrow \id$
is the composition of the three corresponding adjunction counits: 
\begin{align} \label{eqn-appendix-pi1-pi2-morphism}
\pi_{1 *} \left(\pi^!_1\mathcal{O}_{X_1} \otimes \left( E^\vee \otimes
\pi^*_2 \pi_{2 *} \left(E \otimes \pi^*_1 \left(-\right)\right)\right) \right)
\xrightarrow{\pi^*_2 \pi_{2 *} \rightarrow \id}
\pi_{1 *} \left(\pi^!_1\mathcal{O}_{X_1} \otimes \left(E^\vee \otimes
\left( E \otimes \pi^*_1 \left(-\right)\right) \right) \right) \\
\label{eqn-appendix-pi1-trace-morphism}
\pi_{1 *} \left(\pi^!_1\mathcal{O}_{X_1} \otimes \left( E^\vee \otimes \left( E
\otimes \pi^*_1 \left(-\right)\right) \right) \right) 
\xrightarrow{E^\vee \otimes \left( E \otimes (-)\right) \rightarrow \id} 
\pi_{1 *} \left(\pi^!_1\mathcal{O}_{X_1} \otimes \pi^*_1 \left(-\right)\right) \\
\label{eqn-appendix-pi1-trace-adjunction}
\pi_{1 *} \left(\pi^!_1\mathcal{O}_{X_1} \otimes \pi^*_1 \left(-\right)\right) 
\rightarrow
\id
\end{align}
The claim of the theorem is that the composition 
\eqref{eqn-appendix-pi1-pi2-morphism}-\eqref{eqn-appendix-pi1-trace-adjunction} 
is isomorphic to the composition
\eqref{eqn-appendix-transform-derived-restriction}-\eqref{eqn-appendix-transform-adjunction}.

Let us clarify some terminology. We say that 
two morphisms of functors $f \rightarrow g$ and $f' \rightarrow g'$ 
are isomorphic if there exist 
connecting isomorphisms $f \xrightarrow{\sim} f'$ and 
$g \xrightarrow{\sim} g'$ such that the diagram 
\begin{align}
\label{eqn-appendix-f-g-iso-to-f'-g'}
\xymatrix{
f \ar[r] \ar[d]_{\sim} &  
g \ar[d]^{\sim} \\
f' \ar[r] &
g'
}
\end{align}
commutes. Clearly it is an equivalence relation on the set of all
morphisms between all functors between two given categories. In
particular, it is transitive.  

If we further have a morphism of functors $g \rightarrow h$
which is isomorphic to a morphism of functors $g'' \rightarrow h''$
then $f \rightarrow g \rightarrow h$ is isomorphic
to $f' \rightarrow g' \xrightarrow{\sim} g'' \rightarrow h''$, 
where the connecting isomorphism $g' \xrightarrow{\sim} g''$
is the composition of the inverse of the connecting isomorphism
$g \xrightarrow{\sim} g'$ with the connecting isomorphism 
$g \xrightarrow{\sim} g''$.

Our strategy therefore is to 
consecutively replace the morphisms which compose
\eqref{eqn-pi1-pi2-morphism}-\eqref{eqn-pi1-trace-adjunction} 
by isomorphic ones until we obtain 
\eqref{eqn-transform-derived-restriction}-\eqref{eqn-transform-adjunction}.

Observe that the following diagram, whose vertical 
arrows are all isomorphisms, commutes:
\begin{tiny}
\begin{align}
\label{eqn-appendix-transforming-pi2^*-pi_2*-adj-counit-diagram}
\xymatrix{
\pi_{1 *} \left(\pi^!_1\mathcal{O}_{X_1} \otimes \left( E^\vee \otimes
\pi^*_2 \pi_{2 *} \left(E \otimes \pi^*_1 \left(-\right)\right)\right) \right)
\ar[rr]^{\eqref{eqn-appendix-pi1-pi2-morphism}} 
\ar[d]_{\rho^{-1}}
& \quad &
\pi_{1 *} \left(\pi^!_1\mathcal{O}_{X_1} \otimes \left(E^\vee \otimes
\left( E \otimes \pi^*_1 \left(-\right)\right) \right) \right)  
\ar[d]^{\rho^{-1}}
\\
\pi_{1 *} \left(\left(E^\vee \otimes \pi^!_1\mathcal{O}_{X_1}\right) 
\otimes \pi^*_2 \pi_{2 *} \left(E \otimes \pi^*_1 \left(-\right)\right)\right) 
\ar[d]_{\mu}
\ar[rr]^{\id \otimes \left( \pi^*_2 \pi_{2 *} \rightarrow \id \right) }
& \quad &
\pi_{1 *} \left(\left(E^\vee \otimes \pi^!_1\mathcal{O}_{X_1}\right) 
\otimes \left( E \otimes \pi^*_1 \left(-\right)\right)\right)
\ar[d]^{\id \otimes \left(\eta_{\pi_{23}, \Delta} \circ
\zeta^{-1}_{\pi_{12}, \Delta} \right)}
\\
\pi_{1 *} \left(\left(E^\vee \otimes \pi^!_1\mathcal{O}_{X_1}\right) \otimes
\pi_{23*} \pi^*_{12} \left(E \otimes \pi^*_1 \left(-\right)\right)\right) 
\ar[rr]^{\beta_\Delta} 
\ar[d]_{\alpha_{\pi_{23}}}
& \quad &
\pi_{1 *} \left(\left(E^\vee \otimes \pi^!_1\mathcal{O}_{X_1}\right) \otimes 
\pi_{23*} \Delta_{*} \Delta^* \pi^*_{12} \left( E \otimes \pi^*_1 \left(-\right)\right)\right)
\ar[d]^{\nu_\Delta^{-1} \circ \alpha_{\Delta} \circ \alpha_{\pi_{23}}}
\\
\pi_{1 *}\pi_{23_*} \left(\pi^*_{23} \left(E^\vee \otimes \pi^!_1\mathcal{O}_{X_1}\right) \otimes
\pi^*_{12} \left(E \otimes \pi^*_1 \left(-\right)\right)\right) 
\ar[rr]^{\beta_\Delta} 
\ar[d]_{\rho^{-1} \circ \left(\id \otimes \nu_{\pi_{12}}\right)}
& \quad &
\pi_{1 *} \pi_{23_*} \Delta_* \Delta^*
\left(\pi^*_{23} \left(E^\vee \otimes \pi^!_1\mathcal{O}_{X_1}\right) 
\otimes 
\pi^*_{12} \left( E \otimes \pi^*_1 \left(-\right)\right)\right)
\ar[d]^{\rho^{-1} \circ \left(\id \otimes \nu_{\pi_{12}}\right)}
\\
\pi_{1 *}\pi_{23_*}\left(
\left(\pi^*_{23} \left(E^\vee \otimes \pi^!_1\mathcal{O}_{X_1}\right) \otimes
\pi^*_{12} E \right) 
\otimes \pi^*_{12} \pi^*_1 \left(-\right)
\right)
\ar[rr]^{\beta_\Delta} 
\ar[d]_{\id}
& \quad &
\pi_{1 *} \pi_{23_*} \Delta_* \Delta^* \left(
\left(
\pi^*_{23} \left(E^\vee \otimes \pi^!_1\mathcal{O}_{X_1}\right) 
\otimes \pi^*_{12} E \right) 
\otimes \pi^*_{12} \pi^*_1 \left(-\right)
\right)
\ar[d]^{\nu_{\Delta}}
\\
\pi_{1 *}\pi_{23_*}\left(
\left(\pi^*_{23} \left(E^\vee \otimes \pi^!_1\mathcal{O}_{X_1}\right) \otimes
\pi^*_{12} E \right) 
\otimes \pi^*_{12} \pi^*_1 \left(-\right)
\right)
\ar[rr]^{\nu_\Delta \circ \beta_\Delta} 
& \quad &
\pi_{1 *} \pi_{23_*} \Delta_* \left(
\Delta^* \left(
\pi^*_{23} \left(E^\vee \otimes \pi^!_1\mathcal{O}_{X_1}\right) 
\otimes \pi^*_{12} E \right) 
\otimes \Delta^* \pi^*_{12} \pi^*_1 \left(-\right)
\right)
} 
\end{align}
\end{tiny}
The first square in it commutes by functoriality of $\rho^{-1}$, the second 
commutes by Lemma
\ref{lemma-change-of-base-turns-adjunct-into-diag-restrict}, 
the third commutes by Lemma 
\ref{lemma-projection-formula-commutes-with-adjunction-1}, the fourth
commutes by functoriality of $\beta_\Delta$ and the fifth commutes
tautologically. 

We now want to simplify the connecting isomorphism in the right
column of \eqref{eqn-appendix-transforming-pi2^*-pi_2*-adj-counit-diagram}. 
By compatibility of the projection formula with pseudofunctoriality (see diagram
\eqref{eqn-projection-formula-is-compatible-with-pseudofunctoriality})
we have an equality 
$$
\alpha_\Delta \circ \alpha_{\pi_{23}} \circ 
\left(\id \otimes \eta_{\pi_{23},\Delta} \right)
= \left(\zeta^{-1}_{\pi_{23}, \Delta} \otimes \id\right) 
\circ \eta_{\pi_{23}, \Delta} \circ 
\left(\alpha_{\pi_{23} \circ \Delta}\right)
$$ 
of two morphisms 
$$
\left(E^\vee \otimes \pi^!_1\mathcal{O}_{X_1}\right) 
\otimes \Delta^* \pi^*_{12} 
\left( E \otimes \pi^*_1 \left(-\right)\right)
\longrightarrow 
\pi_{23 *} \Delta_* \left(\Delta^* \pi^*_{23} \left(E^\vee \otimes \pi^!_1\mathcal{O}_{X_1}\right) 
\otimes \Delta^* \pi^*_{12} 
\left( E \otimes \pi^*_1 \left(-\right)\right)\right).
$$
Since $\pi_{23} \circ \Delta = \id$, we have 
$\alpha_{\pi_{23} \circ \Delta} = \id$. 
It follows that the right-hand column of  
\eqref{eqn-appendix-transforming-pi2^*-pi_2*-adj-counit-diagram} equals to
\begin{tiny}
\begin{align}
\label{eqn-appendix-transforming-pi2^*-pi_2*-adj-counit-connecting-iso}
\xymatrix{
\pi_{1 *} \left(\pi^!_1\mathcal{O}_{X_1} \otimes \left(E^\vee \otimes
\left( E \otimes \pi^*_1 \left(-\right)\right) \right) \right)  
\ar[d]^{\rho^{-1}}
\\
\pi_{1 *} \left(\left(E^\vee \otimes \pi^!_1\mathcal{O}_{X_1}\right) 
\otimes \left( E \otimes \pi^*_1 \left(-\right)\right)\right)
\ar[d]^{\left(\zeta^{-1}_{\pi_{23}, \Delta} \otimes
\zeta^{-1}_{\pi_{12}, \Delta} \right) \circ \eta_{\pi_{23}, \Delta}}
\\
\pi_{1 *} \pi_{23_*} \Delta_* 
\left(\Delta^* \pi^*_{23} \left(E^\vee \otimes \pi^!_1\mathcal{O}_{X_1}\right) 
\otimes \Delta^* \pi^*_{12} \left( E \otimes \pi^*_1 \left(-\right)\right)\right)
\ar[d]^{
\nu_{\Delta} \circ 
\rho^{-1} \circ 
\left(\id \otimes \nu_{\pi_{12}}\right) \circ 
\nu^{-1}_\Delta 
}
\\
\pi_{1 *} \pi_{23_*} \Delta_* \left(
\Delta^* \left(
\pi^*_{23} \left(E^\vee \otimes \pi^!_1\mathcal{O}_{X_1}\right) 
\otimes \pi^*_{12} E \right) 
\otimes \Delta^* \pi^*_{12} \pi^*_1 \left(-\right)
\right)
} 
\end{align}
\end{tiny}

Note that $\nu^{-1}_{\Delta}$ and $\id \otimes \nu_{\pi_{12}}$ commute
by functoriality. Note further, that by the compatibility of the map 
$\nu_\Delta$ with the associativity of the tensor product (see diagram
\eqref{eqn-nu-is-compatible-with-the-assotiatity-of-the-tensor-product})
we have an equality
$$ 
\nu_{\Delta} \circ \rho^{-1} \circ \nu^{-1}_\Delta
= 
\left( \nu^{-1}_{\Delta} \otimes \id \right)
\circ \rho^{-1} \circ 
\left( \id \otimes \nu_{\Delta}\right)
$$
of two morphisms 
$$
\Delta^* \pi^*_{23} \left(E^\vee \otimes \pi^!_1\mathcal{O}_{X_1}\right) 
\otimes \Delta^* \left(\pi^*_{12} E \otimes \pi^*_{12} \pi^*_1 \left(-\right)\right)
\;\longrightarrow\;
\Delta^* \left(
\pi^*_{23} \left(E^\vee \otimes \pi^!_1\mathcal{O}_{X_1}\right) 
\otimes \pi^*_{12} E \right) 
\otimes \Delta^* \pi^*_{12} \pi^*_1 \left(-\right)
.
$$ 
It follows that composition 
\eqref{eqn-appendix-transforming-pi2^*-pi_2*-adj-counit-connecting-iso} equals
to
\begin{tiny}
\begin{align}
\label{eqn-appendix-transforming-pi2^*-pi_2*-adj-counit-connecting-iso-2}
\xymatrix{
\pi_{1 *} \left(\pi^!_1\mathcal{O}_{X_1} \otimes \left(E^\vee \otimes
\left( E \otimes \pi^*_1 \left(-\right)\right) \right) \right)  
\ar[d]^{\rho^{-1}}
\\
\pi_{1 *} \left(\left(E^\vee \otimes \pi^!_1\mathcal{O}_{X_1}\right) 
\otimes \left( E \otimes \pi^*_1 \left(-\right)\right)\right)
\ar[d]^{\left(\zeta^{-1}_{\pi_{23}, \Delta} \otimes 
\left(
\nu_\Delta \circ \nu_{\pi_{12}} \circ \zeta^{-1}_{\pi_{12}, \Delta} 
\right) 
\right)
\circ \eta_{\pi_{23}, \Delta}}
\\
\pi_{1 *} \pi_{23_*} \Delta_* 
\left(\Delta^* \pi^*_{23} \left(E^\vee \otimes \pi^!_1\mathcal{O}_{X_1}\right) 
\otimes \left( \Delta^* \pi^*_{12} E \otimes \Delta^* \pi^*_{12} \pi^*_1 \left(-\right)\right)\right)
\ar[d]^{\left(\nu^{-1}_\Delta \otimes \id \right) \circ \rho^{-1}}
\\
\pi_{1 *} \pi_{23_*} \Delta_* \left(
\Delta^* \left(
\pi^*_{23} \left(E^\vee \otimes \pi^!_1\mathcal{O}_{X_1}\right) 
\otimes \pi^*_{12} E \right) 
\otimes \Delta^* \pi^*_{12} \pi^*_1 \left(-\right)
\right)
} 
\end{align}
\end{tiny}

By compatibility of $\nu$ with pseudofunctoriality (see diagram
\eqref{eqn-map-nu-is-compatible-with-pseudofunctoriality}) 
we have an equality
$$ 
\nu_\Delta \circ \nu_{\pi_{12}} \circ \zeta^{-1}_{\pi_{12}, \Delta} 
= 
\nu_{\pi_{12} \circ \Delta} \circ \left(\zeta^{-1}_{\pi_{12}, \Delta}
\otimes \zeta^{-1}_{\pi_{12}, \Delta}\right)$$
of morphisms 
$$ E \otimes \pi^*_1(-) \longrightarrow 
\Delta^* \pi^*_{12} E \otimes \Delta^* \pi^*_{12} \pi^*_1(-).$$
Since $\pi_{12} \circ \Delta = \id$ we further have 
$\nu_{\pi_{12} \circ \Delta} = \id$. Therefore 
$$ \left(\zeta^{-1}_{\pi_{23}, \Delta} \otimes 
\left(
\nu_\Delta \circ \nu_{\pi_{12}} \circ \zeta^{-1}_{\pi_{12}, \Delta} 
\right) 
\right)
\circ \eta_{\pi_{23}, \Delta}
= 
\left(\zeta^{-1}_{\pi_{23}, \Delta} \otimes 
\left(
\zeta^{-1}_{\pi_{12}, \Delta} \otimes \zeta^{-1}_{\pi_{12}, \Delta}
\right) 
\right)
\circ \eta_{\pi_{23}, \Delta}
$$
in \eqref{eqn-appendix-transforming-pi2^*-pi_2*-adj-counit-connecting-iso-2}.
Finally, by functoriality of $\rho$ and of 
$\eta_{\pi_{23}, \Delta}$ we have
$$
\rho^{-1} \circ \left(\zeta^{-1}_{\pi_{23}, \Delta} \otimes 
\left(
\zeta^{-1}_{\pi_{12}, \Delta} \otimes \zeta^{-1}_{\pi_{12}, \Delta}
\right) 
\right)
\circ \eta_{\pi_{23}, \Delta}
= 
\left(
\left(
\zeta^{-1}_{\pi_{23}, \Delta} \otimes \zeta^{-1}_{\pi_{12}, \Delta}
\right) \otimes \zeta^{-1}_{\pi_{12}, \Delta} \right)
\circ \eta_{\pi_{23}, \Delta} \circ \rho^{-1}.$$

We conclude that 
\eqref{eqn-appendix-transforming-pi2^*-pi_2*-adj-counit-connecting-iso-2}
equals to 
\begin{tiny}
\begin{align}
\label{eqn-appendix-transforming-pi2^*-pi_2*-adj-counit-connecting-iso-3}
\xymatrix{
\pi_{1 *} \left(\pi^!_1\mathcal{O}_{X_1} \otimes \left(E^\vee \otimes
\left( E \otimes \pi^*_1 \left(-\right)\right) \right) \right)  
\ar[d]^{\rho^{-1} \circ \rho^{-1}}
\\
\pi_{1 *} \left(\left(\left(E^\vee \otimes \pi^!_1\mathcal{O}_{X_1}\right) 
\otimes E\right) \otimes \pi^*_1 \left(-\right)\right)
\ar[d]^{\left(\left( \nu^{-1}_\Delta \circ \left(\zeta^{-1}_{\pi_{23}, \Delta} \otimes 
\zeta^{-1}_{\pi_{12}, \Delta} \right)\right) \otimes \zeta^{-1}_{\pi_{12}, \Delta}
\right)
\circ \eta_{\pi_{23}, \Delta}}
\\
\pi_{1 *} \pi_{23_*} \Delta_* \left(
\Delta^* \left(
\pi^*_{23} \left(E^\vee \otimes \pi^!_1\mathcal{O}_{X_1}\right) 
\otimes \pi^*_{12} E \right) 
\otimes \Delta^* \pi^*_{12} \pi^*_1 \left(-\right)
\right)
} 
\end{align}
\end{tiny}
Recall now that we write 
$Q'$ for $\pi^*_{23} \left(E^\vee \otimes \pi^!_1\mathcal{O}_{X_1}\right) \otimes
\pi^*_{12} E $ and note that 
$\nu_\Delta^{-1} \circ \left(\zeta^{-1}_{\pi_{23}, \Delta}
\otimes \zeta^{-1}_{\pi_{12}, \Delta}\right)$
in \eqref{eqn-appendix-transforming-pi2^*-pi_2*-adj-counit-connecting-iso-3}
is precisely the inverse of isomorphism \eqref{eqn-appendix-Delta-Q'-iso}.  
So what we have shown above is that 
\eqref{eqn-appendix-pi1-pi2-morphism} is isomorphic to 
\begin{align} \label{eqn-appendix-pi_2-pi_13-derived-restriction-prelim}
\pi_{1 *} \pi_{23 *} \left( Q' \otimes \pi^*_{12} \pi^*_1 (-) \right) 
\xrightarrow{\nu_\Delta \circ \beta_\Delta}
\pi_{1 *} \pi_{23 *} \Delta_* \left( \Delta^* Q' \otimes
\Delta^* \pi^*_{12} \pi^*_1 (-) \right)   
\end{align}
with the connecting isomorphism on the RHS being 
$$
\pi_{1 *} \left(\pi^!_1\mathcal{O}_{X_1} \otimes \left(E^\vee \otimes
\left( E \otimes \pi^*_1 \left(-\right)\right) \right) \right)  
\xrightarrow{\left(
\eqref{eqn-appendix-Delta-Q'-iso}^{-1} \otimes
\zeta^{-1}_{\pi_{12},\Delta}\right)
\circ \eta_{\pi_{23},\Delta} \circ \rho^{-1} \circ \rho^{-1}}
\pi_{1 *} \pi_{23 *} \Delta_* 
\left( \Delta^* Q' \otimes \Delta^* \pi^*_{12} \pi^*_1 (-) \right). 
$$
As
$\pi_1 \circ \pi_{23} = \tilde{\pi}_{2} \circ \pi_{13}$
and $\pi_{1} \circ \pi_{12} = \tilde{\pi}_{1} \circ \pi_{13}$ 
(see diagram \eqref{eqn-big-projection-tree}) 
we have the following commutative square
\begin{align*}
\xymatrix{
\pi_{1 *} \pi_{23 *} \left( Q' \otimes \pi^*_{12} \pi^*_1 (-) \right) 
\ar[r]^<<<<<{\eqref{eqn-appendix-pi_2-pi_13-derived-restriction-prelim}}
\ar[d]_{ \left(\id \otimes 
\left(\zeta^{-1}_{\tilde{\pi}_1, \pi_{13}} \circ \zeta_{\pi_1, \pi_{12}}\right)
\right)
\circ \eta_{\tilde{\pi}_{2}, \pi_{13}} \circ
\eta_{\pi_1,\pi_{23}}^{-1}}
& 
\pi_{1 *} \pi_{23 *} \Delta_* \left( \Delta^* Q' \otimes
\Delta^* \pi^*_{12} \pi^*_1 (-) \right)   
\ar[d]^{ \left(\id \otimes 
\left(\zeta^{-1}_{\tilde{\pi}_1, \pi_{13}} \circ \zeta_{\pi_1, \pi_{12}}\right)
\right)
\circ \eta_{\tilde{\pi}_{2}, \pi_{13}} \circ
\eta_{\pi_1,\pi_{23}}^{-1}}
\\
\tilde{\pi}_{2 *} \pi_{13 *} 
\left( Q' \otimes \pi^*_{13} \tilde{\pi}^*_1 (-) \right) 
\ar[r]_<<<<<{\nu_\Delta \circ \beta_\Delta}
& 
\tilde{\pi}_{2 *} \pi_{13 *} \Delta_* 
\left( \Delta^* Q' \otimes \Delta^* \pi^*_{13} \tilde{\pi}^*_1 (-) \right)   
}.
\end{align*}
We finally conclude that \eqref{eqn-appendix-pi1-pi2-morphism}
is isomorphic to
\begin{align} \label{eqn-appendix-pi_2-pi_13-derived-restriction}
\tilde{\pi}_{2 *} \pi_{13 *} 
\left( Q' \otimes \pi^*_{13} \tilde{\pi}^*_1 (-) \right) 
\xrightarrow{\nu_\Delta \circ \beta_\Delta}
\tilde{\pi}_{2 *} \pi_{13 *} \Delta_* 
\left( \Delta^* Q' \otimes \Delta^* \pi^*_{13} \tilde{\pi}^*_1 (-) \right)   
\end{align}
with the connecting isomorphism on the RHS being 
\begin{footnotesize}
\begin{align}
\label{eqn-appendix-pi_2-pi_13-derived-restriction-connecting-iso}
\pi_{1 *} \left(\pi^!_1\mathcal{O}_{X_1} \otimes \left(E^\vee \otimes
\left( E \otimes \pi^*_1 \left(-\right)\right) \right) \right)  
\xrightarrow{\left(
\eqref{eqn-appendix-Delta-Q'-iso}^{-1} \otimes
\zeta^{-1}_{\tilde{\pi}_1, \pi_{13}, \Delta}\right)
\circ \eta_{\tilde{\pi}_{2}, \pi_{13},\Delta} \circ \rho^{-1} \circ \rho^{-1}}
\tilde{\pi}_{2 *} \pi_{13 *} \Delta_* 
\left( \Delta^* Q' \otimes \Delta^* \pi^*_{13} \tilde{\pi}^*_1 (-) \right). 
\end{align}
\end{footnotesize}
Here we have used the fact that by pseudofunctoriality relations 
\eqref{eqn-pseudofunctoriality-relations-identity} and
\eqref{eqn-pseudofunctoriality-relations-associativity}
we have
$$ \eta_{\tilde{\pi}_{2}, \pi_{13}} \circ \eta_{\pi_1,\pi_{23}}^{-1}
\circ \eta_{\pi_{23},\Delta} = \eta_{\tilde{\pi}_{2}, \pi_{13}}
\circ \eta_{\pi_1 \circ \pi_{23}, \Delta} = 
\eta_{\tilde{\pi}_{2}, \pi_{13}}
\circ \eta_{\tilde{\pi}_{2} \circ \pi_{13}, \Delta} 
=
\eta_{\tilde{\pi}_{2}, \pi_{13}, \Delta} $$
and similarly 
$\zeta^{-1}_{\tilde{\pi}_{1}, \pi_{13}} \circ \zeta_{\pi_1,\pi_{12}}
\circ \zeta^{-1}_{\pi_{12},\Delta} = \zeta^{-1}_{\tilde{\pi}_1,
\pi_{13}, \Delta}$. 

Next, we note that the following diagram commutes:
\begin{tiny}
\begin{align}
\label{eqn-appendix-transforming-trace-map-diagram}
\xymatrix{
\pi_{1 *} \left(\pi^!_1\mathcal{O}_{X_1} \otimes \left(E^\vee \otimes
\left( E \otimes \pi^*_1 \left(-\right)\right) \right) \right)  
\ar[rr]^{\eqref{eqn-appendix-pi1-trace-morphism}}
\ar[d]_{\rho^{-1} \circ \rho^{-1}}
& \quad &
\pi_{1 *} \left(\pi^!_1\mathcal{O}_{X_1} \otimes \pi^*_1 \left(-\right)\right)  
\ar[d]^{\id}
\\
\pi_{1 *} \left(\left(\left(\pi^!_1\mathcal{O}_{X_1} \otimes E^\vee\right)  
\otimes E \right)\otimes \pi^*_1 \left(-\right) \right)  
\ar[rr]^{\left(\left(-\right)\otimes E^\vee \right) \otimes E \rightarrow \id}
\ar[d]_{\eqref{eqn-appendix-Delta-Q'-iso}^{-1} \otimes \id}
& \quad &
\pi_{1 *} \left(\pi^!_1\mathcal{O}_{X_1} \otimes \pi^*_1 \left(-\right) \right) 
\ar[d]^{\id}
\\
\pi_{1 *} \left( \Delta^* Q'\otimes \pi^*_1 \left(-\right) \right)  
\ar[rr]^{\eqref{eqn-appendix-Delta-Q'-pi^!_1-O_X1-morphism}}
\ar[d]_{\left( \id \otimes \zeta^{-1}_{\tilde{\pi}_1, \Delta} \right) 
\circ \eta_{\tilde{\pi}_2, \Delta}}
& \quad &
\pi_{1 *} \left( \pi^!_{1}\mathcal{O}_{X_1} \otimes \pi^*_1 \left(-\right) \right)  
\ar[d]^{\left( \id \otimes \zeta^{-1}_{\tilde{\pi}_1, \Delta} \right)
\circ \eta_{\tilde{\pi}_2, \Delta}}
\\
\tilde{\pi}_{2 *} \Delta_*
\pi_{1 *} \left( \Delta^* Q'\otimes  
\pi^*_1 \Delta^* \tilde{\pi}^*_{1}  \left(-\right) \right)  
\ar[rr]^{\eqref{eqn-appendix-Delta-Q'-pi^!_1-O_X1-morphism}}
& \quad & 
\tilde{\pi}_{2 *} \Delta_*
\pi_{1 *} \left( \pi^!_1\mathcal{O}_{X_1} \otimes  
\pi^*_1 \Delta^* \tilde{\pi}^*_{1} \left(-\right) \right)  
} 
\end{align}
\end{tiny}
Here the top square commutes by Lemma
\ref{lemma-tensor-product-adjunction-and-associativity}, the second
square commutes by the definition of map 
\eqref{eqn-appendix-Delta-Q'-pi^!_1-O_X1-morphism} and 
the third square commutes by the functoriality.
Therefore \eqref{eqn-appendix-pi1-trace-morphism} is isomorphic
to 
\begin{align} 
\label{eqn-appendix-pi_2-Delta-pi_1-trace-map}
\tilde{\pi}_{2 *} \Delta_*
\pi_{1 *} \left( \Delta^* Q'\otimes  
\pi^*_1 \Delta^* \tilde{\pi}^*_{1}  \left(-\right) \right)  
\xrightarrow{\eqref{eqn-appendix-Delta-Q'-pi^!_1-O_X1-morphism}}
\tilde{\pi}_{2 *} \Delta_*
\pi_{1 *} \left( \pi^!_1\mathcal{O}_{X_1} \otimes  
\pi^*_1 \Delta^* \tilde{\pi}^*_{1} \left(-\right) \right).
\end{align}

And finally, the following square
\begin{tiny}
\begin{align}
\label{eqn-appendix-transforming-trace-adjunction-diagram}
\xymatrix{
\pi_{1 *} \left( \pi^!_{1}\mathcal{O}_{X_1} \otimes \pi^*_1 \left(-\right) \right)  
\ar[d]_{\left( \id \otimes \zeta^{-1}_{\tilde{\pi}_1, \Delta} \right)
\circ \eta_{\tilde{\pi}_2, \Delta}}
\ar[rr]^{\eqref{eqn-appendix-pi1-trace-adjunction}}
& &
\id
\ar[d]^{\zeta^{-1}_{\tilde{\pi}_1, \Delta} \circ 
\eta_{\tilde{\pi}_2, \Delta}}
\\
\tilde{\pi}_{2 *} \Delta_*
\pi_{1 *} \left( \pi^!_1\mathcal{O}_{X_1}\otimes  
\pi^*_1 \Delta_* \tilde{\pi}^*_{1} \left(-\right) \right)  
\ar[rr]^<<<<<<<<<<<<<<<<{\pi_{1 *}
\left(\pi^!_1\mathcal{O}_{X_1}\otimes\pi^*_1\left(-\right)\right)
\rightarrow \id}
& \quad \quad \quad \quad \quad & 
\tilde{\pi}_{2 *} \Delta_*
\Delta^* \tilde{\pi}^*_{1}  \left(-\right) 
} 
\end{align}
\end{tiny}
commutes by functoriality.  Therefore 
\eqref{eqn-appendix-pi1-trace-adjunction} is isomorphic to 
\begin{align} 
\label{eqn-appendix-pi_2-Delta-pi_1-adjunction}
\tilde{\pi}_{2 *} \Delta_*
\pi_{1 *} \left( \pi^!_1\mathcal{O}_{X_1}\otimes  
\pi^*_1 \Delta_* \tilde{\pi}^*_{1} \left(-\right) \right)  
\xrightarrow{\pi_{1 *}
\left(\pi^!_1\mathcal{O}_{X_1}\otimes\pi^*_1\left(-\right)\right)
\rightarrow \id}
\tilde{\pi}_{2 *} \Delta_*
\Delta^* \tilde{\pi}^*_{1}  \left(-\right). 
\end{align}

We now compute the connecting isomorphisms. Composing 
the inverse of
\eqref{eqn-appendix-pi_2-pi_13-derived-restriction-connecting-iso},
the isomorphism in the right column of 
\eqref{eqn-appendix-transforming-pi2^*-pi_2*-adj-counit-diagram}, 
with the isomorphism in the left column of 
\eqref{eqn-appendix-transforming-trace-map-diagram} we obtain
$$
\tilde{\pi}_{2 *} \pi_{13 *} \Delta_* 
\left( \Delta^* Q' \otimes \Delta^* \pi^*_{13} \tilde{\pi}^*_1 (-) \right)   
\xrightarrow{\; \left( \id \otimes 
\left(\zeta^{-1}_{\tilde{\pi}_1, \Delta} \circ
\zeta_{\tilde{\pi}_{1}, \pi_{13},\Delta} \right)\right)
\circ \eta_{\tilde{\pi}_2, \Delta} \eta^{-1}_{\tilde{\pi}_{2},
\pi_{13},\Delta} \;}
\tilde{\pi}_{2 *} \Delta_*
\pi_{1 *} \left( \Delta^* Q'\otimes  
\pi^*_1 \Delta^* \tilde{\pi}^*_{1}  \left(-\right) \right)
$$
and 
by pseudofunctoriality relations 
\eqref{eqn-pseudofunctoriality-relations-identity} and
\eqref{eqn-pseudofunctoriality-relations-associativity} this
is equal to
\begin{align}
\label{eqn-appendix-pi_2-pi_13-Delta-to-Delta-pi_1-relabelling}
\tilde{\pi}_{2 *} \pi_{13 *} \Delta_* 
\left( \Delta^* Q' \otimes \Delta^* \pi^*_{13} \tilde{\pi}^*_1 (-) \right)   
\xrightarrow{\left( \id \otimes 
\left(\zeta^{-1}_{\Delta, \pi_1} \circ
\zeta_{\pi_{13},\Delta} \right)\right)
\circ \eta_{\Delta, \pi_1} \circ \eta^{-1}_{\pi_{13},\Delta}}
\tilde{\pi}_{2 *} \Delta_*
\pi_{1 *} \left( \Delta^* Q'\otimes  
\pi^*_1 \Delta^* \tilde{\pi}^*_{1}  \left(-\right) \right).
\end{align}
On the other hand, the composition of the inverse of the isomorphism
in the right column of 
\eqref{eqn-appendix-transforming-trace-map-diagram} with the isomorphism 
in the left column of  
\eqref{eqn-appendix-transforming-trace-adjunction-diagram} is clearly $\id$. 

We can now conclude that the adjunction counit 
$\Phi^{\text{ladj}}_E \Phi_E \rightarrow \id$, 
being the composition of \eqref{eqn-appendix-pi1-pi2-morphism},
\eqref{eqn-appendix-pi1-trace-morphism} and
\eqref{eqn-appendix-pi1-trace-adjunction},
is isomorphic 
to the composition of \eqref{eqn-appendix-pi_2-pi_13-derived-restriction}, 
\eqref{eqn-appendix-pi_2-pi_13-Delta-to-Delta-pi_1-relabelling},
\eqref{eqn-appendix-pi_2-Delta-pi_1-trace-map} and
\eqref{eqn-appendix-pi_2-Delta-pi_1-adjunction}.
The claim of the theorem then follows from the fact that the 
following diagram commutes:
\begin{align} \label{eqn-appendix-big-commutative-ladder}
\xymatrix{
\tilde{\pi}_{2 *} \left(\pi_{13 *} Q' \otimes \tilde{\pi}^*_1(-) \right)
\ar[r]^{\sim} \ar[d]_{\eqref{eqn-appendix-transform-derived-restriction}} &
\tilde{\pi}_{2 *} \pi_{13 *} \left(Q' \otimes \pi^*_{13} \tilde{\pi}^*_1(-) \right)
\ar[d]^{\eqref{eqn-appendix-pi_2-pi_13-derived-restriction}} \\
\tilde{\pi}_{2 *} \left(\pi_{13 *} \Delta_* \Delta^* Q' \otimes \tilde{\pi}^*_1(-) \right)
\ar[r]^{\sim}
\ar[d]_{\eqref{eqn-appendix-transform-pi_13-Delta-to-Delta-pi_1-relabelling}} &
\tilde{\pi}_{2 *} \pi_{13 *} \Delta_* \left(\Delta^* Q' \otimes \Delta^* \pi^*_{13} \tilde{\pi}^*_1(-) \right)
\ar[d]^{\eqref{eqn-appendix-pi_2-pi_13-Delta-to-Delta-pi_1-relabelling}} \\
\tilde{\pi}_{2 *} \left(\Delta_* \pi_{1 *} \Delta^* Q' \otimes \tilde{\pi}^*_1(-) \right)
\ar[r]^{\sim} \ar[d]_{\eqref{eqn-appendix-transform-trace-map}} & 
\tilde{\pi}_{2 *} \Delta_* \pi_{1 *} \left(\Delta^* Q' \otimes \pi^*_{1} \Delta^* \tilde{\pi}^*_1(-) \right)
\ar[d]^{\eqref{eqn-appendix-pi_2-Delta-pi_1-trace-map}} \\
\tilde{\pi}_{2 *} \left(\Delta_* \pi_{1 *} \pi^{!}_1\mathcal{O}_{X_1} \otimes \tilde{\pi}^*_1(-) \right)
\ar[r]^{\sim} \ar[d]_{\eqref{eqn-appendix-transform-adjunction}} & 
\tilde{\pi}_{2 *} \Delta_* \pi_{1 *} \left(\pi^{!}_1\mathcal{O}_{X_1} \otimes \pi^*_{1} \Delta^* \tilde{\pi}^*_1(-) \right)
\ar[d]^{\eqref{eqn-appendix-pi_2-Delta-pi_1-adjunction}} \\
\tilde{\pi}_{2 *} \left(\Delta_* \mathcal{O}_{X_1} \otimes \tilde{\pi}^*_1(-) \right)
\ar[r]^{\sim} &
\tilde{\pi}_{2 *} \Delta_* \Delta^* \tilde{\pi}^*_1(-)
} 
\end{align}
where the horizontal isomorphisms are all due to the projection formula.
To see that diagram \eqref{eqn-appendix-big-commutative-ladder} indeed 
commutes, 
observe that its topmost square commutes by Lemma
\ref{lemma-projection-formula-commutes-with-adjunction-1}, the middle 
two commute by functoriality and the lowermost square commutes by Lemma 
\ref{lemma-projection-formula-commutes-with-adjunction-2}.
\end{proof}

\bibliography{../../references}
\bibliographystyle{amsalpha}

\end{document}